\documentclass[12pt]{article}
\usepackage{amssymb,amsmath,bm,mathrsfs,makeidx,amsfonts,graphicx,amsthm,color}
\usepackage[square, comma, sort&compress, numbers]{natbib}
\usepackage{epsfig}
\usepackage{rotating}
\usepackage{mathtools}
\usepackage{bm}
\usepackage{extarrows}
\usepackage[colorlinks,linkcolor=red,anchorcolor=blue,citecolor=blue]{hyperref}
\usepackage{booktabs}
\usepackage{xr}
\usepackage{chngpage}

\usepackage[title]{appendix}

\numberwithin{equation}{section}
\def\beginn{\begin{eqnarray*}}
\def\endn{\end{eqnarray*}}
\def\beginy{\begin{eqnarray}}
\def\endy{\end{eqnarray}}
\def\begine{\begin{enumerate}}
\def\ende{\end{enumerate}}

\def\be{\begin{equation}}
\def\ee{\end{equation}}
\def\bea{\begin{eqnarray}}
\def\eea{\end{eqnarray}}

\numberwithin{equation}{section}
\theoremstyle{plain}
\newtheorem{thm}{Theorem}[section]
\newtheorem{prop}{Proposition}[section]
\newtheorem{lem}{Lemma}[section]
\newtheorem{coro}{Corollary}[section]
\newtheorem{rmk}{Remark}[section]

\newtheorem{deff}{Definition}[section]

\newtheorem{cond}{Condition}[section]
\newcommand{\non}{\nonumber \\}
\newcommand{\bbA}{{\bf A}}

\newcommand{\bbB}{{\bf B}}
\newcommand{\bbC}{{\bf C}}

\newcommand{\bbD}{{\bf D}}

\newcommand{\bbe}{{\bf e}}

\newcommand{\bbI}{{\bf I}}

\newcommand{\bbr}{{\bf r}}

\newcommand{\bbS}{{\bf S}}

\newcommand{\bbV}{{\bf V}}
\newcommand{\bbv}{{\bf v}}

\newcommand{\bbX}{{\bf X}}

\newcommand{\bbx}{{\bf x}}

\newcommand{\bbY}{{\bf Y}}
\newcommand{\bby}{{\bf y}}
\newcommand{\bbZ}{{\bf Z}}
\newcommand{\bbz}{{\bf z}}

\newcommand{\bbzero}{{\bf 0}}
\newcommand{\bbone}{{\bf 1}}
\newcommand{\bbmu}{{\boldsymbol\mu}}

\newcommand{\tr}{\mathrm{tr}}
\newcommand{\E}{\mathbb{E}}
\newcommand{\D}{\mathcal{D}}
\newcommand{\s}{\mathcal{S}}
\newcommand{\C}{\mathcal{C}}

\newcommand{\pr}{\mathbb{P}}

\newcommand{\diff}{\mathrm{Diff}}

\begin{document}

\title{
Quadratic Discriminant Analysis under Moderate Dimension}

 \author{Qing Yang\footnotemark[1]\;\;  and  Guang Cheng\footnotemark[2]}
\renewcommand{\thefootnote}{\fnsymbol{footnote}}
\footnotetext[1]{Postdoctoral Researcher, Department of Statistics, Purdue University, West Lafayette, IN 47906. (Email: qyang1@e.ntu.edu.sg).}
\footnotetext[2]{Corresponding Author. Professor, Department of Statistics, Purdue University, West Lafayette, IN 47906. (E-mail: chengg@purdue.edu). Research Sponsored by NSF DMS-1712907, DMS-1811812, DMS-1821183, and Office of Naval Research, (ONR N00014-18-2759).}

\date{}
\maketitle

\begin{abstract}
Quadratic discriminant analysis (QDA) is a simple method to classify a subject into two populations, and was proven to perform as well as the Bayes rule when the data dimension $p$ is fixed. The main purpose of this paper is to examine the empirical and theoretical behaviors of QDA where $p$ grows proportionally to the sample sizes without imposing any structural assumption on the parameters. The first finding in this moderate dimension regime is that QDA can perform as poorly as random guessing even when the two populations deviate significantly. This motivates a generalized version of QDA that automatically adapts to dimensionality. Under a finite fourth moment condition, we derive misclassification rates for both the generalized QDA and the optimal one. A direct comparison reveals one ``easy'' case where the difference between two rates converges to zero and one ``hard'' case where that converges to some strictly positive constant. For the latter, a divide-and-conquer approach over dimension (rather than sample) followed by a screening procedure is proposed to narrow the gap. Various numerical studies are conducted to back up the proposed methodology.
\end{abstract}

\vspace{10pt}


\newpage

\section{Introduction}\label{sec:intr}
Suppose we have two $p$-variate classes with mean vectors and covariance matrices $(\bbmu_1,\Sigma_1)$ (class 1) and $(\bbmu_2,\Sigma_2)$ (class 2) respectively. The aim is to identify to which class a new observation $\bbz$ belongs, on the basis of two sets of training samples $\{\bbx_1,\cdots, \bbx_{n_1}\}$ and $\{\bby_1,\cdots, \bby_{n_2}\}$. This problem has been well studied in the fixed dimensional setting, see \cite{ander84} for example. Some recent high dimensional studies allow the data dimension $p$ to be much larger than the sample sizes, but heavily rely on the sparsity or other structural assumptions on the population parameters $\bbmu_i$ and $\Sigma_i$ $(i=1,2)$. For example, $(\bbmu_2-\bbmu_1)$ is sparse or $\Sigma_i$s satisfy special structures such that better estimators (e.g. thresholding, diagonalization) can be constructed --
many improvements have been made over the classical classification rules. One may refer to \cite{cai11, Fan12, Fan13, Hao15, shao11, Tib11} among others for such improvements over the well-known Fisher's linear discriminant analysis (LDA) when $\Sigma_1=\Sigma_2$; and see \cite{Cheng04, Fan15, shao15} and  \cite{Qin18}, to list but a few for the modifications over the quadratic discriminant analysis (QDA) when $\Sigma_1\neq\Sigma_2$.

Despite these recent progress, there has been relatively fewer development for classification in the moderate dimension regime, by which we mean that $p/n_i\rightarrow c_i\in (0, 1)$ and there is no structural assumption on the parameters $\bbmu_i$ and $\Sigma_i$. One exception is moderate-dimensional LDA by \cite{wang17}. Clearly, their analyses do not apply to QDA that mainly relies on a quadratic form. Another somehow related line is the study of logistic regression that can be used as a classification rule; see \cite{candes17, candes18}. However, these works mostly focused on the inference results on the parameter itself, e.g., log-likelihood ratio test, rather than classification performances.  One can also refer to \cite{Dob18, Dono16, Kar13, Kar17, Karo13, Guo18, Jan17} and \cite{Lei18} for other moderate dimensional results in the linear regression models.

It is well known that the state of the art quadratic discriminant rule classifies $\bbz$ to class 1 if and only if
\begin{equation}\label{yq1}
d_1(\bbz)+\log |\Sigma_1|<d_2(\bbz)+\log |\Sigma_2|,
\end{equation}
where ``$|A|$'' denotes the determinant of the matrix $A$ and
\begin{equation}\label{yq2}
d_i(\bbz)=(\bbz-\bbmu_i)^T\Sigma_i^{-1}(\bbz-\bbmu_i),\quad i=1,2.
\end{equation}
We call it ``optimal QDA,'' which serves as a performance benchmark in our paper. In practice, the population parameters need to be estimated from training samples, leading to what we call as ``sample QDA'' as follows. That is to say, the new observation $\bbz$ is classified to class 1 if and only if \begin{equation}\label{yq4}
D_1(\bbz)+\log |S_1|<D_2(\bbz)+\log |S_2|,
\end{equation}
where
\begin{eqnarray}\label{yq5}
&&D_1(\bbz)=(\bbz-\bar\bbx)^TS_1^{-1}(\bbz-\bar\bbx),\quad
D_2(\bbz)=(\bbz-\bar\bby)^TS_2^{-1}(\bbz-\bar\bby),\nonumber\\
&&S_1=\frac{1}{n_1-1}\sum\limits_{i=1}^{n_1}(\bbx_i-\bar\bbx)(\bbx_i-\bar\bbx)^T,\quad
S_2=\frac{1}{n_2-1}\sum\limits_{i=1}^{n_2}(\bby_i-\bar\bby)(\bby_i-\bar\bby)^T.
\end{eqnarray}

\begin{figure}[!htp]
    \centering
    \includegraphics[width=7in,height=1.8in]{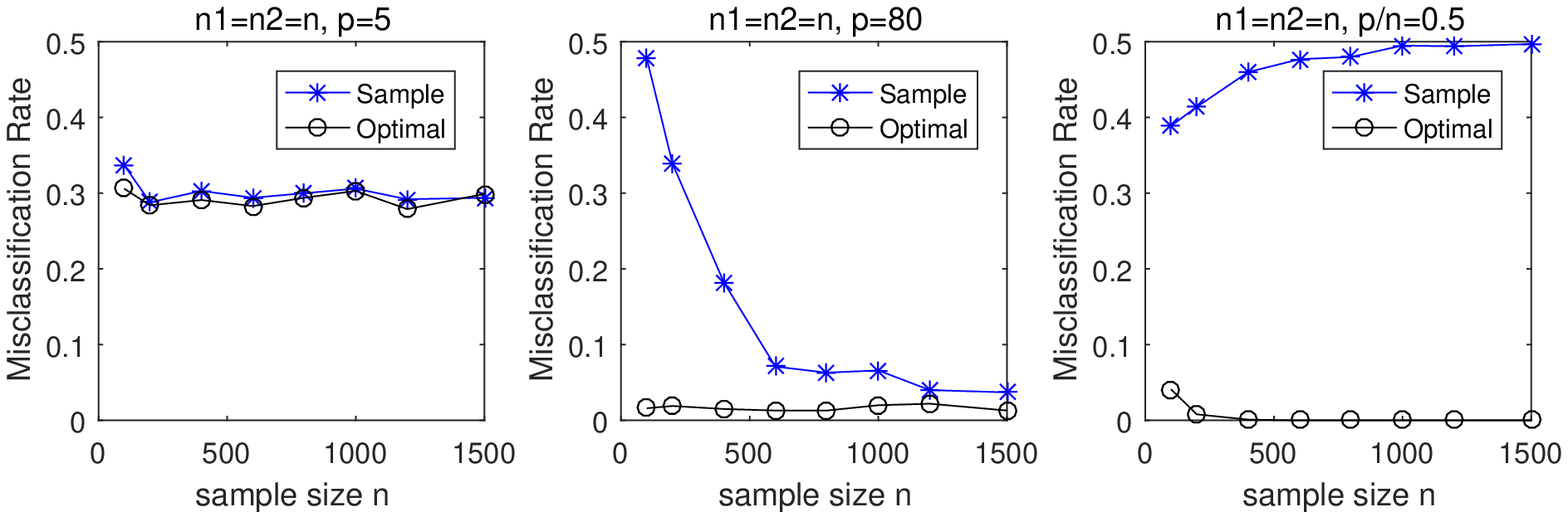}
    \caption{\small{\it Misclassification rates for two classes $N_p(\bbzero,\bbI_p)$ and $N_p(\bbzero,2\bbI_p)$ based on 1000 replications.}}
    \label{eg1}
  \end{figure}

We next conduct a simple experiment to empirically examine the classification performances of sample QDA when $p$ is moderately large compared with the sample size. Specifically, we generate $n_1$ and $n_2$ samples from two classes $N_p(\bbzero, \bbI_p)$ and $N_p(\bbzero,2\bbI_p)$, respectively. Set $n_1=n_2=n$ with $n$ varying from $100$ to $1500$. The misclassification rate is defined as
\begin{equation}\label{1804.2}
R=\frac{1}{2}\left[\pr_{2|1}+\pr_{1|2}\right],
\qquad
\pr_{i|j}\triangleq\pr\{\text{classify}\ \bbz\ \text{to}\ \text{class }\ i\ |\bbz\in \text{class }\ j\}.
\end{equation}
Denote $R^S$ and $R^O$ as misclassification rates for sample QDA and optimal QDA\footnote{Throughout this paper, superscripts ``O'', ``S'' and ``G'' are adopted to denote the values corresponding to the optimal QDA (\ref{yq1}), the sample QDA (\ref{yq4}) and our generalized QDA (\ref{11.1}) introduced later.}, respectively. From Figure \ref{eg1} (based on 1000 replications), we observe that the sample QDA matches well with the optimal one when $p$ is very small compared with $n$ (the left plot), but the gap between them becomes notable when $n$ is not significantly large compared with $p$ (the middle plot) and $R^S$ even converges to 0.5 (random guessing) when $p$ is proportional to $n$ (the right plot). This is an indication of sample QDA's failure in the moderate-dimensional setting.

The rationale behind the above moderate-dimension phenomenon is simple: the sample mean vectors and sample covariance matrices are no longer consistent in terms of the $L_2$ norm and the spectral norm, respectively. In fact, a careful analysis of two terms $D_i(\bbz)$ and $\log|S_i|$ reveals that an extra scaling factor to the former and an extra shift factor to latter need to be introduced to adapt to dimensional effect; see Section~\ref{yq18.1}. This leads to a generalized QDA that is built on these two corrected terms; see Section~\ref{yq18.3}.

We employ various techniques in the random matrix theory to derive misclassification rates of the generalized QDA, the optimal QDA and the sample QDA in Sections \ref{rescal} and \ref{optim}, without imposing any structural or parametric distribution assumption. This is radically different from the QDA results in the literature that were developed based on normal distribution. In comparison, we only need a finite fourth moment condition on the data, which we believe to be the weakest moment assumption. This can be achieved mainly because we figure out the limiting point of an average squared diagonal entries of a general inverse sample covariance matrix, as formulated in Proposition \ref{prop3}. Such a result is new to our knowledge. By comparing the misclassification rates between the generalized QDA and the optimal one, we find two cases -- the ``easy'' case where the rate difference converges to zero and the ``hard'' case where the rate difference converges to some strictly positive constant.

Various simulation studies in Section \ref{simu} support our theoretical conclusions. To narrow the gap in the ``hard'' case, two divide-and-conquer approaches are proposed in Section \ref{dc}. In contrast with the conventional partition over samples (which is shown not to work in Supplement \ref{negativ}), our methods are conducted over the dimension followed by a screening procedure. Both show significant improvement over the original generalized QDA. Unless otherwise noted, all the main proofs are relegated to  Appendix \ref{profs}. Extra theoretical and numerical results are included in the supplementary material.

\section{Methodology -- Tame the Dimensionality Effect}\label{sec:2}

In this section, we study the performances of $D_i(\bbz)$ and $\log|S_i|$ in (\ref{yq4}) under moderate dimension, see Section \ref{yq18.1}, which motivate a dimension adaptive version of QDA in Section \ref{yq18.3}. A set of mild conditions are needed throughout this paper.
\begin{cond}\label{cond1} [Population]
The class 1 has the form $\bbx_i=\Sigma_1^{\frac{1}{2}}\bbx_i^0+\bbmu_1$, $i=1,\cdots, n_1$, where $\bbx_i^0=(X_{i1},\cdots,X_{ip})^T$ has $p$ i.i.d centered and standardized components with finite fourth moments, i.e. $m_4:=\E|X_{ij}|^4< \infty$. The same form applies to class 2, i.e. $\bby_i=\Sigma_2^{\frac{1}{2}}\bby_i^0+\bbmu_2$. Moreover, there exist $c$ and $C$ s.t. $0<c\leq \lambda_{\min}(\Sigma_i)\leq \lambda_{\max}(\Sigma_i)\leq C<\infty$, $i=1,2$.
\end{cond}
\begin{cond}\label{cond2}[Dimensionality]
$p/n_1\rightarrow c_1\in (0,1)$ and $p/n_2\rightarrow c_2\in (0,1)$.
\end{cond}
\begin{cond}\label{cond3}[Covariance Matrix]
The following limits exist
\[
\frac{1}{p}\tr(\Sigma_1\Sigma_2^{-1})\rightarrow M_1, \qquad
\frac{1}{p}\sum_{i=1}^p[(\Sigma_1^{\frac{1}{2}}\Sigma_2^{-1}\Sigma_1^{\frac{1}{2}})_{ii}]^2\rightarrow M_2,
\]
\[
\frac{1}{p}\tr(\Sigma_2\Sigma_1^{-1})\rightarrow M_3, \qquad
\frac{1}{p}\sum_{i=1}^p[(\Sigma_2^{\frac{1}{2}}\Sigma_1^{-1}\Sigma_2^{\frac{1}{2}})_{ii}]^2\rightarrow M_4,
\]
\[
\frac{1}{p}\tr(\Sigma_1\Sigma_2^{-1})^2\rightarrow M_5,\qquad
\frac{1}{p}\tr(\Sigma_2\Sigma_1^{-1})^2\rightarrow M_6.
\]
\end{cond}

\begin{rmk}  A simple example to understand the limits in Condition \ref{cond3}:
if $\Sigma_1=\kappa\Sigma_2$, $\kappa>0$, then $M_1=\kappa$, $M_3=\frac{1}{\kappa}$, $M_2=M_5=\kappa^2$ and $M_4=M_6=\frac{1}{\kappa^2}$.
\end{rmk}

\subsection{Dimension distortion of $D_i(\bbz)$ and $\log|S_i|$}\label{yq18.1}

One observes that the mean values of $d_1(\bbz)$ and $d_1(\bbz)$ in (\ref{yq2}) are equal when $\bbz$ is correctly classified irrespective of data dimension, i.e., if $\bbz$ belongs to class 1, $\E d_1(\bbz)=p$ and if $\bbz$ belongs to class 2, $\E d_2(\bbz)=p$. However, such an equivalence does not hold for their sample version $D_i(\bbz)$ in (\ref{yq5}) when $p$ diverges proportionally to $n$.

To be more precise, we derive the following limiting distributions for $D_i(\bbz)$, $i=1,2$. This weak convergence result is also needed in deriving the misclassification rate in Section \ref{sec:thy}.

\begin{thm}\label{them1}
Under Conditions \ref{cond1} and  \ref{cond2}, if $\bbz$ belongs to class 1, we have
\begin{equation}\label{1031.1}
\frac{1}{\sqrt{p}}D_1(\bbz)-\sqrt{p}s_{0n}\xrightarrow{D} N\Big(0,(m_4-3)s_0^2+2s'_0\Big).
\end{equation}
If $\bbz$ belongs to class 2, we have
\begin{equation}\label{1031.2}
\frac{1}{\sqrt{p}}D_2(\bbz)-\sqrt{p}m_{0n}\xrightarrow{D} N\Big(0,(m_4-3)m_0^2+2m'_0\Big).
\end{equation}
Here ``$\xrightarrow{D}$'' denotes convergence in distribution, and
{\small
\begin{equation}\label{1221.1}
s_{0n}=\frac{1}{1-p/n_1},\ m_{0n}=\frac{1}{1-p/n_2},\ s_0=\frac{1}{1-c_1},\ s'_0=\frac{1}{(1-c_1)^3},\ m_0=\frac{1}{1-c_2},\ m'_0=\frac{1}{(1-c_2)^3},
\end{equation}}
where $c_1$ and $c_2$ are given in Condition \ref{cond2}.
\end{thm}

Theorem \ref{them1} implies that under moderate dimension, $D_1(\bbz)\approx ps_{0n}$ when $\bbz$ belongs to class 1 and $D_2(\bbz)\approx pm_{0n}$ when $\bbz$ belongs to class 2 (covering the fixed dimensional case, i.e., $s_{0n}=m_{0n}=1$, as a special case). Hence, to counteract the effects of moderate dimension, we need to rescale the quadratic terms in the sample QDA (\ref{yq4}) as follows
\begin{equation}\label{1221.2}
\frac{1}{s_{0n}}D_1(\bbz)\qquad \text{and}\qquad \frac{1}{m_{0n}}D_2(\bbz).
\end{equation}

A ``re-centering'' type of dimension correction will be applied to $\log|S_i|$ in (\ref{yq4}). Define
\begin{equation}\label{1805.2}
\bbX^0=(\bbx_1^0,\cdots,\bbx_{n_1}^0),\quad \bar\bbx^0=\frac{1}{n_1}\sum_{i=1}^{n_1}\bbx_{i}^0,\quad \bar\bbX^0=\bar\bbx^0\cdot \bbone_{n_1}^T,\quad
S_1^0=\frac{1}{n_1-1}(\bbX^0-\bar\bbX^0)(\bbX^0-\bar\bbX^0)^T.
\end{equation}
Here $\bbone_{n_1}$ is a $n_1$-variate vector with all entries being ones. Notations for class 2 can be defined similarly. For example,  $S_2^0=\frac{1}{n_2-1}(\bbY^0-\bar\bbY^0)(\bbY^0-\bar\bbY^0)^T$. Moreover, from \cite{Bai04} and \cite{pan14}, we know that with probability 1, \[\frac{1}{p}\log|S_i^0|- l_{in}\rightarrow 0,\quad l_{in}=\frac{p/n_i-1}{p/n_i}\log(1-p/n_i)-1<0.\] Let $S_1=\Sigma_1^{\frac{1}{2}}S_1^0\Sigma_1^{\frac{1}{2}}$ and
$S_2=\Sigma_2^{\frac{1}{2}}S_2^0\Sigma_2^{\frac{1}{2}}$.
The relation between $S_i$ and $S_i^0$ thus implies that
\begin{equation}\label{1221.3}
\frac{1}{p}\log|S_i|-\left[\frac{1}{p}\log|\Sigma_i|+l_{in}\right]\rightarrow 0,\quad i=1,2.
\end{equation}
Note that (\ref{1221.3}) recovers the classical setting that $p/n_i\rightarrow 0$, i.e., $l_{1n}=l_{2n}=0$. However, when $p/n_i\nrightarrow 0$, we have to re-center $\log|S_i|$ in (\ref{yq4}) as follows:
\begin{equation}\label{1804.1}
\log|S_1|-p l_{1n}\qquad \text{and}\qquad \log|S_2|-p l_{2n}.
\end{equation}
In light of (\ref{1221.2}) with (\ref{1804.1}), we propose a generalized QDA that is adaptive to data dimension.

\subsection{A generalized version of QDA}\label{yq18.3}

We define the generalized QDA rule: classify a new observation $\bbz$ to class 1 if and only if
\begin{equation}\label{11.1}
\frac{1}{s_{0n}}D_1(\bbz)+\log |S_1|-p l_{1n}<\frac{1}{m_{0n}}D_2(\bbz)+\log |S_2|-p l_{2n},
\end{equation}
where recall that
{\small
\begin{equation}\label{1805.1}
s_{0n}=\frac{1}{1-p/n_1},\ m_{0n}=\frac{1}{1-p/n_2},\ l_{1n}=\frac{p/n_1-1}{p/n_1}\log(1-p/n_1)-1,\  l_{2n}=\frac{p/n_2-1}{p/n_2}\log(1-p/n_2)-1.
\end{equation}}

\begin{figure}[!htp]
    \centering
    \includegraphics[width=7in]{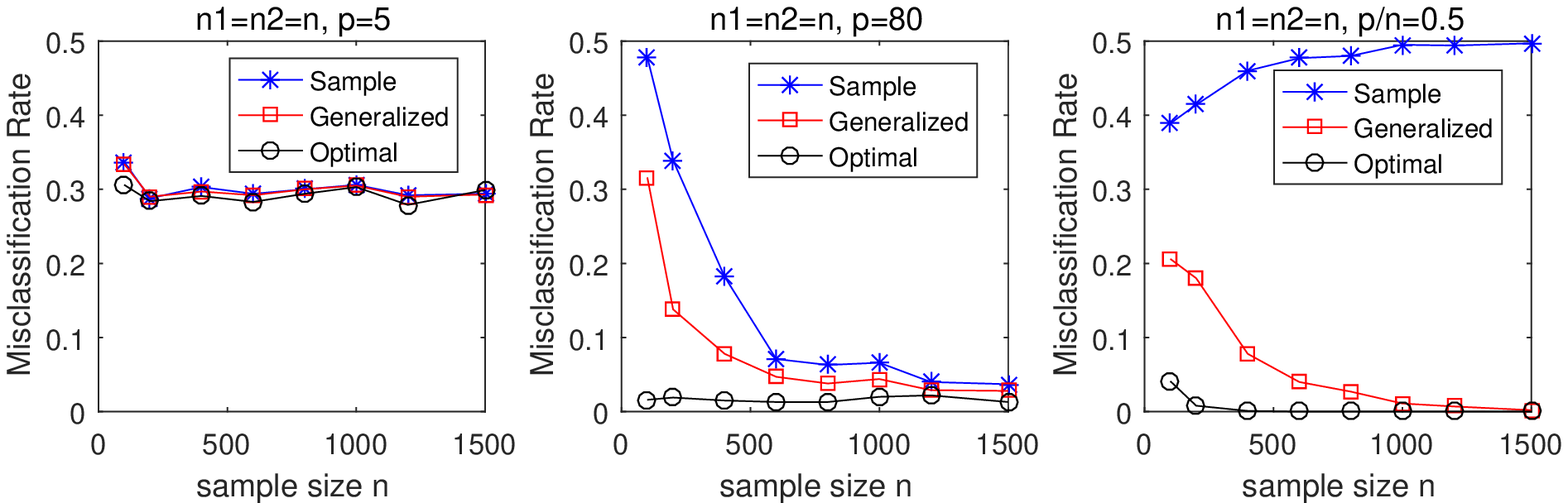}
    \caption{\small{\it Misclassification rates for two classes $N_p(\bbzero,\bbI_p)$ and $N_p(\bbzero,2\bbI_p)$ based on 1000 replications.}}
    \label{eg2}
  \end{figure}

Now we test the classification performances of (\ref{11.1}) by re-visiting the example considered in Section~\ref{sec:intr} (more comprehensive numerical analysis will be conducted in Section~\ref{simu}). Specifically, Figure \ref{eg2} demonstrates that the generalized QDA maintains comparable misclassification rates to the sample QDA (without dimension correction) in the low-dimensional regime. While it significantly diminishes the misclassification rate comparing to the sample QDA, and eventually converges to the optimal one in the moderate-dimensional case.

\subsection{Another interpretation of generalized QDA}\label{interp}

In this section, we offer a more intuitive way to justify generalized QDA, which leads to more rigorous analysis on the misclassification rate in the next section.

Our analysis starts from the optimal rule (\ref{yq1}). In order to correctly classify the new observation $\bbz$ to class 1, we need
\begin{equation*}
\mathcal{T}_{2|1}^O\triangleq d_1(\bbz)+\log |\Sigma_1|-d_2(\bbz)-\log |\Sigma_2|<0.
\end{equation*}
From the proof of Proposition \ref{prop1} below, we know that $$\frac{1}{\sqrt{p}}\left[\mathcal{T}_{2|1}^O-\mathcal{E}_{2|1}^O\right]\xrightarrow{D}N(0, \sigma^2),$$ where $\mathcal{E}_{2|1}^O$ is the expectation of $\mathcal{T}_{2|1}^O$ when $\bbz$ belongs to class 1, and
\begin{equation*}
\mathcal{E}_{2|1}^O=\tr(\bbI_p-\Sigma_1\Sigma_2^{-1})+\log|\Sigma_1\Sigma_2^{-1}|
-(\bbmu_1-\bbmu_2)^T\Sigma_2^{-1}(\bbmu_1-\bbmu_2).
\end{equation*}
Note that one can write
$
\tr(\bbI_p-\Sigma_1\Sigma_2^{-1})+\log|\Sigma_1\Sigma_2^{-1}|=\sum_{i=1}^p\big(1-\lambda_i+\log\lambda_i\big),
$
where $\lambda_i>0$ are the eigenvalues of the matrix $(\Sigma_1\Sigma_2^{-1})$. Define the function
$f(x)=1-x+\log x$, $x>0$.
It is easy to see that $f(x)\leq 0$ and $f(x)=0$ if and only if $x=1$. Therefore, as long as $\Sigma_1\neq \Sigma_2$,
$
\tr(\bbI_p-\Sigma_1\Sigma_2^{-1})+\log|\Sigma_1\Sigma_2^{-1}|<0.
$
Together with the fact that
$-(\bbmu_1-\bbmu_2)^T\Sigma_2^{-1}(\bbmu_1-\bbmu_2)\leq 0,$
we can see
$\mathcal{E}_{2|1}^O<0,$
which prompts the correct identification $\mathcal{T}_{2|1}^O<0$.

We next apply similar analysis to the sample QDA (\ref{yq4}). From the proofs of Lemmas \ref{lema2} and \ref{lema3} below, the counterparts to $\mathcal{T}_{2|1}^O$ and $\mathcal{E}_{2|1}^O$ are
$$\mathcal{T}_{2|1}^S= D_1(\bbz)+\log |S_1|-D_2(\bbz)-\log |S_2|\quad \text{and}\non$$
$$\mathcal{E}_{2|1}^S=s_{0n}\tr\bbI_p-m_{0n}\tr(\Sigma_1\Sigma_2^{-1})+
\log|\Sigma_1\Sigma_2^{-1}|
+pl_{1n}-pl_{2n}-m_{0n}(\bbmu_1-\bbmu_2)^T\Sigma_2^{-1}(\bbmu_1-\bbmu_2),
$$
which reduces to $\mathcal{E}_{2|1}^O<0$ when $p$ is fixed, i.e., $s_{0n}=m_{0n}=1$ and $l_{1n}=l_{2n}=0$. However, when $p$ diverges, the sign of $\mathcal{E}_{2|1}^S$ is no longer determined. For example, let $\bbmu_1=\bbmu_2$, $\Sigma_1=\kappa\Sigma_2$, $0<\kappa<1$ and $n_1=n_2$, then $\mathcal{E}_{2|1}^S=s_{0n}p-s_{0n}p\kappa+p\log\kappa$.
\begin{itemize}
\item[(i)]If $s_{0n}>\frac{-\log k}{1-k}$, then $\mathcal{E}_{2|1}^S>0$.
\item[(ii)]If $s_{0n}=\frac{-\log k}{1-k}$, then $\mathcal{E}_{2|1}^S=0$.
\item[(iii)]\footnote{Note that $-\log \kappa>1-\kappa$ when $0<\kappa<1$.}If $1<s_{0n}<\frac{-\log k}{1-k}$, then $\mathcal{E}_{2|1}^S<0$.
\end{itemize}
However, by adopting the generalized QDA in (\ref{11.1}), we have (based on proofs of Lemmas \ref{lema2} and \ref{lema3}, again)
$$\mathcal{E}_{2|1}^G=\mathcal{E}_{2|1}^O<0,$$ which justifies the similar misclassification rates as the optimal QDA.

\section{Misclassification Rate Results}\label{sec:thy}

In this section, we derive the asymptotic misclassification rates of the generalized QDA and the optimal QDA. By comparing rate difference, we specify two cases -- the ``easy'' case when the difference converges to zero and the ``hard'' case when a non-degenerate difference exists.

\subsection{Misclassification rate of the generalized QDA}\label{rescal}

The misclassification rate of the generalized QDA is written as
$$R^G=\frac{1}{2}\left[\pr_{2|1}^G+\pr_{1|2}^G\right],$$ where Lemmas \ref{lema2} and \ref{lema3} give the limits of $\pr_{2|1}^G$ and $\pr_{1|2}^G$, respectively.

\begin{lem}\label{lema2}
When $\bbz$ belongs to class 1, under Conditions \ref{cond1}-\ref{cond3}, we have
\begin{eqnarray*}
\pr_{2|1}^G\xrightarrow{i.p} 1-\Phi\left(\frac{T}{\psi}\right),
\end{eqnarray*}
where
\begin{equation}\label{1220.1}
T=t_2+t_3,\quad t_2=-\lim_{p\rightarrow \infty}T_2,\quad t_3=-\lim_{p\rightarrow \infty}T_3
\end{equation}
and
\begin{equation*}
T_2=\frac{1}{\sqrt{p}}\tr(\bbI_p-\Sigma_1\Sigma_2^{-1})+\frac{1}{\sqrt{p}}\log|\Sigma_1\Sigma_2^{-1}|<0,\quad
T_3=-\frac{1}{\sqrt{p}}(\bbmu_1-\bbmu_2)^T\Sigma_2^{-1}(\bbmu_1-\bbmu_2)\leq 0.
\end{equation*}
Here ``$\xrightarrow{i.p}$'' denotes convergence in probability. The parameter $\psi>0$ is given by
\[
\psi^2=(m_4-3)(1-2M_1+M_2)+2\left(\frac{1}{1-c_1}-2M_1+M_5+\frac{c_2}{1-c_2}M_1^2\right),
\]
where $M_1,M_2,M_5$ are given in Condition \ref{cond3}.
\end{lem}

Similarly, Lemma~\ref{lema3} holds by swapping $\Sigma_1$ with $\Sigma_2$ and $n_1$ with $n_2$ in Lemma \ref{lema2}.
\begin{lem}\label{lema3}
When $\bbz$ belongs to class 2, under Conditions \ref{cond1}-\ref{cond3}, we have
\begin{eqnarray*}
\pr_{1|2}^G\xrightarrow{i.p} 1-\Phi\left(\frac{\widetilde{T}}{\widetilde{\psi}}\right),
\end{eqnarray*}
where $\widetilde{T}$ and $\widetilde{\psi}$  are calculated by swapping $\Sigma_1$ with $\Sigma_2$, $M_1$ with $M_3$, $M_2$ with $M_4$, $M_5$ with $M_6$ and $c_1$ with $c_2$ in the expressions of $T$ and $\psi$ in Lemma \ref{lema2}.
\end{lem}

Combing the above two lemmas, we can directly conclude the asymptotic property of the misclassification rate of (\ref{11.1}) in Theorem \ref{themnew}.
\begin{thm}\label{themnew}
Under Conditions \ref{cond1}-\ref{cond3}, the misclassification rate of our generalized QDA
\[
R^G=\frac{1}{2}\left[\pr_{2|1}^G+\pr_{1|2}^G\right]\xrightarrow{i.p}
1-\frac{1}{2}\left[\Phi\left(\frac{T}{\psi}\right)+\Phi\left(\frac{\widetilde{T}}{\widetilde{\psi}}\right)\right],
\]
where $(T, \psi)$ and $(\widetilde{T},\widetilde{\psi})$ are given in Lemma \ref{lema2} and Lemma \ref{lema3} respectively.
\end{thm}
To better understand $R^G$, we consider specific cases in Corollary \ref{coronew} below. Let $\lambda_i>0$ $(i=1,\cdots,p)$ be the eigenvalues of the matrix $(\Sigma_1\Sigma_2^{-1})$, and thus $\lambda_i^{-1}$ $(i=1,\cdots,p)$ are the eigenvalues of the matrix $(\Sigma_2\Sigma_1^{-1})$. Denote
\[
s=\#\{\lambda_i\neq 1, \quad i=1,\cdots,p\},\quad
s(\epsilon)=\#\{|\lambda_i-1|>\epsilon>0, \quad i=1,\cdots,p\},
\]
where $\epsilon>0$ is any positive constant. There exist an $\gamma=\gamma(\epsilon)>0$ such that
\[
\#\{|\frac{1}{\lambda_i}-1|>\gamma>0, \quad i=1,\cdots,p\}=s(\epsilon).
\]
The parameters $s$ and $s(\epsilon)$ describe the deviation between $\Sigma_1$ and $\Sigma_2$. For example, if $\Sigma_1=\Sigma_2$, then $s=s(\epsilon)=0$; and if $\Sigma_1=\kappa \Sigma_2$, $\kappa>0$ is a constant, then $s=s(\epsilon)=p$.

\begin{coro}\label{coronew}
Suppose Conditions \ref{cond1}-\ref{cond3} hold. Consider the following three terms:
\begin{equation}\label{0103.1}
\zeta_1=\frac{1}{\sqrt{p}}\|\bbmu_1-\bbmu_2\|^2,\qquad  \zeta_2=\frac{s}{\sqrt{p}},\qquad
\zeta(\epsilon)=\frac{s(\epsilon)}{\sqrt{p}}.
\end{equation}
\begin{itemize}
\item[(i)] When either $\zeta_1$ or $\zeta(\epsilon)$ (given any $\epsilon>0$) diverges to infinity, we have  $$R^G\xrightarrow{i.p}0.$$
\item[(ii)] When both $\zeta_1$ and $\zeta_2$ degenerate to zero, we have  $$R^G\xrightarrow{i.p}\frac{1}{2},\quad \text{random guessing}.$$
\item[(iii)] If both $\zeta_1$ and $\zeta_2$ are bounded, and at least one of $\zeta_1$ and $\zeta(\epsilon)$ does not degenerate to zero, we have $$R^G\xrightarrow{i.p}\widetilde{C}\in (0,\frac{1}{2}).$$
\end{itemize}
\end{coro}

\begin{rmk}
Corollary \ref{coronew} implies that when the Euclidean norm of the mean difference between the two classes is of a larger order than $p^{1/4}$, regardless of the covariance matrices, the misclassification rate converges to zero. On the other hand, if the difference between the two covariance matrices is significant in the sense that $\zeta(\epsilon)\rightarrow \infty$, then regardless of the mean vectors,  the misclassification rate also tends to zero. However, if the two classes are too close to each other in the sense that both $\zeta_1$ and $\zeta_2$ degenerate to zero,  the classification rule behaves as random guessing, which is reasonable. In between these two extreme cases, we show that $R^G$ tends to some constant between $(0,1/2)$, which is not surprising.
\end{rmk}

\begin{rmk}
Based on the proof of Theorem \ref{themnew}, we also derive the asymptotic misclassification rate for the sample QDA \eqref{yq4} under moderate dimension. The detailed theoretical results and some plots for easy comparison are deferred to Supplement \ref{compare2}.
Note that $\mathcal{E}_{2|1}^S$ in Section \ref{interp} is related to the parameter $T_S$ in Proposition \ref{prop2}, specifically,  $\mathcal{E}_{2|1}^S=-\sqrt{p}T_S$.
\end{rmk}

\subsection{Comparison with the optimal QDA}\label{optim}
In this section, we compare the generalized QDA (\ref{11.1}) with its oracle version (\ref{yq1}) in terms of misclassification rates. In particular, we find that their limits are the same in the first two cases of Corollary \ref{coronew}, but different in the last case.

We first derive the misclassification rate for the optimal QDA.
\begin{prop}\label{prop1}
Under Conditions \ref{cond1}-\ref{cond3}, the misclassification rate of the optimal QDA (\ref{yq1})
\[
R^O=\frac{1}{2}\left[\pr_{2|1}^O+\pr_{1|2}^O\right]\xrightarrow{i.p}
1-\frac{1}{2}\left[\Phi\left(\frac{T}{\psi_0}\right)+
\Phi\left(\frac{\widetilde{T}}{\widetilde{\psi}_0}\right)\right],
\]
where
$T$ and $\widetilde{T}$ are the same as the ones in Lemma \ref{lema2} and Lemma \ref{lema3} respectively. The parameters $\psi_0$ and $\widetilde{\psi}_0$ are nonnegative constants given by
\[
\psi_0^2=(m_4-3)(1-2M_1+M_2)+2\left(1-2M_1+M_5\right),
\]
\[
\widetilde{\psi}_0^2=(m_4-3)(1-2M_3+M_4)+2\left(1-2M_3+M_6\right).
\]

\end{prop}

Proposition \ref{prop1} implies the following analogue of Corollary \ref{coronew}.
\begin{coro}\label{coro3}
Under the same conditions and cases as in Corollary \ref{coronew}, we have

$(i)$ $R^O\xrightarrow{i.p}0.$
$(ii)$ $R^O\xrightarrow{i.p}\frac{1}{2}$, random guessing.
$(iii)$ $R^G\xrightarrow{i.p}C\in (0,\frac{1}{2})$.
\end{coro}

Based on Theorem \ref{themnew} and Proposition \ref{prop1}, we study the difference of misclassification rates between our generalized QDA and the optimal one, denoted as
\[
\diff\triangleq  R^G-R^O.
\]
We find that the difference converges to zero when the two underlying populations are either close enough or deviate enough, or converges to some strictly positive constant. The former is called as ``easy case'' (corresponding to $(i)$-$(ii)$ in Corollary \ref{coronew}), while the latter as ``hard case'' (corresponding to $(iii)$ in Corollary \ref{coronew}).

\begin{thm}\label{them3}
Suppose Conditions \ref{cond1}-\ref{cond3} hold, and $\psi_0^2$ and $\widetilde{\psi}_0^2$ in the limit of $R^O$ are nonzero. 	

\begin{itemize}
\item[(i)]\textbf{\emph{[``easy'' case]}} If conditions in (i) or (ii) of Corollary \ref{coronew} hold, we have $$\diff\xrightarrow{i.p} 0.$$
\item[(ii)] \textbf{\emph{[``hard'' case]}} If conditions in (iii) of Corollary \ref{coronew} hold, we have
{\small
$$
\diff\xrightarrow{i.p} \widetilde{C}-C=\frac{1}{2}\left[\Phi\bigg(\frac{T}{\psi_0}\bigg)
-\Phi\bigg(\frac{T}{\psi}\bigg)+\Phi\bigg(\frac{\widetilde{T}}{\widetilde{\psi}_0}\bigg)
-\Phi\bigg(\frac{\widetilde{T}}{\widetilde{\psi}}\bigg)\right]>0,
$$
}
where $T,\widetilde{T}$ are defined in Lemmas \ref{lema2} and \ref{lema3}, $(\psi,\widetilde{\psi})$ and $(\psi_0,\widetilde{\psi}_0)$ are given in
Theorem \ref{themnew} and Proposition \ref{prop1}, respectively.
\end{itemize}
\end{thm}
Note that in the hard case of Theorem \ref{them3}, the  difference of the misclassification rates Diff tends to a positive constant. This is due to the observation that
\begin{equation}\label{1807.10}
\psi^2-\psi_0^2=2\left(\frac{c_1}{1-c_1}+\frac{c_2}{1-c_2}M_1^2\right)>0,
\quad
\widetilde{\psi}^2-\widetilde{\psi}_0^2=2\left(\frac{c_2}{1-c_2}+\frac{c_1}{1-c_1}M_3^2\right)>0
\end{equation}
under moderate dimension. And the smaller $c_1$ and $c_2$ are, the smaller Diff is. We would like to mention that in the low-dimensional setting, both $c_1$ and $c_2$ are zeros and thus $\psi^2=\psi_0^2$, $\widetilde{\psi}^2=\widetilde{\psi}_0^2$, which implies that Diff also converges to zero in the hard case. Hence, this non-vanishing gap represents one of ``moderate dimension phenomena.'' Please also see Figures \ref{figg5}, \ref{figg6} for numerical evidence in the simulation section.

\section{Numerical studies}\label{nume}
In Section \ref{simu}, various simulations are conducted to compare the numerical performance of the three QDA rules mentioned above, while Section \ref{dc} aims at reducing the gap between our generalized QDA and the optimal one for the ``hard'' case.

\subsection{Performance of the QDA rules}\label{simu}
 Throughout this section, standardized $t(5)$ distribution is used to generate $\bbx_i^0$ ($i=1,\cdots, n_1$) and $\bby_j^0$ ($j=1,\cdots, n_2$). Note that its moments of order 5 or higher do not exist.

Denote a block-diagonal matrix composed by matrices $\bbA$ and $\bbB$ as {\bf blk}$(\bbA,\bbB)$. Let $\Sigma_1=\bbI_p$, the identity matrix. We adopt six choices for $\Sigma_2$ to investigate these rules' behavior under different alternatives.
\begin{itemize}
\item {\bf Case 1:} $\Sigma_2=2\bbI_p$; \qquad \qquad \quad\ {\bf Case 2:} $\Sigma_2=3\bbI_p$;
\item {\bf Case 3:} $\Sigma_2=U\Lambda_1 U^T$;\qquad\qquad  {\bf Case 4:} $\Sigma_2=U\Lambda_2 U^T$;
\end{itemize}
Here $\Lambda_1$ and $\Lambda_2$ are diagonal matrices with diagonal elements drawn uniformly from $(1.5,2.5)$ and $(2.5,3.5)$, respectively.   $U$ is an orthogonal matrix. In the simulation, we generate it by selecting the eigenvector matrix of $\frac{1}{n_1}\bbZ\bbZ^T$, where the entries of $\bbZ_{p\times n_1}$ are i.i.d $N(0,1)$.
\begin{itemize}
\item {\bf Case 5:} $\Sigma_2=\text{\bf blk}(4\bbI_{3\lfloor\sqrt{p}\rfloor},\bbI_{p-3\lfloor\sqrt{p}\rfloor})$;
\qquad {\bf Case 6:} $\Sigma_2=\text{\bf blk}(5\bbI_{3\lfloor\sqrt{p}\rfloor},\bbI_{p-3\lfloor\sqrt{p}\rfloor})$. Here $\lfloor a\rfloor$ rounds $a$ to the nearest integer less than or equal to $a$.
\end{itemize}
Cases 1 and 2 consider the scenario when the difference $(\Sigma_2-\Sigma_1)$ is sparse (each row only has one nonzero element), while cases 3 and 4 consider the situation when the difference $(\Sigma_2-\Sigma_1)$ is not sparse (each row has $p$ nonzero elements). These four cases correspond to the ``easy'' case in Theorem \ref{them3}, where $s=s(\epsilon)=p$ and $\zeta_1=0$, $\zeta_2=\zeta(\epsilon)=\sqrt{p}$. While
 cases 5 and 6 correspond to the ``hard'' case, where $s=s(\epsilon)=3\sqrt{p}$ and $\zeta_1=0$, $\zeta_2=\zeta(\epsilon)=3$. The difference between two classes increases from Case 1 (3, 5) to Case 2 (4, 6, resp.).

\begin{figure}[!htp]
    \centering
    \includegraphics[height = 4in,width=5.5in]{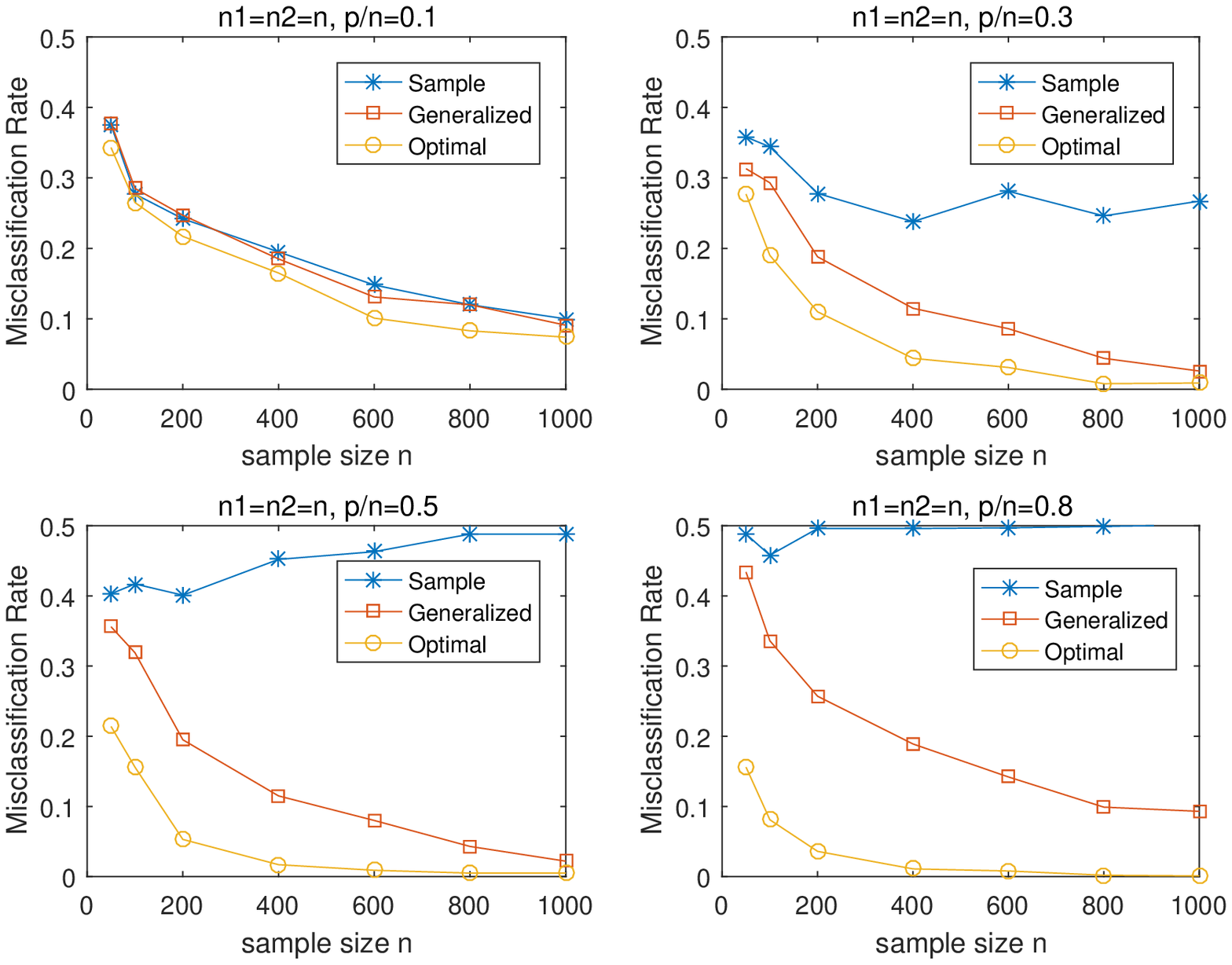}
    \caption{\small{\it Misclassification rates for Case 1 with $\bbmu_1=\bbmu_2=0$, $n_1=n_2=n$ and 1000 replications.}}
    \label{figg1}
    \includegraphics[height = 4in,width=5.5in]{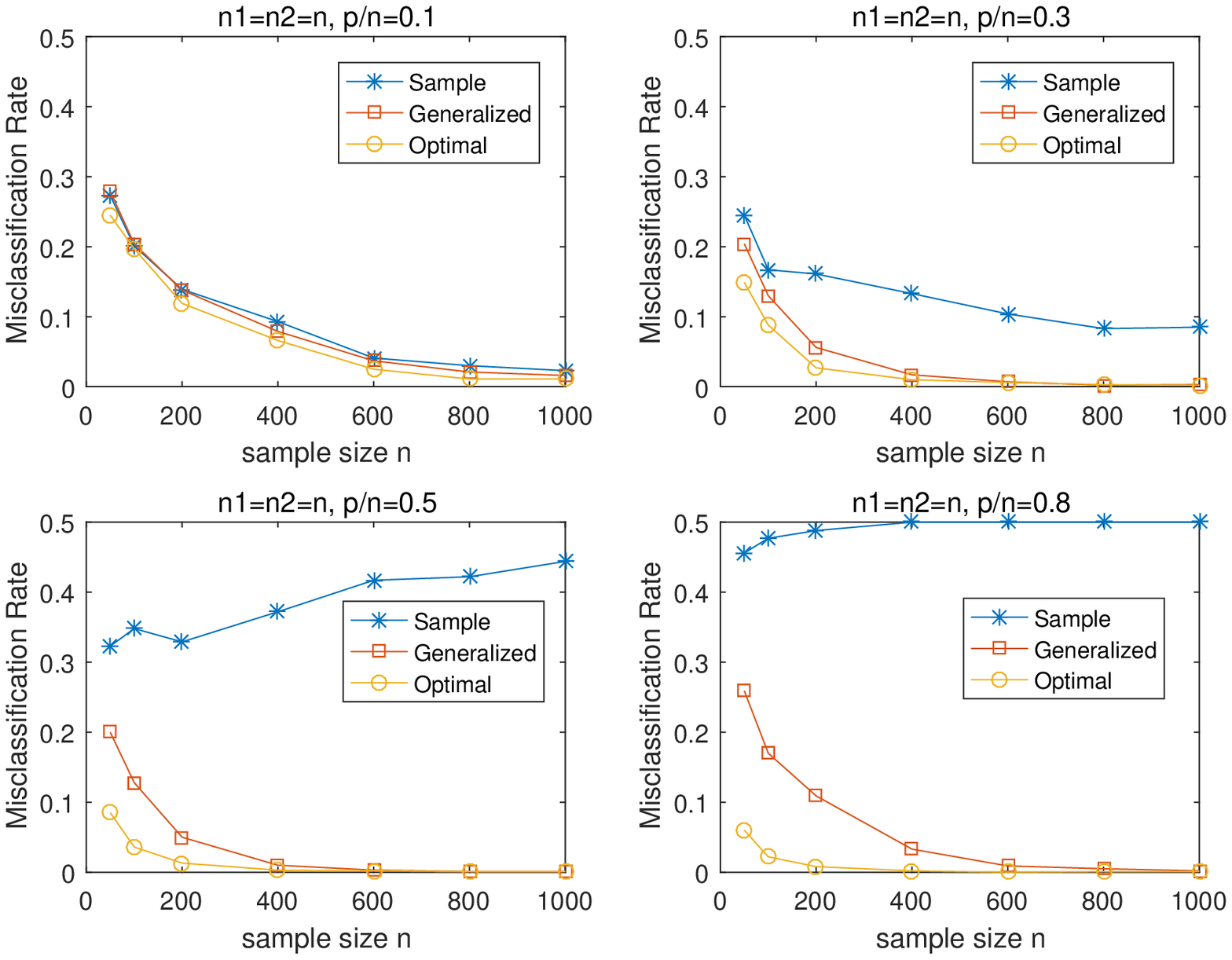}
    \caption{\small{\it Misclassification rates for Case 2 with $\bbmu_1=\bbmu_2=0$, $n_1=n_2=n$ and 1000 replications.}}
    \label{figg2}
  \end{figure}

\begin{figure}[!htp]
    \centering
    \includegraphics[height = 4in,width=5.5in]{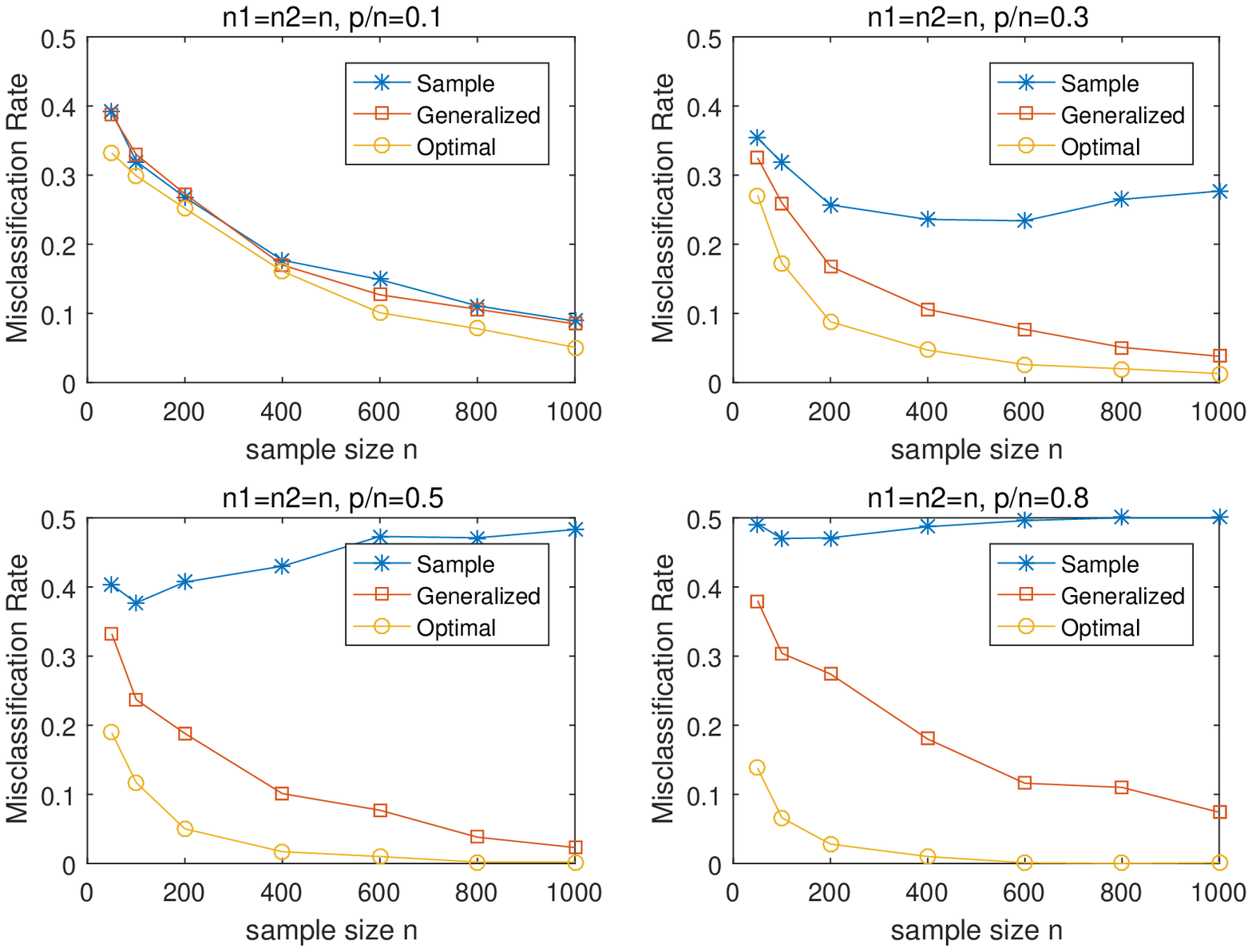}
    \caption{\small{\it Misclassification rates for Case 3 with $\bbmu_1=\bbmu_2=0$, $n_1=n_2=n$ and 1000 replications.}}
    \label{figg3}
    \includegraphics[height = 4in,width=5.5in]{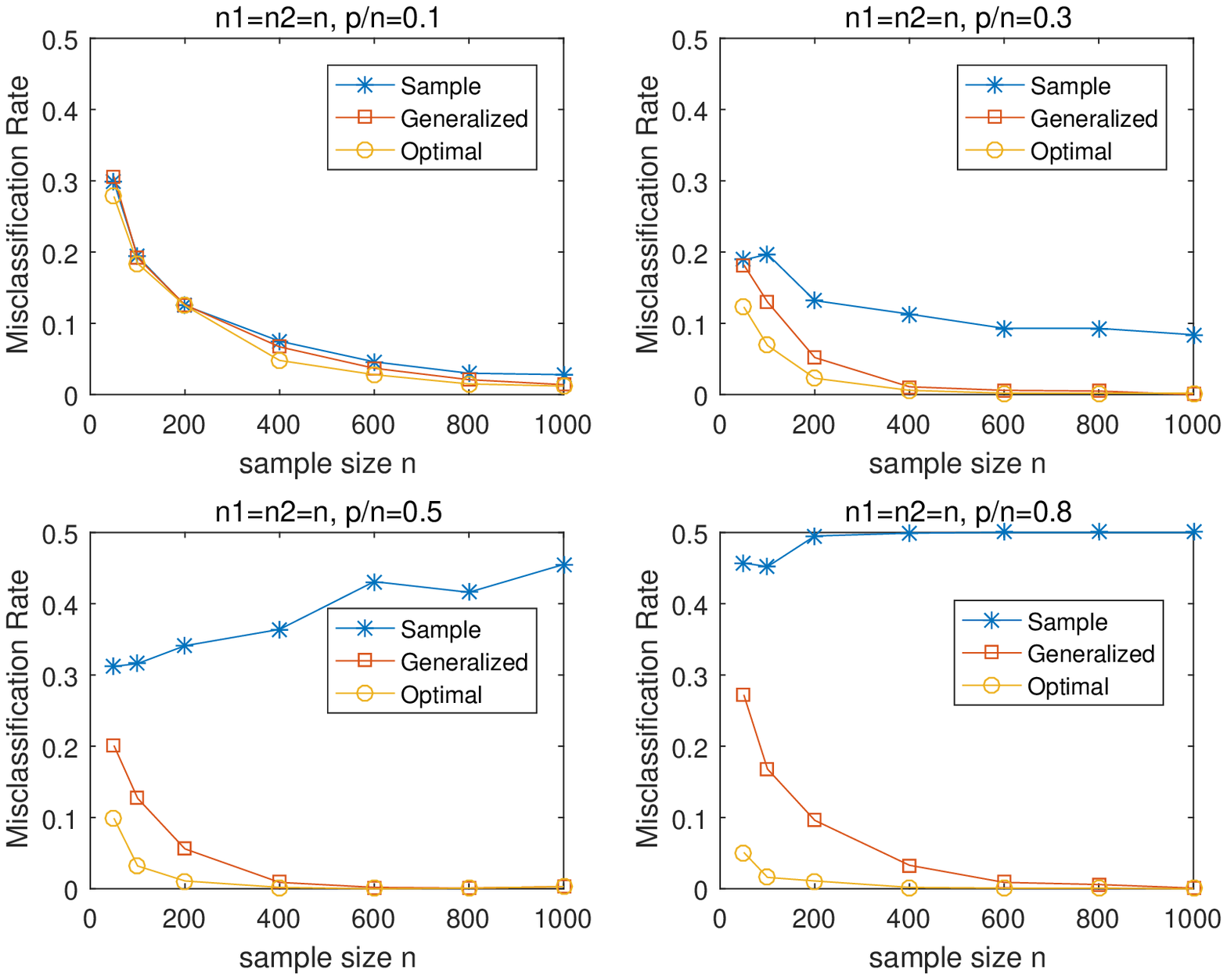}
    \caption{\small{\it Misclassification rates for Case 4 with $\bbmu_1=\bbmu_2=0$, $n_1=n_2=n$ and 1000 replications.}}
    \label{figg4}
  \end{figure}

\begin{figure}[!htp]
    \centering
    \includegraphics[height = 4in,width=5.5in]{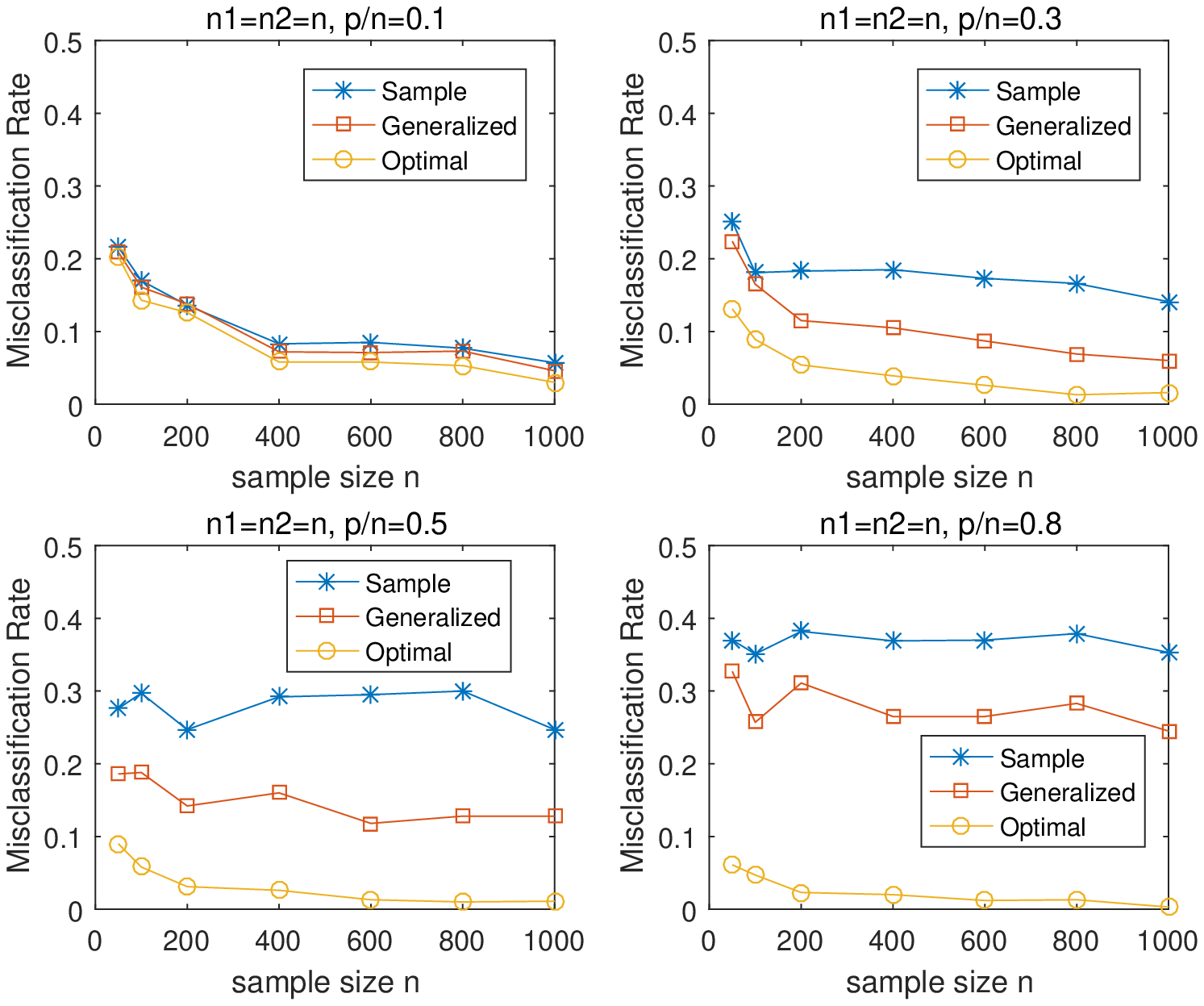}
    \caption{\small{\it Misclassification rates for Case 5 (``hard'' case) with $\bbmu_1=\bbmu_2=0$, $n_1=n_2=n$ and 1000 replications.}}
    \label{figg5}
    \includegraphics[height = 4in,width=5.5in]{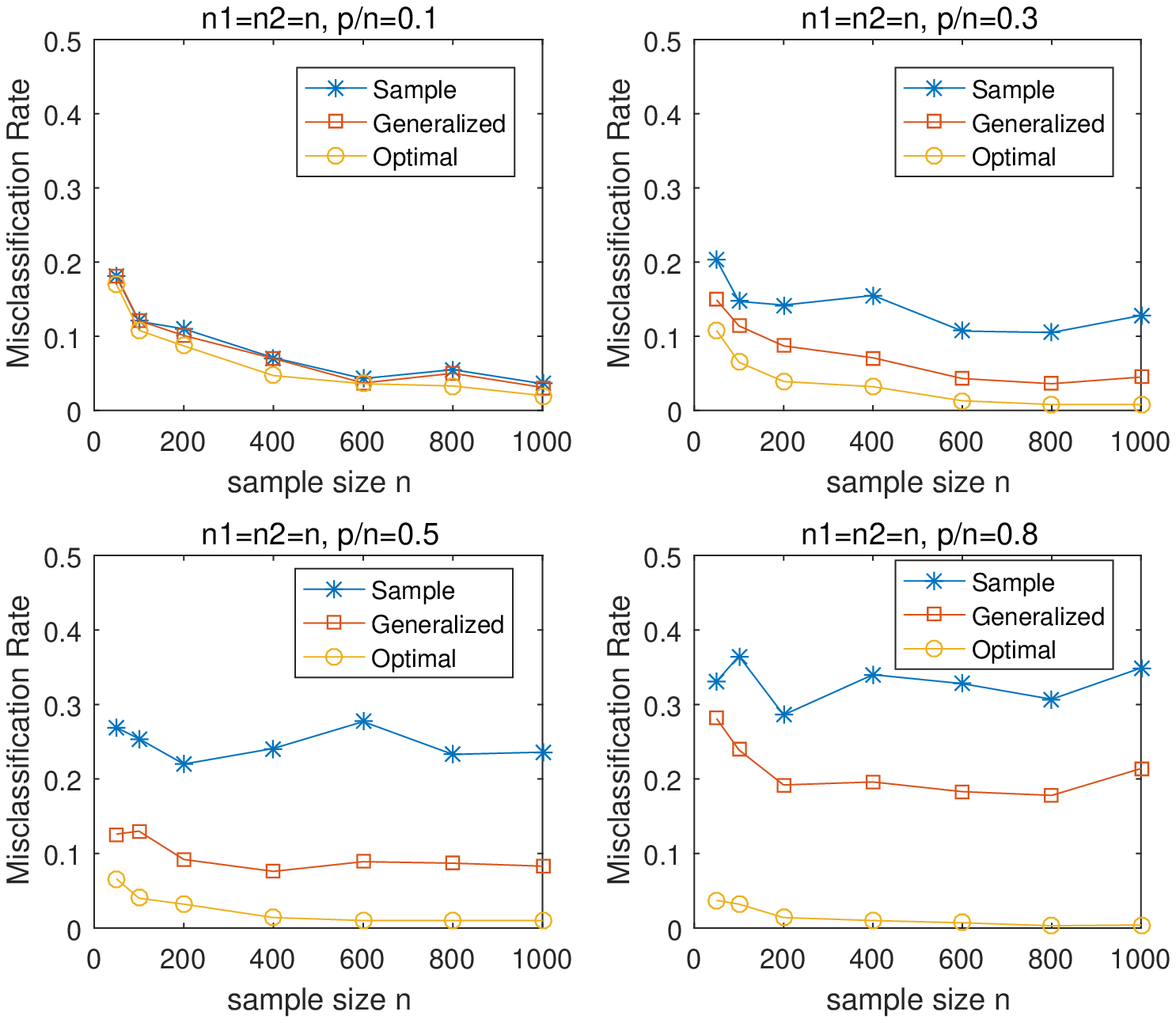}
    \caption{\small{\it Misclassification rates for Case 6 (``hard'' case) with $\bbmu_1=\bbmu_2=0$, $n_1=n_2=n$ and 1000 replications.}}
    \label{figg6}
  \end{figure}

First let the mean vectors and sample sizes be equivalent. Based on 1000 replications, the estimated misclassification rates for six cases are displayed in Figures \ref{figg1}-\ref{figg6}, where the sample size varies from 50 to 1000. The numerical performance from the six figures can be summarized as follows:
\begin{itemize}
 \item[(1)]The larger the ratio $p/n$ is, the worse the sample QDA performs. When the sample QDA behaves as well as the optimal one (see $p/n=0.1$), our generalized QDA also maintains a similar nice property. When $p/n$ equals to 0.5 or 0.8, our generalized QDA significantly improves over the sample QDA which behaves like random guessing.
 \item[(2)] When the difference between $\Sigma_1$ and $\Sigma_2$ increases from Figure \ref{figg1} (Figure \ref{figg3} or Figure \ref{figg5}) to Figure \ref{figg2} (Figure \ref{figg4} or Figure \ref{figg6}, resp.), both $R^O$ and $R^G$ decrease significantly in all plots of $p/n$. This is nature since the two classes are more separated.  However, when the ratio is large (see $p/n=0.8$), $R^S$ does not share this trend and may be even worse.
 \item[(3)] In the ``easy'' cases (Figures \ref{figg1}-\ref{figg4}), the distance between $R^G$ and $R^O$ approaches to zero as $n$ becomes larger. While in the ``hard'' cases (Figures \ref{figg5} and \ref{figg6}), there exists a non-vanishing gap. This  is consistent with the Theorem \ref{them3}.

 \end{itemize}

\begin{figure}[!htp]
    \centering
    \includegraphics[height = 4in,width=5.5in]{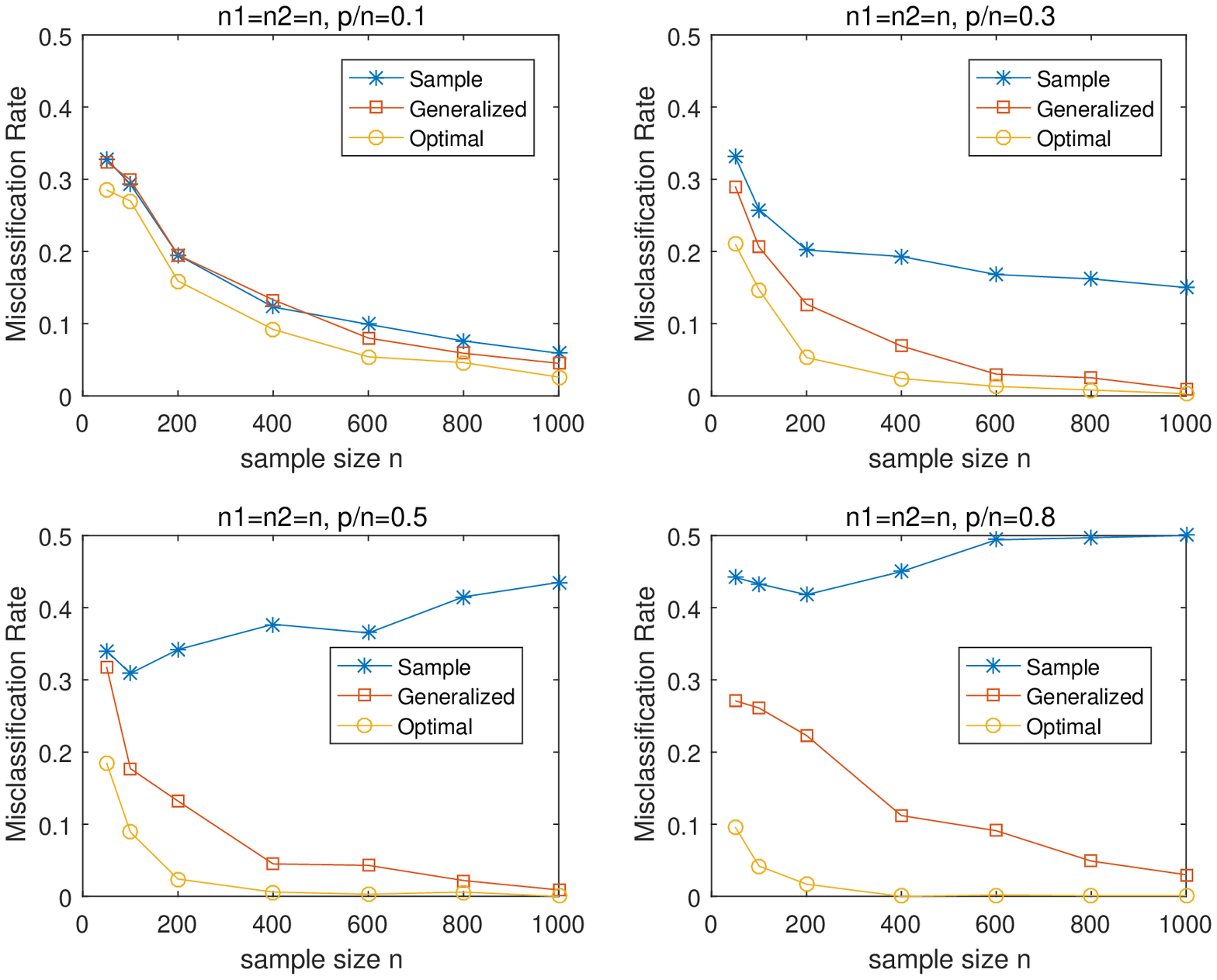}
    \caption{\small{\it (Unequal means). Misclassification rates for Case 1 with $\bbmu_1=0$, $\bbmu_2=\bbmu$, $n_1=n_2=n$ and 1000 replications.}}
    \label{figg7}
    \includegraphics[height = 4in,width=5.5in]{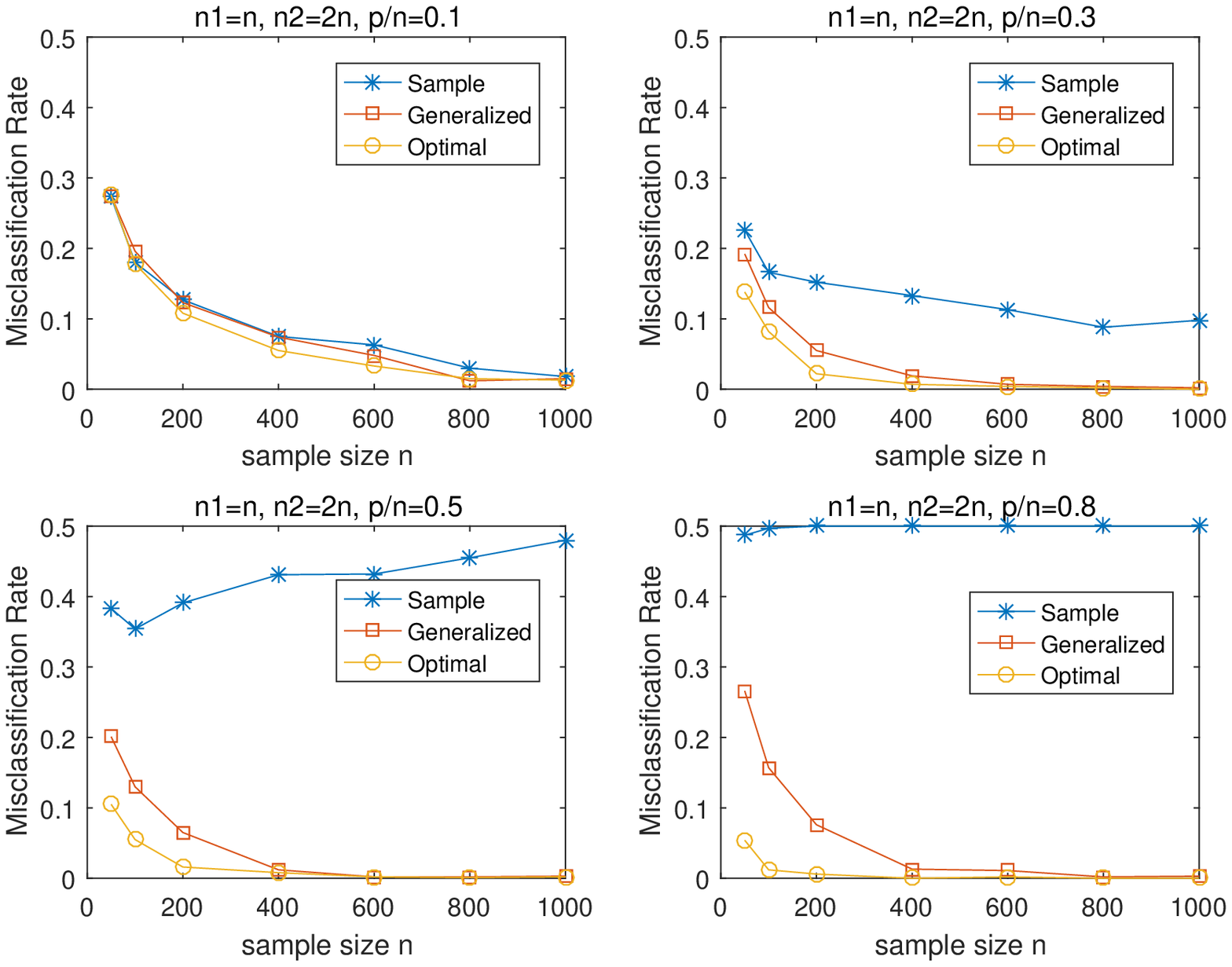}
    \caption{\small{\it (Unequal sizes). Misclassification rates for Case 2 with $\bbmu_1=\bbmu_2=0$, $n_2=2n_1=2n$ and 1000 replications.}}
    \label{figg8}
  \end{figure}

We next consider unequal mean vectors and unequal sample sizes. We only show the plots of one covariance matrix case for each of them and others behave similarly. In the setting with unequal mean vectors for Case 1 in Figure \ref{figg7},  $\bbmu_1$ is still a zero vector while the entries of $\bbmu_2=\bbmu\neq \bbzero$ are drawn uniformly from $(-0.6,0.6)$. Comparing it with Figure \ref{figg1}, in general we can observe that the misclassification rates in Figure \ref{figg7} are slightly smaller than the corresponding ones in Figure \ref{figg1}. One may numerically compare the values from Table \ref{tab1}, which records the simulated values. This is reasonable because with unequal mean vectors, the two classes are more separated from each other. In the setting with unequal sample sizes for Case 2 in Figure \ref{figg8}, we choose $n_2=2n_1$. It behaves in a similar way to Figure \ref{figg2}.

\begin{figure}[!htp]
    \centering
    \includegraphics[height = 4in,width=5.5in]{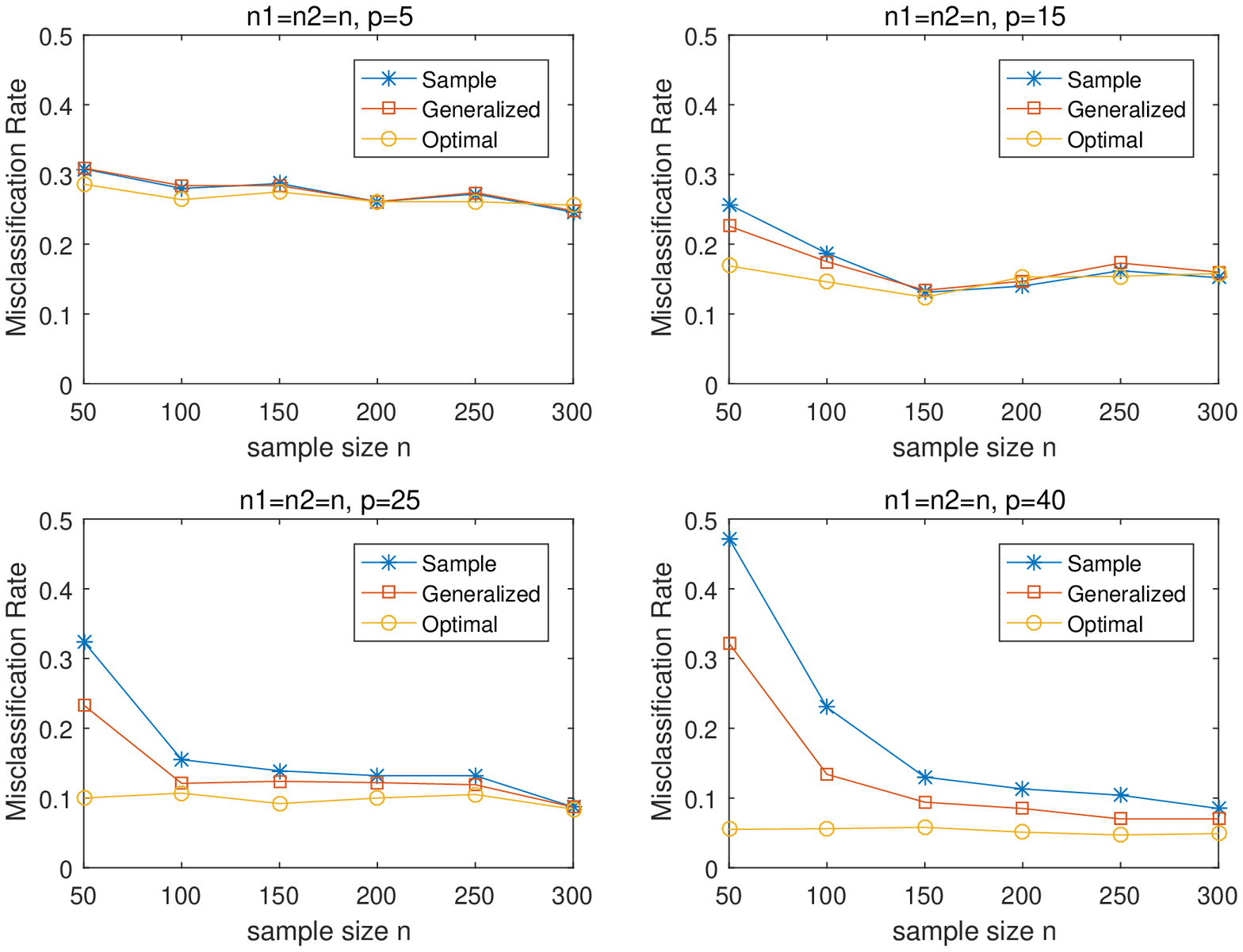}
    \caption{\small{\it ($p$ is fixed). Misclassification rates for Case 4 with $\bbmu_1=\bbmu_2=0$, $n_1=n_2=n$ and 1000 replications.}}
    \label{figg9}
  \end{figure}

Finally, although our theoretical results are based on the assumption that $p$ tends to infinity together with the sample sizes, one may be interested in the performance under the fixed dimension. To this end, in Figure \ref{figg9}, the estimated misclassification rates when $p$ is fixed at $5$, $15$, $25$ or $40$ are plotted. When $p=5$, the three curves are essentially coincident with each other. As $p$ increases, the generalized QDA shows remarkable improvement than the sample one, especially when the sample sizes are not significantly larger than $p$.

\subsection{Proposals for the ``hard'' case}\label{dc}
In this section, we intend to narrow the gap in the ``hard'' case by adapting the divide-and-conquer method. Conventionally, divide-and-conquer is done over samples, that is, for each class,
the samples are divided into non-overlapping subgroups and the final decision is made by averaging or majority voting over these subgroups. As an initial attempt, this conventional sample splitting trick is proven not to work in our case, either empirically or theoretically. For the constraint of space, we defer them to Supplement \ref{negativ} for details.

Unlike the conventional approach, we propose a new modification -- {\bf divide-and-conquer over dimension} -- in the following two ways:

\vspace{10pt}
\begin{itemize}
\item {\bf Method 1: Subgroup screening over dimension}

\begin{itemize}
\item Step 1: Divide the dimension $p$ into $H$ subgroups, each with $p_0$ components, i.e. $H=\lfloor \frac{p}{p_0}\rfloor$.
\item Step 2: For each $i=1,\cdots, H$, consider the generalized QDA \eqref{11.1} with data dimension in the $i$-th subgroup and denote its left side by $T_{1i}$, right side by $T_{2i}$. Note that the values in \eqref{1805.1} are calculated by replacing $p$ with $p_0$.
\item Step 3: Locate the index $I$ such that $I=\arg\max\limits_{1\leq i\leq H}|T_{1i}-T_{2i}|$.
\item Step 4: Identify the class label of the new observation based on our generalized QDA with dimension in the $I$-th group.
\end{itemize}
\end{itemize}

\begin{itemize}
\item {\bf Method 2: Component-wise screening over dimension}
\begin{itemize}
\item Step 1: For each component $j=1,\cdots, p$, consider the generalized QDA \eqref{11.1} with data entries in the $j$th component (e.g, $\bar\bbx$ is replaced by $\bbe_j^T\bar\bbx$, where $\bbe_j$ is a $p$-variate unit vector with the $j$th entry being one) and denote its left side by $T_{1j}$, right side by $T_{2j}$. In this case the values in \eqref{1805.1} are calculated by replacing $p$ with $1$.
\item Step 2: Define the index set $I'$ as
$$I'=\{1\leq j\leq p: |T_{1j}-T_{2j}| \ \text{is among the first $p_0$ largest of all}\}.$$
\item Step 3: Identify the class label of the new observation based on our generalized QDA with dimension in the index set $I'$.
\end{itemize}
\end{itemize}

\begin{figure}[!htp]
    \centering
    \includegraphics[width=5.5in,height=1.85in]{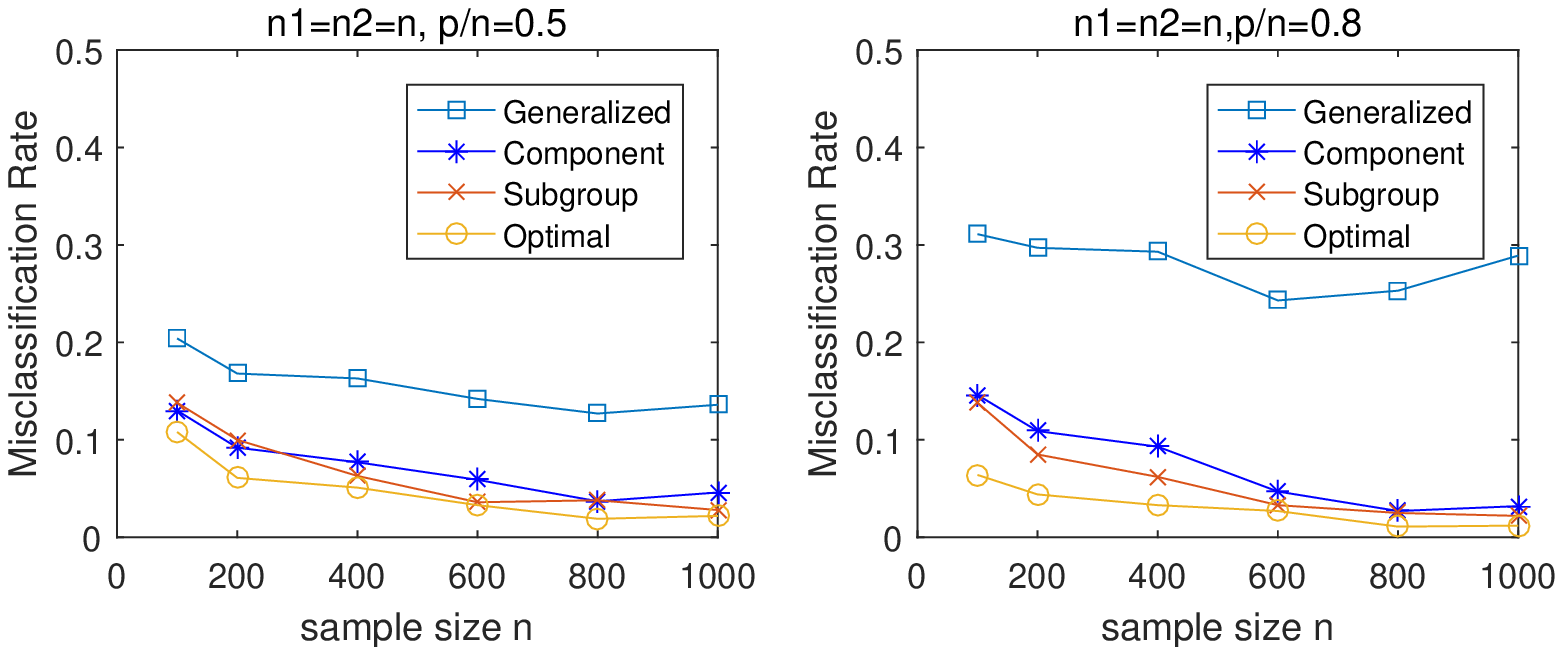}
    \caption{\small{\it Misclassification rates for Case 5. Apply divide-and-conquer over dimension (``Subgroup''(Method 1) and ``Component''(Method 2)) with $p_0=3\lfloor \sqrt{p}\rfloor$.}}
    \label{figg10}
\end{figure}

\begin{figure}[!htp]
    \centering
    \includegraphics[width=5.5in,height=1.85in]{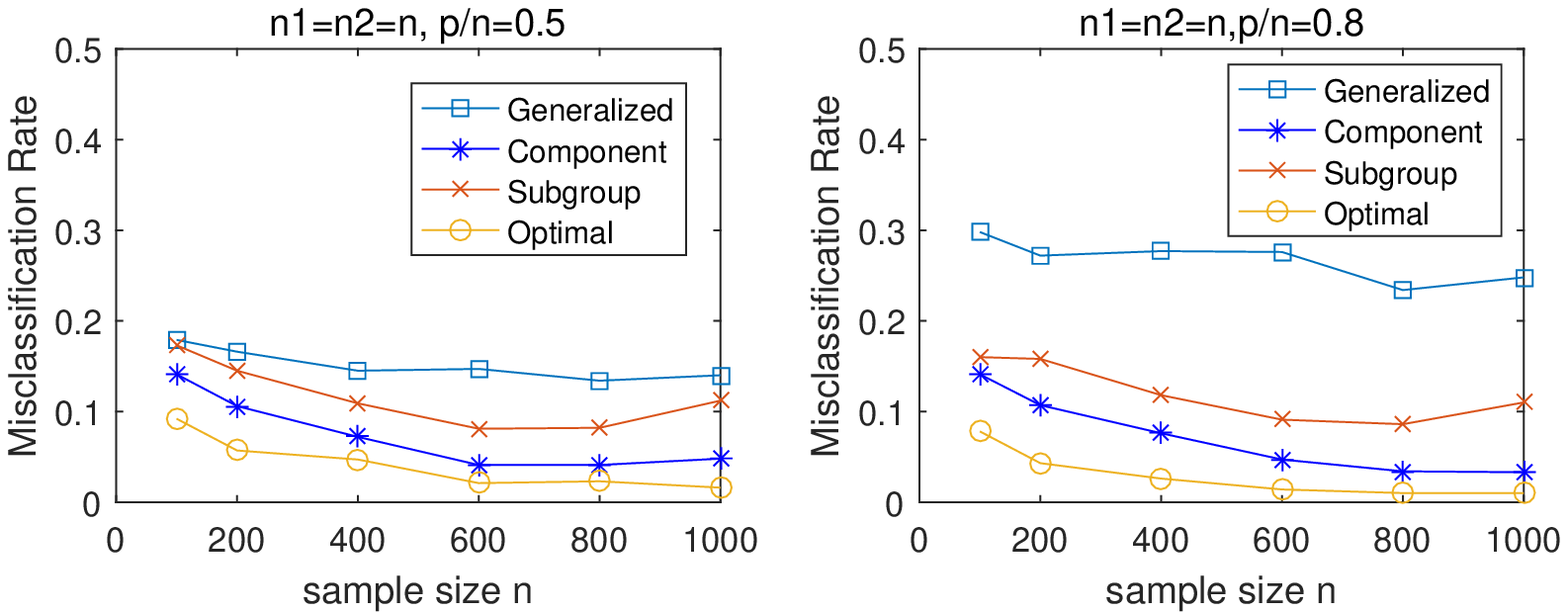}
    \caption{\small{\it Misclassification rates for Case 7. Apply divide-and-conquer over dimension (``Subgroup''(Method 1) and ``Component''(Method 2)) with $p_0=3\lfloor \sqrt{p}\rfloor$.}}
    \label{figg11}
      \end{figure}

Figure \ref{figg10} visualizes the performance after applying these two modification methods (marked as ``Subgroup''(Method 1) and ``Component''(Method 2)) by re-visiting the hard case 5 above.
 Both methods display a significant improvement over our original generalized QDA rule and the gap with the optimal one becomes negligible for large sample sizes. Moreover, under this hard case 5, one may observe that Method 1 is slightly better than Method 2. This is reasonable since in case 5, the different entries between the two covariance matrices are clustered together. Then it is more possible for the subgroup screening to select most of the significant components. To remove the clustering property of case 5,  we randomly select $3\lfloor \sqrt{p}\rfloor$ out of the $p$ diagonal entries of $\Sigma_1$ and assign the same value 4 -- as in case 5 -- to them. This alternative one is named as case 7 and Figure \ref{figg11} plots the misclassification rates. Similar to case 5, both the two modification methods improve the performance, but different from case 5, the Method 2 now outperforms Method 1.

\begin{figure}[!htp]
    \centering
    \includegraphics[width=5.5in,height=1.85in]{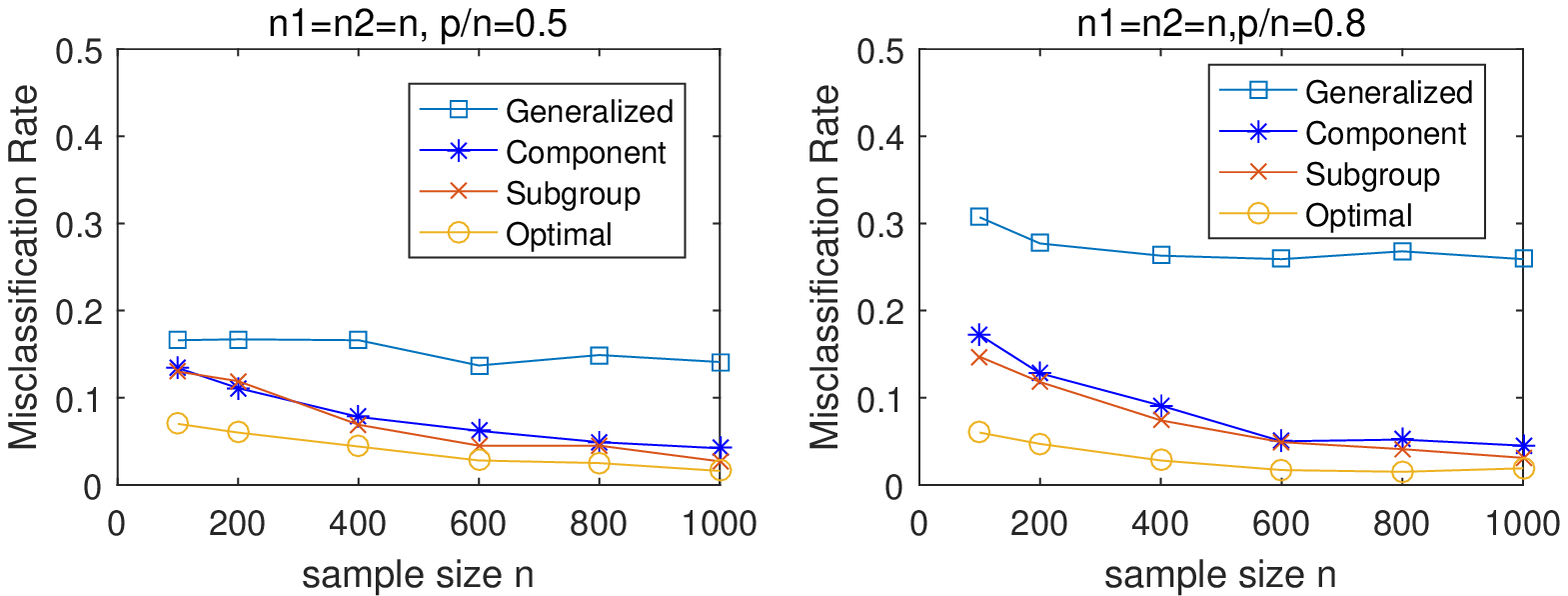}
    \includegraphics[width=5.5in,height=1.85in]{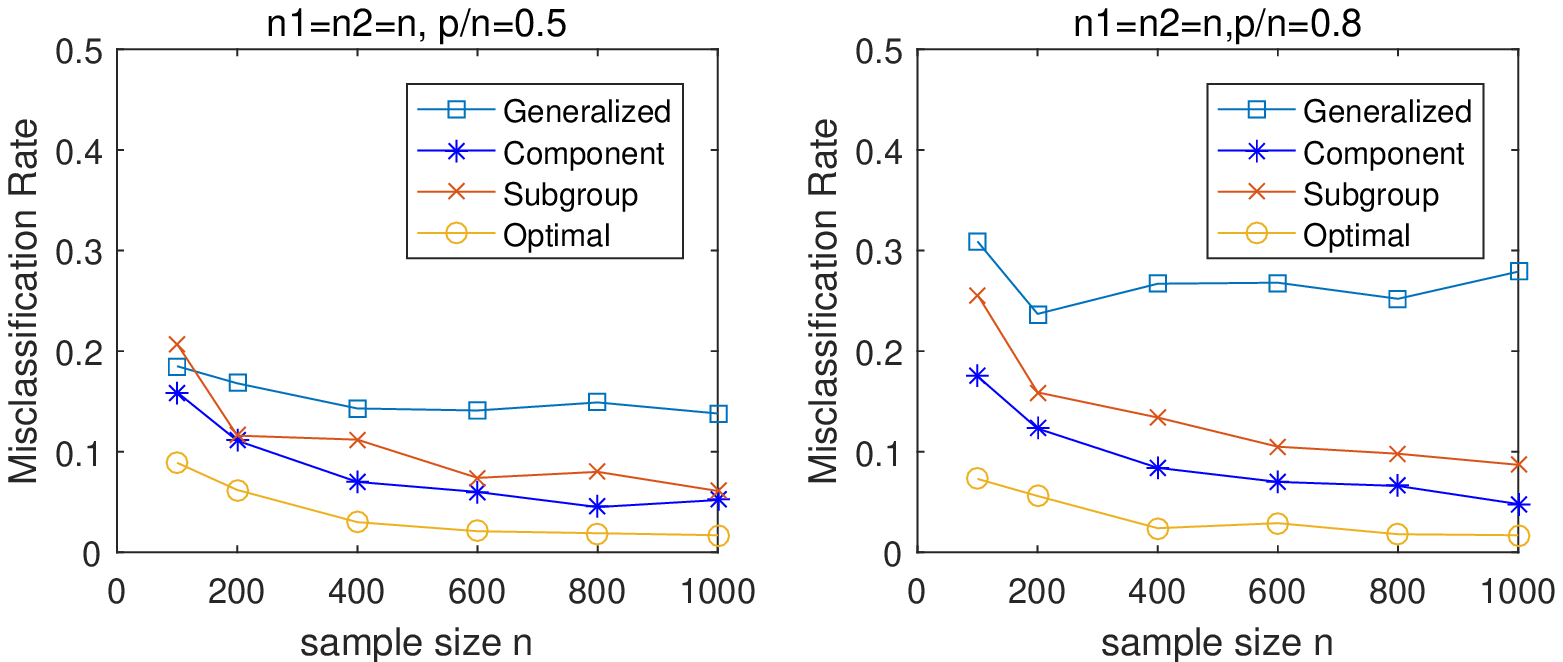}
    \caption{\small{\it Misclassification rates for Case 5 (upper row) and Case 7 (lower row). Apply divide-and-conquer over dimension (``Subgroup''(Method 1) and ``Component''(Method 2)) with $p_0=5\lfloor \sqrt{p}\rfloor$.}}
    \label{figg12}
  \end{figure}

Note that in the above two figures, we select $p_0=3\lfloor \sqrt{p}\rfloor$. This order is due to the condition of ``hard'' case (case $(iii)$ in Corollary \ref{coronew}), which roughly implies a significant $O(\sqrt{p})$ fraction out of the $p$ components. As for the constant $3$, it is related to the magnitude of $p$. One may choose a constant value that makes both the subgroup size ($p_0$) and the number of subgroups ($H$) not be too small. Figure \ref{figg12} presents the performance if we change the coefficient from 3 to 5. It shows similar phenomenon as before.

\begin{rmk}(An intuitive interpretation for the positive result of divide-and-conquer): In view of the hard case's condition in Theorem \ref{them3} - roughly speaking - the difference between the mean vectors or covariance matrices only appear in a small fraction of the $p$ dimensions. Considering the most significant fraction has two consequences: (i) discards much noisy information; (ii) decreases the ratio from $p/n$ to $p_0/n$. As a result of these two effects, the gap with the optimal one has been narrowed according to \eqref{1807.10}.
\end{rmk}

\begin{table}[!htp]
{\tiny
\begin{center}
\caption{\label{tab1}\it Misclassification rates for Case 1 with equal or unequal mean vectors, $n_1=n_2=n$}
\begin{tabular}{c|c c c c c c c||c c c c c c c}
\hline
&\multicolumn{14}{c}{$n_1=n_2=n$}\\
\hline
${\bf \frac{p}{n}}$
&\multicolumn{7}{c||}{$\bbmu_1=\bbmu_2$ (Figure \ref{figg1})}&\multicolumn{7}{c}{$\bbmu_1\neq\bbmu_2$ (Figure \ref{figg7})}\\
\hline
{\bf 0.1}&50&100&200&400&600&800&1000&50&100&200&400&600&800&1000\\
\hline
$R^S$&0.3750 &   0.2770 &   0.2420 &   0.1950  &  0.1480   & 0.1200  &  0.1000 &0.3280  &  0.2930   & 0.1950  &  0.1230 &   0.0990 &   0.0760  &  0.0590\\
$R^G$&0.3780 &   0.2850 &   0.2470  &  0.1850   & 0.1310  &  0.1200 &   0.0910 &0.3240 &   0.2990  &  0.1950  &  0.1330  &  0.0800&    0.0590  &  0.0450\\
$R^O$&0.3430  &  0.2650&    0.2170 &   0.1650&    0.1010 &   0.0830  &  0.0740&  0.2850 &   0.2700  &  0.1590  &  0.0920 &   0.0540  &  0.0460&    0.0260\\
\hline
\hline
{\bf 0.3}&50&100&200&400&600&800&1000&50&100&200&400&600&800&1000\\
\hline
$R^S$& 0.3580&    0.3440  &  0.2780&    0.2380  &  0.2810   & 0.2460 &   0.2670 &0.3310  &  0.2580   & 0.2020   & 0.1930    &0.1680  &  0.1620  &  0.1500\\
$R^G$& 0.3130  &  0.2930  &  0.1880 &   0.1150&    0.0860   & 0.0440  &  0.0260 & 0.2900 &   0.2060  &  0.1270 &   0.0690 &   0.0300  &  0.0250  &  0.0090\\
$R^O$& 0.2770 &   0.1900   & 0.1100 &   0.0440   & 0.0310  &  0.0080&    0.0090& 0.2100   & 0.1470   & 0.0530 &   0.0240   & 0.0130 &   0.0080  &  0.0030\\
\hline
\hline
{\bf 0.5}&50&100&200&400&600&800&1000&50&100&200&400&600&800&1000\\
\hline
$R^S$& 0.4030  &  0.4170  &  0.4010  &  0.4520  &  0.4630 &   0.4880 &   0.4880 &0.3400  &  0.3090 &   0.3420    &0.3770  &  0.3650   & 0.4150  &  0.4350\\
$R^G$&0.3560  &  0.3190 &   0.1950 &   0.1150 &   0.0800  &  0.0430   & 0.0220  &0.3180  &  0.1770   & 0.1320   & 0.0450    &0.0430   & 0.0220    &0.0090\\
$R^O$& 0.2140 &   0.1560  &  0.0530   & 0.0170 &   0.0090 &   0.0050 &   0.0050& 0.1850  &  0.0890 &   0.0240 &   0.0060 &   0.0030    &0.0060  &       0\\
\hline
\hline
{\bf 0.8}&50&100&200&400&600&800&1000&50&100&200&400&600&800&1000\\
\hline
$R^S$& 0.4880 &   0.4580   & 0.4960  &  0.4960  &  0.4970   & 0.4990 &   0.5010& 0.4420&    0.4330    &0.4180   & 0.4500  &  0.4940   & 0.4970    &0.5000\\
$R^G$&0.4340  &  0.3360 &   0.2570 &   0.1890 &   0.1420   & 0.0990 &   0.0930& 0.2710  &  0.2610   & 0.2230  &  0.1120    &0.0910 &   0.0490  &  0.0300\\
$R^O$&0.1560  &  0.0810&    0.0360 &   0.0110  &  0.0080   & 0.0020   & 0.0010 &0.0960 &   0.0420  &  0.0170&         0  &  0.0020  &  0.0010   & 0.0010\\
\hline
\end{tabular}
\end{center}
}
\end{table}


\newpage
\begin{appendices}
\section{Proof of the main results}\label{profs}
\subsection{Preliminary knowledge}
We first introduce some basic definitions in the random matrix theory.
\begin{deff}\label{def1}
For any $n\times n$ symmetric matrix $\bbA$ with real eigenvalues $\lambda_n(\bbA)\le...\le \lambda_2(\bbA) \le \lambda_1(\bbA)$,
the empirical spectral distribution (ESD) of $\bbA$ is defined by
$$F^{\bbA}(x)=\frac{1}{n}\sum_{i=1}^n \mathbf{1}_{\{\lambda_i(\bbA)\le x\}}.$$
The limit distribution of ESD is called the limiting spectral distribution (LSD).
\end{deff}

\begin{deff}\label{def2}
For any cumulative distribution function (c.d.f.) $F$, its Stieltjes transform is defined by
\[
s_F(z)=\int\frac{1}{x-z}dF(x),\qquad \Im z\neq 0.
\]
\end{deff}


Then let's look at a lemma that is used frequently in the proofs.
\begin{lem}[Theorem 7.2 in \cite{Bai08}]\label{lema1} Let $\{A_n=[a_{ij(n)}]\}_n$ be a sequence of $n\times n$ real symmetric matrices, $(x_i)_{i\in \mathbb{N}}$ be a sequence of i.i.d. $K$-dimensional real random vectors.  Write $x_i=(x_{1i},\cdots,x_{Ki})^T$. Assume $\E(x_i)=0$, $\E(x_1x_1^T)=(\gamma_{ij})$, $1\leq i,j\leq K$, and $\E[|x_{j1}|^4]<\infty$, $j=1,\cdots, K$. Let $X(l)=(x_{l1},\cdots,x_{ln})^T$, $1\leq l\leq K$. Assume the following limits exist
\begin{eqnarray*}
&&w=\lim_{n\rightarrow\infty}\frac{1}{n}\sum_{u=1}^n a_{uu}^2(n),\non
&&\theta=\lim_{n\rightarrow\infty}\frac{1}{n}\tr A_n^2.
\end{eqnarray*}
Then, the $M$-dimensional random vectors
\[
Z_n=(Z_{n,l}),\quad Z_{n,l}=\frac{1}{\sqrt{n}}[X(l)^TA_nX(l)-\gamma_{ll}\tr A_n],\quad 1\leq l\leq K,
\]
converge weakly to a zero-mean Gaussian vector with covariance matrix $D=D_1+D_2$ where
\begin{eqnarray*}
&&D_1=w(\E[x_{l1}^2x_{l'1}^2]-\gamma_{ll}\gamma_{l'l'}),\quad 1\leq l,l'\leq K,\non
&&D_2=(\theta-w)(\gamma_{ll'}\gamma_{l'l}+\gamma_{ll'}^2),\quad 1\leq l,l'\leq K.
\end{eqnarray*}
\end{lem}

As the last step before proceeding to the proofs of the results in the main paper, we develop and prove a proposition, which is crucial to the main results.
\begin{prop}\label{prop3}
Let $\bbC=\bbV^T\s^{-1}\bbV=(\C_{ij})_{p\times p}$, where $\bbV=(\bbv_1,\cdots, \bbv_p)$, $\s=\frac{1}{n-1}(\bbX-\bar\bbX)(\bbX-\bar\bbX)^T$, the entries in the matrix $\bbX=(X_{ij})_{p\times n}$ are i.i.d with mean zero, variance 1 and finite fourth moment and $\bar\bbx$ is the sample mean vector of the $n$ columns of $\bbX$, $\bar\bbX=\bar\bbx\cdot \bbone_{n_1}^T$. Assume that $p/n\rightarrow c\in (0,1)$, $\|\bbv_j\|< \infty$, $j=1,\cdots,p$ and
\[
\frac{1}{p}\sum\limits_{i=1}^p(\bbv_i^T\bbv_i)^2\rightarrow M.
\] Then
\[
\frac{1}{p}\sum\limits_{i=1}^p \C_{ii}^2\xrightarrow {i.p} \frac{1}{(1-c)^2}\cdot M.
\]
\end{prop}

Throughout the proof,  $\|A\|$ indicates the Euclidean norm if $A$ is a vector, and denotes the spectral norm if $A$ is a matrix.

\subsection{Proof of Proposition \ref{prop3}}
It is obvious that $\C_{ii}=\bbv_i^T\s^{-1}\bbv_i$, $i=1,\cdots,p$. Define $\bbS=\frac{1}{n-1}\bbX\bbX^T$, $C_{ii}=\bbv_i^T\bbS^{-1}\bbv_i$ and $C_{ii}(z)=\bbv_i^T\bbD^{-1}(z)\bbv_i$, where $\bbD(z)=\bbS+z\bbI$, $z>0$, $z=z(p)\rightarrow 0$ not too fast (slower than $p^{-l}$, $l$ can be any positive value, say $z=1/\log p$, $z=1/p$, etc).
The proof of Proposition \ref{prop3} is separated into five steps.
\begin{itemize}
\item {\bf Step 1}. Prove  $\frac{1}{p}\sum\limits_{i=1}^p C_{ii}^2-\frac{1}{p}\sum\limits_{i=1}^p C_{ii}^2(z)\xrightarrow{i.p} 0$.
\end{itemize}
We can write
\[
    C_{ii}-C_{ii}(z)=\bbv_i^T(\bbS^{-1}-\bbD^{-1}(z))\bbv_i=-\bbv_i^T\bbD^{-1}(z)(\bbS-\bbD(z))\bbS^{-1}\bbv_i=
    z\cdot \bbv_i^T\bbD^{-1}(z)\bbS^{-1}\bbv_i.
\]
Then
\begin{eqnarray*}
&&\left|\frac{1}{p}\sum\limits_{i=1}^p C_{ii}^2-\frac{1}{p}\sum\limits_{i=1}^p C_{ii}^2(z)\right|
=\left|\frac{1}{p}\sum\limits_{i=1}^p(C_{ii}+C_{ii}(z))(C_{ii}-C_{ii}(z))\right|\non
&=&z\cdot\frac{1}{p}\sum\limits_{i=1}^p\bbv_i^T(\bbS^{-1}+\bbD^{-1}(z))\bbv_i
\bbv_i^T\bbD^{-1}(z)\bbS^{-1}\bbv_i\non
&\leq&z\cdot\|\bbS^{-1}+\bbD^{-1}(z)\|\|\bbD^{-1}(z)\|\|\bbS^{-1}\|\cdot\frac{1}{p}\sum\limits_{i=1}^p\|\bbv_i\|^4.
\end{eqnarray*}
The conclusion is thus achieved from the conditions that $z\rightarrow 0$, $\|\bbv_i\|<\infty$ and the observation that (see \cite{BS10})
\begin{equation}\label{add1}
    \|\bbS^{-1}\|=O_p(1),\quad \|\bbD^{-1}(z)\|=O_p(1).
\end{equation}

Step 1 ensures that hereafter we can investigate $\frac{1}{p}\sum\limits_{i=1}^p C_{ii}^2(z)$ instead of $\frac{1}{p}\sum\limits_{i=1}^p C_{ii}^2$.

Next we replace $C_{ii}(z)$ by its  truncated and centralized version.
Define $\widehat{\bbS}=\frac{1}{n-1}\widehat{\bbX}\widehat{\bbX}^T$ with $\widehat{\bbX}_{p\times n}$  having $(i,j)$th entry $\widehat{X}_{ij}=X_{ij}I_{\{|X_{ij}|<\delta_p\sqrt{p}\}}$ and let
$\widetilde{\bbS}=\frac{1}{n-1}\widetilde{\bbX}\widetilde{\bbX}^T$ with $\widetilde{\bbX}_{p\times n}$  having $(i,j)$th entry $\widetilde{X}_{ij}=(\widehat{X}_{ij}-\E \widehat{X}_{ij})/\sigma_n$ and $\sigma_n^2=\E|\widehat{X}_{ij}-\E \widehat{X}_{ij}|^2$. Here $\delta_p\rightarrow 0$ so that
$$\delta_p^{-4}\int_{\{|X_{11}|\geq \delta_p\sqrt{p}\}}|X_{11}|^4\rightarrow 0.$$
One may refer to \cite{Bai04} for detailed illustrations of such truncation under the existence of fourth moments.  The notations $\widehat{\bbD}(z)$, $\widetilde{\bbD}(z)$ ($\widehat{C}_{ii}(z)$, $\widetilde{C}_{ii}(z)$) indicate the analogues of $\bbD(z)$ ($C_{ii}(z)$) with the matrix $\bbS$ replaced by $\widehat{\bbS}$ and $\widetilde{\bbS}$, respectively.

\begin{itemize}
\item {\bf Step 2}. Prove  $\frac{1}{p}\sum\limits_{i=1}^p C_{ii}^2(z)-\frac{1}{p}\sum\limits_{i=1}^p \widetilde{C}_{ii}^2(z)\xrightarrow{i.p} 0$.
\end{itemize}
We first have
\begin{eqnarray}\label{add2}
&&\pr\left\{\frac{1}{p}\sum\limits_{i=1}^p C_{ii}^2(z)\neq \frac{1}{p}\sum\limits_{i=1}^p \widehat{C}_{ii}^2(z)\right\}\leq \pr\{\bbS\neq \widehat{\bbS}\}\non
&\leq& np\cdot\pr(|X_{11}|>\delta_p\sqrt{p})\leq \frac{np}{\delta_p^4 p^2}\int_{|X_{11}|>\delta_p\sqrt{p}}|X_{11}|=o(1).
\end{eqnarray}

Then compare $\widehat{C}_{ii}(z)$ with $\widetilde{C}_{ii}(z)$. As in Step 1, we can write
\[
\widehat{C}_{ii}^2(z)-\widetilde{C}_{ii}^2(z)=-\bbv_i^T(\widehat{\bbD}^{-1}(z)+\widetilde{\bbD}^{-1}(z))\bbv_i\bbv_i^T\widetilde{\bbD}^{-1}(z)(\widehat{\bbS}-\widetilde{\bbS})\widehat{\bbD}^{-1}(z)\bbv_i
\]
and similar to \eqref{add1}, both $\|\widehat{\bbD}^{-1}(z)\|$ and $\|\widetilde{\bbD}^{-1}(z)\|$ are bounded with probability 1.
For the term $(\widehat{\bbS}-\widetilde{\bbS})$, we have
\begin{eqnarray*}
\|\widehat{\bbS}-\widetilde{\bbS}\|^2
&=&\frac{1}{(n-1)^2}\|\widehat{\bbX}\widehat{\bbX}^T
-\widetilde{\bbX}\widetilde{\bbX}^T\|^2
\leq \frac{2}{(n-1)^2}\|(\widehat{\bbX}-\widetilde{\bbX})(\widehat{\bbX}-\widetilde{\bbX})^T\|^2
+\frac{4}{(n-1)^2}\|\widetilde{\bbX}(\widehat{\bbX}-\widetilde{\bbX})^T\|^2\non
&\leq&\frac{2}{(n-1)^2}\Big(\tr(\widehat{\bbX}-\widetilde{\bbX})(\widehat{\bbX}-\widetilde{\bbX})^T\Big)^2
+\frac{4}{(n-1)^2}\|\widetilde{\bbX}\|^2\tr(\widehat{\bbX}-\widetilde{\bbX})(\widehat{\bbX}-\widetilde{\bbX})^T
=o_p(\frac{1}{n}),
\end{eqnarray*}
where the last step uses the result that $\frac{1}{n-1}\tr(\widehat{\bbX}-\widetilde{\bbX})(\widehat{\bbX}-\widetilde{\bbX})^T=o_p(n^{-1})$ by checking derivations on page 560 in \cite{Bai04} and the observation that $\frac{1}{n-1}\|\widetilde{\bbX}\|^2$ are bounded with probability 1. Therefore,
\[
\left|\frac{1}{p}\sum\limits_{i=1}^p \widehat{C}_{ii}^2(z)-\frac{1}{p}\sum\limits_{i=1}^p \widetilde{C}_{ii}^2(z)\right|\leq
\|\widehat{\bbD}^{-1}(z)+\widetilde{\bbD}^{-1}(z))\|\|\widetilde{\bbD}^{-1}(z)\|
\|\widehat{\bbS}-\widetilde{\bbS}\|\|\widehat{\bbD}^{-1}(z)\|
\cdot\frac{1}{p}\sum\limits_{i=1}^p\|\bbv_i\|^4=o_p(\frac{1}{\sqrt{n}}).
\]

Combing with \eqref{add2}, we get
$$
\frac{1}{p}\sum\limits_{i=1}^p C_{ii}^2(z)-\frac{1}{p}\sum\limits_{i=1}^p \widetilde{C}_{ii}^2(z)\xrightarrow{i.p} 0.
$$

Step 2 guarantees that  we can assume the underlying variables are truncated at $\delta_p\sqrt{p}$, centralized and renormalized. With these assumptions on $\bbX$, in the sequal, we still use $\bbX$, $\bbS$, $\bbD(z)$ and $C_{ii}(z)$ to ease notations.

\begin{itemize}
\item {\bf Step 3}.  Prove $\frac{1}{p}\sum\limits_{i=1}^p C_{ii}^2(z)-\frac{1}{p}\sum\limits_{i=1}^p(\E C_{ii}(z))^2\xrightarrow{i.p} 0$.
\end{itemize}
Denote $\bbX=(\bbx_1,\cdots,\bbx_n)$ and $\bbr_j=\frac{1}{\sqrt{n-1}}\bbx_j$, $j=1,\cdots,n$. Then $\bbS=\sum\limits_{j=1}^n \bbr_j\bbr_j^T$. Define
\[\bbD_j(z)=\bbD(z)-\bbr_j\bbr_j^T,\quad
\beta_j(z)=\frac{1}{1+\bbr_j^T\bbD_j^{-1}(z)\bbr_j},\quad
\E_j(\cdot)=\E(\cdot|\bbr_1,\cdots,\bbr_j),\quad
\E_0(\cdot)=\E(\cdot).
\]
We have $0<\beta_j(z)\leq 1$ and
\[
\bbD^{-1}(z)-\bbD_j^{-1}(z)=-\beta_j(z)\bbD_j^{-1}(z)\bbr_j\bbr_j^T\bbD_j^{-1}(z).
\]
Note that
\begin{equation}\label{add3}
C_{ii}^2(z)-(\E C_{ii}(z))^2=(C_{ii}(z)-\E C_{ii}(z))^2+2\E C_{ii}(z)(C_{ii}(z)-\E C_{ii}(z)).
\end{equation}
Since we already did truncation in Step 2, according to (1.9b) in \cite{Bai04}, for any positive $l$, whenever $0<a<\lim\inf_n\lambda_{\min}^T I_{(0,1)}(c)(1-\sqrt{c})^2$,
\[
\pr(\lambda_{\min}^{\bbS}\leq a)=o(p^{-l}).
\]
Therefore, we have
\begin{equation}\label{add4}
\E\|\bbD^{-m}(z)\|\leq \frac{1}{a^m}\pr(\lambda_{\min}^{\bbS}> a)+\frac{1}{|z|^m}\pr(\lambda_{\min}^{\bbS}\leq a) \quad \text{is bounded},\quad m=1,2,\cdots
\end{equation}
\begin{equation}\label{add5}
\text{and}\quad |\E C_{ii}(z)|\leq \|\bbv_i\|^2\left[\frac{1}{a}\pr(\lambda_{\min}^{\bbS}> a)+\frac{1}{|z|}\pr(\lambda_{\min}^{\bbS}\leq a)\right]\quad \text{is bounded}.
\end{equation}

We then calculate
\begin{eqnarray*}
C_{ii}(z)-\E C_{ii}(z)&=&\bbv_i^T\bbD^{-1}(z)\bbv_i-\E\bbv_i^T\bbD^{-1}(z)\bbv_i
=\tr\bbD^{-1}(z)\bbv_i\bbv_i^T-\E\tr\bbD^{-1}(z)\bbv_i\bbv_i^T\non
&=&\sum_{j=1}^n\Big[\tr\E_j[\bbD^{-1}(z)\bbv_i\bbv_i^T-\bbD_j^{-1}(z)\bbv_i\bbv_i^T]
-\tr\E_{j-1}[\bbD^{-1}(z)\bbv_i\bbv_i^T-\bbD_j^{-1}(z)\bbv_i\bbv_i^T]\Big]\non
&=&-\sum_{j=1}^n(\E_j-\E_{j-1})\beta_j(z)\bbr_j^T\bbD_j^{-1}(z)\bbv_i\bbv_i^T\bbD_j^{-1}(z)\bbr_j.
\end{eqnarray*}
Therefore,
\begin{eqnarray*}
&&\E|C_{ii}(z)-\E C_{ii}(z)|^2=\E\Big|\sum_{j=1}^n(\E_j-\E_{j-1})\beta_j(z)\bbr_j^T\bbD_j^{-1}(z)\bbv_i\bbv_i^T\bbD_j^{-1}(z)\bbr_j\Big|^2\non
&=&\sum_{j=1}^n\E\Big|(\E_j-\E_{j-1})\beta_j(z)\bbr_j^T\bbD_j^{-1}(z)\bbv_i\bbv_i^T\bbD_j^{-1}(z)\bbr_j\Big|^2\non
&\leq&2\sum_{j=1}^n\Big[\E\E_j|\beta_j(z)\bbr_j^T\bbD_j^{-1}(z)\bbv_i\bbv_i^T\bbD_j^{-1}(z)\bbr_j|^2
+\E\E_{j-1}|\beta_j(z)\bbr_j^T\bbD_j^{-1}(z)\bbv_i\bbv_i^T\bbD_j^{-1}(z)\bbr_j|^2\Big]\non
&\leq&4\sum_{j=1}^n\E|\bbr_j^T\bbD_j^{-1}(z)\bbv_i\bbv_i^T\bbD_j^{-1}(z)\bbr_j|^2.
\end{eqnarray*}
What's more, by Lemma 2.2 in \cite{Bai04}, we know that
\begin{eqnarray}\label{add9}
&&\E|\bbr_j^T\bbD_j^{-1}(z)\bbv_i\bbv_i^T\bbD_j^{-1}(z)\bbr_j-\tr\frac{1}{n-1}\bbD_j^{-1}(z)\bbv_i\bbv_i^T\bbD_j^{-1}(z)|^2\non
&=&\E\E_{j-1}|\bbr_j^T\bbD_j^{-1}(z)\bbv_i\bbv_i^T\bbD_j^{-1}(z)\bbr_j-\tr\frac{1}{n-1}\bbD_j^{-1}(z)\bbv_i\bbv_i^T\bbD_j^{-1}(z)|^2\non
&\leq& K\cdot \frac{1}{(n-1)^2}\E\tr(\bbD_j^{-1}(z)\bbv_i\bbv_i^T\bbD_j^{-1}(z))^2=
K\cdot \frac{1}{(n-1)^2}\E(\bbv_i^T\bbD_j^{-2}(z)\bbv_i\bbv_i^T\bbD_j^{-2}(z)\bbv_i)\non
&=&O(\frac{\|\bbv_i\|^4}{n^2}),
\end{eqnarray}
where $K>0$ is a constant and the last step is due to \eqref{add4}. Above two inequalities together show that
\begin{eqnarray}\label{add6}
&&\E|C_{ii}(z)-\E C_{ii}(z)|^2\leq 8\sum_{j=1}^n\Big[\E(\tr\frac{1}{n-1}\bbD_j^{-1}(z)\bbv_i\bbv_i^T\bbD_j^{-1}(z))^2+O(\frac{\|\bbv_i\|^4}{n^2})\Big]\non
&=&8\sum_{j=1}^n\Big[\frac{1}{(n-1)^2}\E(\bbv_i^T\bbD_j^{-2}(z)\bbv_i\bbv_i^T\bbD_j^{-2}(z)\bbv_i)+O(\frac{\|\bbv_i\|^4}{n^2})\Big]=O(\frac{\|\bbv_i\|^4}{n}).
\end{eqnarray}

Combing \eqref{add3}, \eqref{add5} and \eqref{add6}, we get that
\[
\E|C_{ii}^2(z)-(\E C_{ii}(z))^2|=O(\frac{\|\bbv_i\|^4}{\sqrt{n}})
\]
and thus
$$\frac{1}{p}\sum\limits_{i=1}^p C_{ii}^2(z)-\frac{1}{p}\sum\limits_{i=1}^p(\E C_{ii}(z))^2\xrightarrow{i.p} 0.$$

\begin{itemize}
\item {\bf Step 4}. Limit of $\E C_{ii}(z)$, $i=1,\cdots,p$.
\end{itemize}
Note that
\begin{eqnarray}\label{add7}
&&\bbD(z)=\bbS+z\bbI=\sum\limits_{j=1}^n \bbr_j\bbr_j^T+z\bbI \Longrightarrow
\bbv_i^T\bbv_i=\sum\limits_{j=1}^n \bbv_i^T\bbr_j\bbr_j^T\bbD^{-1}(z)\bbv_i+z\bbv_i^T\bbD^{-1}(z)\bbv_i^T\non
&\Longrightarrow&\bbv_i^T\bbv_i=\sum\limits_{j=1}^n \frac{\bbv_i^T\bbr_j\bbr_j^T\bbD_j^{-1}(z)\bbv_i}{1+\bbr_j^T\bbD_j^{-1}(z)\bbr_j}+z\bbv_i^T\bbD^{-1}(z)\bbv_i^T\non
&\Longrightarrow&\bbv_i^T\bbv_i=\sum\limits_{j=1}^n \E\frac{\bbv_i^T\bbr_j\bbr_j^T\bbD_j^{-1}(z)\bbv_i}{1+\bbr_j^T\bbD_j^{-1}(z)\bbr_j}+z\E\bbv_i^T\bbD^{-1}(z)\bbv_i^T.
\end{eqnarray}

We claim that for each $j=1,\cdots,n$,
\begin{equation}\label{add8}
\E\frac{\bbv_i^T\bbr_j\bbr_j^T\bbD_j^{-1}(z)\bbv_i}{1+\bbr_j^T\bbD_j^{-1}(z)\bbr_j}
=\frac{1}{n-1}\cdot\frac{\E\bbv_i^T\bbD_j^{-1}(z)\bbv_i}{1+(n-1)^{-1}\E\tr\bbD_j^{-1}(z)}+o(\frac{1}{n}).
\end{equation}
If this is true, according to \eqref{add7}, one can see that
\begin{eqnarray*}
\bbv_i^T\bbv_i&=&\sum\limits_{j=1}^n\frac{1}{n-1}\cdot\frac{\E\bbv_i^T\bbD_j^{-1}(z)\bbv_i}{1+(n-1)^{-1}\E\tr\bbD_j^{-1}(z)}
+z\E\bbv_i^T\bbD^{-1}(z)\bbv_i^T+o(1)\non
&=&\frac{n}{n-1}\cdot\frac{\E\bbv_i^T\bbD_1^{-1}(z)\bbv_i}{1+(n-1)^{-1}\E\tr\bbD_1^{-1}(z)}
+z\E\bbv_i^T\bbD^{-1}(z)\bbv_i^T+o(1)
\end{eqnarray*}
\begin{eqnarray*}
\Longrightarrow\E\bbv_i^T\bbD_1^{-1}(z)\bbv_i
&=&\frac{n-1}{n}(1+\frac{1}{n-1}\E\tr\bbD_1^{-1}(z))(\bbv_i^T\bbv_i-z\E\bbv_i^T\bbD^{-1}(z)\bbv_i^T+o(1))\non
&\xrightarrow[p/n\rightarrow c]{z\rightarrow 0}&\left(1+c\cdot\frac{1}{1-c}\right)\bbv_i^T\bbv_i=\frac{1}{1-c}\bbv_i^T\bbv_i,
\end{eqnarray*}
where the value in the penultimate step may refer to (3.3.5) in \cite{BS10}.
Note that the only difference between $D_1(z)$ and $D(z)$ is that the sample size is increased from $(n-1)$ to $n$, which does not influence the value $c$. Therefore, we get
\[
\E C_{ii}(z)=\E\bbv_i^T\bbD^{-1}(z)\bbv_i\rightarrow \frac{1}{1-c}\bbv_i^T\bbv_i,\quad i=1,\cdots,p.
\]

The thing left in this step is to verify the claim \eqref{add8}. To this end, we first calculate
\begin{eqnarray*}
&&\E\left|\frac{\bbv_i^T\bbr_j\bbr_j^T\bbD_j^{-1}(z)\bbv_i}{1+\bbr_j^T\bbD_j^{-1}(z)\bbr_j}
-\frac{\bbv_i^T\bbr_j\bbr_j^T\bbD_j^{-1}(z)\bbv_i}{1+(n-1)^{-1}\tr\bbD_j^{-1}(z)}\right|\non
&=&\E\left|\frac{\bbv_i^T\bbr_j\bbr_j^T\bbD_j^{-1}(z)\bbv_i\cdot(\bbr_j^T\bbD_j^{-1}(z)\bbr_j-(n-1)^{-1}\tr\bbD_j^{-1}(z))}
{(1+\bbr_j^T\bbD_j^{-1}(z)\bbr_j)(1+(n-1)^{-1}\tr\bbD_j^{-1}(z))}\right|\non
&\leq&\Big[\E|\bbr_j^T\bbD_j^{-1}(z)\bbv_i\bbv_i^T\bbr_j|^2\cdot\E|\bbr_j^T\bbD_j^{-1}(z)\bbr_j-(n-1)^{-1}\tr\bbD_j^{-1}(z)|^2\Big]^{1/2}.
\end{eqnarray*}
Refer to \eqref{add9}, one can see that
\begin{eqnarray*}
&&\E|\bbr_j^T\bbD_j^{-1}(z)\bbr_j-\tr\frac{1}{n-1}\bbD_j^{-1}(z)|^2
=\E\E_{j-1}|\bbr_j^T\bbD_j^{-1}(z)\bbr_j-\tr\frac{1}{n-1}\bbD_j^{-1}(z)|^2\non
&\leq& K\cdot \frac{1}{(n-1)^2}\E\tr(\bbD_j^{-2}(z))=O(\frac{1}{n}).
\end{eqnarray*}
Similarly,
\begin{eqnarray*}
\E|\bbr_j^T\bbD_j^{-1}(z)\bbv_i\bbv_i^T\bbr_j|^2
\leq 2\Big[\E|\frac{1}{n-1}\bbv_i^T\bbD_j^{-1}(z)\bbv_i|^2+
\E|\bbr_j^T\bbD_j^{-1}(z)\bbv_i\bbv_i^T\bbr_j-\frac{1}{n-1}\bbv_i^T\bbD_j^{-1}(z)\bbv_i|^2\Big]=O(\frac{1}{n^2}).
\end{eqnarray*}
Therefore,
\begin{equation}\label{add10}
\E\left|\frac{\bbv_i^T\bbr_j\bbr_j^T\bbD_j^{-1}(z)\bbv_i}{1+\bbr_j^T\bbD_j^{-1}(z)\bbr_j}
-\frac{\bbv_i^T\bbr_j\bbr_j^T\bbD_j^{-1}(z)\bbv_i}{1+(n-1)^{-1}\tr\bbD_j^{-1}(z)}\right|=O(\frac{1}{n}\sqrt{\frac{1}{n}}).
\end{equation}
It is easy to see that
\begin{equation}\label{add11}
\E\frac{\bbv_i^T\bbr_j\bbr_j^T\bbD_j^{-1}(z)\bbv_i}{1+(n-1)^{-1}\tr\bbD_j^{-1}(z)}=\E\E_{j-1}\frac{\bbv_i^T\bbr_j\bbr_j^T\bbD_j^{-1}(z)\bbv_i}{1+(n-1)^{-1}\tr\bbD_j^{-1}(z)}
=\E\frac{(n-1)^{-1}\bbv_i^T\bbD_j^{-1}(z)\bbv_i}{1+(n-1)^{-1}\tr\bbD_j^{-1}(z)}.
\end{equation}
According to \eqref{add4}, we know that
\begin{eqnarray}\label{add12}
&&\E\left|\frac{(n-1)^{-1}\bbv_i^T\bbD_j^{-1}(z)\bbv_i}{1+(n-1)^{-1}\tr\bbD_j^{-1}(z)}
-\frac{(n-1)^{-1}\bbv_i^T\bbD_j^{-1}(z)\bbv_i}{1+(n-1)^{-1}\E\tr\bbD_j^{-1}(z)}\right|\non
&\leq&
\Big[\E|\frac{1}{n-1}\bbv_i^T\bbD_j^{-1}(z)\bbv_i|^2\cdot \E|\frac{1}{n-1}\tr\bbD_j^{-1}(z)-\frac{1}{n-1}\tr\E\bbD_j^{-1}(z)|^2\Big]^{1/2}\non
&=&\Big[O(\frac{1}{n^2})\cdot \frac{1}{(n-1)^2}\E|\tr\bbD_j^{-1}(z)-\tr\E\bbD_j^{-1}(z)|^2\Big]^{1/2}
=O(\frac{1}{n^2}),
\end{eqnarray}
where in the last step the conclusion $\E|\tr\bbD_j^{-1}(z)-\tr\E\bbD_j^{-1}(z)|^2=O(1)$ refers to the proof of the convergence in distribution of the random part $M_n^{1}(z)$ in \cite{Bai04}. Claim \eqref{add8} is then a direct conclusion from \eqref{add10}-\eqref{add12}.

Above four steps conclude Proposition \ref{prop3} if $\s$ is replaced by $\bbS$. To complete the proof, we further need to verify that such replacement does not influence the result.

\begin{itemize}
\item {\bf Step 5}. Prove $\frac{1}{p}\sum\limits_{i=1}^p \C_{ii}^2(z)-\frac{1}{p}\sum\limits_{i=1}^p C_{ii}^2(z)\xrightarrow{i.p} 0$.
\end{itemize}
 As in above steps, let $\C_{ii}(z)=\bbv_i^T\D^{-1}(z)\bbv_i$, where $\D(z)=\s+z\bbI$. Similar to Step 1 and Step 2, we can also work on $\C_{ii}(z)$ with truncated variables instead of $\C_{ii}$. The details are omitted here. Note that $\D(z)-\bbD(z)=\s-\bbS=-\frac{n}{n-1}\bar\bbx\bar\bbx^T$ is of rank 1. So we can calculate
\begin{equation}\label{add14}
\D^{-1}(z)-\bbD^{-1}(z)=\frac{\frac{n}{n-1}\bbD^{-1}(z)\bar\bbx\bar\bbx^T\bbD^{-1}(z)}
{1-\frac{n}{n-1}\tr(\bbD^{-1}(z)\bar\bbx\bar\bbx^T)}=\frac{\frac{n}{n-1}\bbD^{-1}(z)\bar\bbx\bar\bbx^T\bbD^{-1}(z)}
{1-\frac{n}{n-1}\bar\bbx^T\bbD^{-1}(z)\bar\bbx}
\end{equation}
and
\begin{equation}\label{add13}
\C_{ii}(z)-C_{ii}(z)=\bbv_i^T (\D^{-1}(z)-\bbD^{-1}(z))\bbv_i
=\frac{\frac{n}{n-1}\bar\bbx^T\bbD^{-1}(z)\bbv_i\bbv_i^T\bbD^{-1}(z)\bar\bbx}
{1-\frac{n}{n-1}\bar\bbx^T\bbD^{-1}(z)\bar\bbx}.
\end{equation}
Then we have
\begin{eqnarray}\label{add16}
&&\Big|\frac{1}{p}\sum\limits_{i=1}^p \C_{ii}^2(z)-\frac{1}{p}\sum\limits_{i=1}^p C_{ii}^2(z)\Big|\leq\frac{1}{p}\sum\limits_{i=1}^p |\C_{ii}(z)+C_{ii}(z)||\C_{ii}(z)-C_{ii}(z)|\non
&\leq&\frac{1}{p}\sum\limits_{i=1}^p \|\bbv_i\|^2(\|\D^{-1}(z)\|+\|\bbD^{-1}(z)\|)
\left|\frac{1}{1-\frac{n}{n-1}\bar\bbx^T\bbD^{-1}(z)\bar\bbx}\right|\cdot \frac{n}{n-1}\bar\bbx^T\bbD^{-1}(z)\bbv_i\bbv_i^T\bbD^{-1}(z)\bar\bbx\non
&=&(\|\D^{-1}(z)\|+\|\bbD^{-1}(z)\|)
\left|\frac{1}{1-\frac{n}{n-1}\bar\bbx^T\bbD^{-1}(z)\bar\bbx}\right| \frac{n}{n-1}\cdot\frac{1}{p}\sum\limits_{i=1}^p \|\bbv_i\|^2\bar\bbx^T\bbD^{-1}(z)\bbv_i\bbv_i^T\bbD^{-1}(z)\bar\bbx\non
&&
\end{eqnarray}
By \eqref{add14}, we get
\begin{eqnarray*}
&&\frac{n}{n-1}\bar\bbx^T(\D^{-1}(z)-\bbD^{-1}(z))\bar\bbx
=\frac{(\frac{n}{n-1}\bar\bbx^T\bbD^{-1}(z)\bar\bbx)^2}
{1-\frac{n}{n-1}\bar\bbx^T\bbD^{-1}(z)\bar\bbx}\non
&\Longrightarrow&
1+\frac{n}{n-1}\bar\bbx^T\D^{-1}(z)\bar\bbx=\frac{1}{1-\frac{n}{n-1}\bar\bbx^T\bbD^{-1}(z)\bar\bbx}.
\end{eqnarray*}
Notice that with probability 1, $\|\bar\bbx\|^2$, $\|\D^{-1}(z)\|$ and $\|\bbD^{-1}(z)\|$ are
bounded and thus $[1-\frac{n}{n-1}\bar\bbx^T\bbD^{-1}(z)\bar\bbx]^{-1}=O_p(1)$. What's more, $\|\bbv_i\|^2<\infty$, then according to \eqref{add16}, to prove this step's conclusion, it suffices to show that for each $i$, $\E\bar\bbx^T\bbD^{-1}(z)\bbv_i\bbv_i^T\bbD^{-1}(z)\bar\bbx=o(1)$.
Let $\bar\bbx_j=\frac{1}{n}\sum\limits_{i\neq j}\bbx_i$ and $b_1(z)=[1+(1/(n-1))\E\tr \bbD^{-1}(z)]^{-1}$. We write
\begin{eqnarray*}
&&\E\bar\bbx^T\bbD^{-1}(z)\bbv_i\bbv_i^T\bbD^{-1}(z)\bar\bbx\non
&=&\frac{1}{n}\sum_{j=1}^n\E\bbx_j^T\bbD^{-1}(z)\bbv_i\bbv_i^T\bbD^{-1}(z)\bar\bbx\beta_j(z)
-\frac{1}{n^2}\sum_{j=1}^n\E\bbx_j^T\bbD^{-1}(z)\bbv_i\bbv_i^T\bbD^{-1}(z)\bbx_j\bbx_j^T\bbD^{-1}(z)
\bar\bbx\beta_j^2(z)\non
&=&q_{n1}+q_{n2}+q_{n3}+q_{n4},
\end{eqnarray*}
where
\begin{eqnarray*}
&&q_{n1}=\frac{1}{n}\sum_{j=1}^n\E\bbx_j^T\bbD^{-1}(z)\bbv_i\bbv_i^T\bbD^{-1}(z)\bar\bbx\beta_j(z),
\quad q_{n2}=\frac{1}{n^2}\sum_{j=1}^n\E\bbx_j^T\bbD^{-1}(z)\bbv_i\bbv_i^T\bbD^{-1}(z)\bbx_j\beta_j(z),\non
&&q_{n3}=-\frac{1}{n^2}\sum_{j=1}^n\E\bbx_j^T\bbD^{-1}(z)\bbv_i\bbv_i^T\bbD^{-1}(z)\bbx_j\bbx_j^T\bbD^{-1}(z)
\bar\bbx_j\beta_j^2(z),\non
&&q_{n4}=-\frac{1}{n^3}\sum_{j=1}^n\E\bbx_j^T\bbD^{-1}(z)\bbv_i\bbv_i^T\bbD^{-1}(z)\bbx_j\bbx_j^T\bbD^{-1}(z)
\bbx_j\beta_j^2(z).
\end{eqnarray*}
Similar to Section 2.4 of \cite{pan14}, we can get that $q_{n1}=o(1)$, $q_{n3}=o(1)$ and
\begin{eqnarray*}
&&q_{n2}=b_1(z)\E\Big[\frac{1}{n}\tr\bbD^{-1}(z)\bbv_i\bbv_i^T\bbD^{-1}(z)\Big]+O(\frac{1}{\sqrt{n}}),
\non
&&
q_{n4}=-b_1^2(z)\E\Big[\frac{1}{n}\tr\bbD^{-1}(z)\Big]
\E\Big[\frac{1}{n}\tr\bbD^{-1}(z)\bbv_i\bbv_i^T\bbD^{-1}(z)\Big]+O(\frac{1}{\sqrt{n}}).
\end{eqnarray*}
Note that $|b_1(z)|\leq 1$,
$
\E\Big[\frac{1}{n}\tr\bbD^{-1}(z)\bbv_i\bbv_i^T\bbD^{-1}(z)\Big]
=\frac{1}{n}\E\bbv_i^T\bbD^{-1}(z)\bbv_i=O(\frac{1}{n})$ and
$\E\Big[\frac{1}{n}\tr\bbD^{-1}(z)\Big]= O(1)$. Then $q_{n2}=o(1)$, $q_{n4}=o(1)$ and thus
\[
\E\bar\bbx^T\bbD^{-1}(z)\bbv_i\bbv_i^T\bbD^{-1}(z)\bar\bbx=o(1).
\]
Therefore, $\frac{1}{p}\sum\limits_{i=1}^p \C_{ii}^2(z)-\frac{1}{p}\sum\limits_{i=1}^p C_{ii}^2(z)\xrightarrow{i.p} 0$. Combing the above five steps, we complete the proof of Proposition \ref{prop3}.

\subsection{Proof of Theorem \ref{them1}}

For class 1, recall the notation in (\ref{1805.2}) and
denote
\begin{equation}\label{1030.1}
A=(S_1^0)^{-1}=\left[\frac{1}{n_1-1}(\bbX^0-\bar\bbX^0)(\bbX^0-\bar\bbX^0)^T
\right]^{-1}=(a_{ij})_{p\times p}.
\end{equation}
Similarly, for class 2,
let
\begin{equation}\label{1201.1}
B=(S_2^0)^{-1}=\left[\frac{1}{n_2-1}(\bbY^0-\bar\bbY^0)(\bbY^0-\bar\bbY^0)^T
\right]^{-1}=(b_{ij})_{p\times p}.
\end{equation}
Then
\[
S_1^{-1}=\Sigma_1^{-\frac{1}{2}}A\Sigma_1^{-\frac{1}{2}},\qquad
S_2^{-1}=\Sigma_2^{-\frac{1}{2}}B\Sigma_2^{-\frac{1}{2}}.
\]
When $z$ belongs to class 1, $\bbz=\Sigma_1^{\frac{1}{2}}\bbz^0+\bbmu_1$,
\begin{eqnarray*}
D_1(\bbz)&=&(\bbz-\bar\bbx)^TS_1^{-1}(\bbz-\bar\bbx)
=(\bbz^0-\bar\bbx^0)^T\Sigma_1^{\frac{1}{2}}S_1^{-1}\Sigma_1^{\frac{1}{2}}(\bbz^0-\bar\bbx^0)
=(\bbz^0-\bar\bbx^0)^TA(\bbz^0-\bar\bbx^0)\non
&\triangleq& D_{11}-2D_{12}+D_{13},
\end{eqnarray*}
where
\[
D_{11}=(\bbz^0)^TA\bbz^0, \quad D_{12}=(\bar\bbx^0)^TA(\bbz^0),\quad D_{13}=(\bar\bbx^0)^TA\bar\bbx^0.
\]
By Lemma \ref{lema1}, when $A$ is fixed,
\begin{equation}\label{yq10.1}
\frac{1}{\sqrt{p}}(D_{11}-\tr A)\xrightarrow{D} N(0,\sigma^2),
\end{equation}
where $\sigma^2=(m_4-1)w+2(\theta-w)$ and
\begin{equation}\label{yq10.2}
w=\lim_{p\rightarrow \infty}\frac{1}{p}\sum_{i=1}^p a_{ii}^2,\quad \theta=\lim_{p\rightarrow \infty}\frac{1}{p}\tr A^2.
\end{equation}
Choosing $\bbV=\bbI_p$, $\s=S_1^0$ in Proposition \ref{prop3} gives $\frac{1}{p}\sum\limits_{i=1}^p(\bbv_i^T\bbv_i)^2=1$  and therefore $w=\frac{1}{(1-c_1)^2}$.

Let $s_0$ and $m_0$ be the values of the Stieltjes transforms of $S_1^0$ and $S_2^0$ at the point zero, respectively. And $s'_0$ and $m'_0$ are the corresponding first derivatives at the point zero. Since the LSD of $S_1^0$ and $S_2^0$ tend to the standard MP law (\cite{pan14}) and
the Stieltjes transform $s(z)$ of the standard MP law satisfies the equation
\[
c_izs^2(z)-(1-c_i-z)s(z)+1=0,\quad i=1,2,
\]
it holds then
\begin{equation}\label{yq1004.1}
s_0=\frac{1}{1-c_1},\quad s'_0=\frac{1}{(1-c_1)^3},\quad m_0=\frac{1}{1-c_2},\quad m'_0=\frac{1}{(1-c_2)^3}.
\end{equation}
Let $s_{0n}=\frac{1}{1-p/n_1}$ and $m_{0n}=\frac{1}{1-p/n_2}$.
The central limit theorem (CLT) of linear spectral statistics (LSS) for sample covariance matrices (\cite{Bai04} and \cite{pan14}) implies that
\begin{equation}\label{yq10.6}
\frac{1}{\sqrt{p}}\tr A-\sqrt{p}s_{0n}=O_p(\frac{1}{\sqrt{p}}),\qquad
\frac{1}{p}\tr A^2=s'_0+o_p(1).
\end{equation}
Combing \eqref{yq10.1}-\eqref{yq10.6}, we get that
\begin{equation}\label{1030.3}
\frac{1}{\sqrt{p}}D_{11}-\sqrt{p}s_{0n}\xrightarrow{D} N\Big(0,(m_4-3)s_0^2+2s'_0\Big).
\end{equation}

Next we look at the terms $D_{12}$ and $D_{13}$.
For $D_{12}$, given $\bar\bbx^0$ and $A$, we have
\begin{eqnarray*}
&&\E(\frac{1}{\sqrt{p}}D_{12})=0,\non
&&\E(\frac{1}{p}D_{12}^2)=\frac{1}{p}(\bar\bbx^0)^TA^2\bar\bbx^0\leq \frac{1}{p}\lambda_1(A^2)(\bar\bbx^0)^T\bar\bbx^0.
\end{eqnarray*}
Let $\epsilon>0$ be a sufficiently small constant. With probability 1, $\lambda_1(A^2)\leq (1-\sqrt{c_1}-\epsilon)^{-4}$ and thus
$\frac{1}{p}\lambda_1(A^2)(\bar\bbx^0)^T\bar\bbx^0\leq \frac{1}{p}(1-\sqrt{c_1}-\epsilon)^{-4}(\bar\bbx^0)^T\bar\bbx^0$. Moreover, since
$\E|(\bar\bbx^0)^T\bar\bbx^0|=\E(\bar\bbx^0)^T\bar\bbx^0=\frac{p}{n_1}=c_1$, then $(\bar\bbx^0)^T\bar\bbx^0=O_p(1)$ and $\frac{1}{p}\lambda_1(A^2)(\bar\bbx^0)^T\bar\bbx^0=O_p(\frac{1}{p})$. Therefore, we have
\begin{equation}\label{1030.4}
\frac{1}{\sqrt{p}}D_{12}=O_p(\frac{1}{\sqrt{p}}).
\end{equation}
For $D_{13}$, note that
\[
\frac{1}{\sqrt{p}}D_{13}=\frac{1}{\sqrt{p}}(\bar\bbx^0)^TA\bar\bbx^0
\leq \frac{1}{\sqrt{p}}\lambda_1(A)(\bar\bbx^0)^T\bar\bbx^0.
\]
Since $\lambda_1(A)\leq (1-\sqrt{c_1}-\epsilon)^{-2}$ with probability 1 and $(\bar\bbx^0)^T\bar\bbx^0=O_p(1)$, then
\begin{equation}\label{1030.5}
\frac{1}{\sqrt{p}}D_{13}=O_p(\frac{1}{\sqrt{p}}).
\end{equation}

Combing (\ref{1030.3}), (\ref{1030.4}) and (\ref{1030.5}), we can conclude that
\begin{equation}\label{1030.6}
\frac{1}{\sqrt{p}}D_1(\bbz)-\sqrt{p}s_{0n}\xrightarrow{D} N\Big(0,(m_4-3)s_0^2+2s'_0\Big).
\end{equation}

The weak convergence (\ref{1031.2}) when $\bbz$ belongs to class 2 can be proved in a similar way.

\subsection{Proof of Lemma \ref{lema2}}\label{proflem}
 Denote $\Gamma=\Sigma_2^{-\frac{1}{2}}\Sigma_1^{\frac{1}{2}}=(\gamma_1,\cdots,\gamma_p)$. When $\bbz$ belongs to class 1, we rewrite the difference of two rescaled quadratic terms as
\begin{eqnarray*}
&&\frac{1}{s_{0n}}D_1(\bbz)-\frac{1}{m_{0n}}D_2(\bbz)\non
&=&\frac{1}{s_{0n}}(\bbz-\bar\bbx)^TS_1^{-1}(\bbz-\bar\bbx)-\frac{1}{m_{0n}}(\bbz-\bar\bby)^TS_2^{-1}(\bbz-\bar\bby)\non
&=&\frac{1}{s_{0n}}(\bbz^0-\bar\bbx^0)^TA(\bbz^0-\bar\bbx^0)-\non
&&\frac{1}{m_{0n}}[(\Sigma_1^{\frac{1}{2}}\bbz^0
-\Sigma_2^{\frac{1}{2}}\bar\bby^0)+(\bbmu_1-\bbmu_2)]^T
\Sigma_2^{-\frac{1}{2}}B\Sigma_2^{-\frac{1}{2}}[(\Sigma_1^{\frac{1}{2}}\bbz^0-\Sigma_2^{\frac{1}{2}}\bar\bby^0)
+(\bbmu_1-\bbmu_2)]\non
&=&\left[(\bbz^0)^T\big(\frac{1}{s_{0n}}A-\frac{1}{m_{0n}}\Gamma^T B\Gamma\big)\bbz^0\right]\non
&&+\left[-\frac{2}{s_{0n}}(\bar\bbx^0)^TA\bbz^0+\frac{1}{s_{0n}}(\bar\bbx^0)^TA\bar\bbx^0
-\frac{2}{m_{0n}}(\bar\bby^0)^TB\Gamma\bbz^0+\frac{1}{m_{0n}}(\bar\bby^0)^TB\bar\bby^0\right]\non
&&+\left[-\frac{2}{m_{0n}}(\Gamma\bbz^0-\bar\bby^0)^TB\Sigma_2^{-\frac{1}{2}}(\bbmu_1-\bbmu_2)\right]+\left[
-\frac{1}{m_{0n}}(\bbmu_1-\bbmu_2)^T\Sigma_2^{-\frac{1}{2}}B\Sigma_2^{-\frac{1}{2}}(\bbmu_1-\bbmu_2)\right]\non
&\triangleq&Q_1+Q_2+Q_3+Q_4.
\end{eqnarray*}

For $Q_1$, let $C=\frac{1}{s_{0n}}A-\frac{1}{m_{0n}}\Gamma^T B\Gamma=(c_{ij})_{p\times p}$. By Lemma \ref{lema1}, when $A$ and $B$ are fixed,
\begin{equation}\label{y1031.3}
\frac{1}{\sqrt{p}}(Q_1-\tr C)\xrightarrow{D} N(0,\psi^2),
\end{equation}
where $\psi^2=(m_4-1)\xi+2(\eta-\xi)$ and
\begin{equation*}
\xi=\lim_{p\rightarrow \infty}\frac{1}{p}\sum_{i=1}^p c_{ii}^2,\quad \eta=\lim_{p\rightarrow \infty}\frac{1}{p}\tr C^2.
\end{equation*}
To find the value $\xi$,  write
\[
c_{ii}=\frac{1}{s_{0n}}a_{ii}-\frac{1}{m_{0n}}\gamma_i^TB\gamma_i, \quad i=1,\cdots,p.
\]
Choosing $\bbV=\bbI_p$, $\s=S_1^0$ in Proposition \ref{prop3} gives
\begin{equation}\label{add20}
\frac{1}{p}\sum_{i=1}^p (\frac{1}{s_{0n}}a_{ii})^2\xrightarrow{i.p} 1,
\end{equation}
and selecting $\bbV=\Gamma$, $\s=S_2^0$ tells that
\begin{equation}\label{add21}
\frac{1}{p}\sum_{i=1}^p(\frac{1}{m_{0n}}\gamma_i^TB\gamma_i)^2- \frac{1}{p}\sum_{i=1}^p(\gamma_i^T\gamma_i)^2\xrightarrow{i.p} 0.
\end{equation}
Then we claim that
\begin{equation}\label{add22}
\frac{1}{p}\sum_{i=1}^p c_{ii}^2-\frac{1}{p}\sum_{i=1}^p(1-\gamma_i^T\gamma_i)^2\xrightarrow{i.p} 0.
\end{equation}
To verify this claim, according to \eqref{add20} and \eqref{add21}, we only need to check the cross term
$$
\frac{1}{p}\sum_{i=1}^p \left[ (\frac{1}{s_{0n}}a_{ii})(\frac{1}{m_{0n}}\gamma_i^TB\gamma_i)-1\cdot \gamma_i^T\gamma_i\right]\xrightarrow{i.p} 0.
$$
This can be proved in a similar way as the steps in proving Proposition \ref{prop3}.  We thus only summarize several key points here. Denote $A(z)=(S_1^0+z\bbI)^{-1}=(a_{ij}(z))_{p\times p}$ and $B(z)=(S_2^0+z\bbI)^{-1}=(b_{ij}(z))_{p\times p}$. Then it suffices to verify
\begin{equation}\label{add23}
\frac{1}{p}\sum_{i=1}^p \left[ (\frac{1}{s_{0n}}a_{ii}(z))(\frac{1}{m_{0n}}\gamma_i^TB(z)\gamma_i)-
\E(\frac{1}{s_{0n}}a_{ii}(z))\E(\frac{1}{m_{0n}}\gamma_i^TB(z)\gamma_i)\right]\xrightarrow{i.p} 0.
\end{equation}
To this end, we calculate
\begin{eqnarray*}
&&\E\left|(\frac{1}{s_{0n}}a_{ii}(z))(\frac{1}{m_{0n}}\gamma_i^TB(z)\gamma_i)-
\E(\frac{1}{s_{0n}}a_{ii}(z))\E(\frac{1}{m_{0n}}\gamma_i^TB(z)\gamma_i)\right|\non
&=&\E\bigg|\Big[(\frac{1}{s_{0n}}a_{ii}(z))-\E(\frac{1}{s_{0n}}a_{ii}(z))\Big](\frac{1}{m_{0n}}\gamma_i^TB(z)\gamma_i)
\non
&&+\E(\frac{1}{s_{0n}}a_{ii}(z))
\Big[(\frac{1}{m_{0n}}\gamma_i^TB(z)\gamma_i)
-\E(\frac{1}{m_{0n}}\gamma_i^TB(z)\gamma_i)\Big]
\bigg|\non
&\leq& K\cdot\left\{\E\Big[(\frac{1}{s_{0n}}a_{ii}(z))-\E(\frac{1}{s_{0n}}a_{ii}(z))\Big]^2
+\E\Big[(\frac{1}{m_{0n}}\gamma_i^TB(z)\gamma_i)
-\E(\frac{1}{m_{0n}}\gamma_i^TB(z)\gamma_i)\Big]^2\right\}^{1/2}=o(1),
\end{eqnarray*}
where $K$ is some constant and the last step is concluded by the same method as in deriving \eqref{add6}. Then claim \eqref{add22} is done and
together with Condition \ref{cond3}, we have
\begin{eqnarray}\label{1118.1}
\xi&=&\lim_{p\rightarrow \infty}\frac{1}{p}\sum_{i=1}^p c_{ii}^2
=\lim_{p\rightarrow \infty}\frac{1}{p}\sum_{i=1}^p (1-\gamma_i^T\gamma_i)^2
=\lim_{p\rightarrow \infty}\left[1-\frac{2}{p}\tr(\Sigma_1\Sigma_2^{-1})+\frac{1}{p}\sum_{i=1}^p[(\Sigma_1^{\frac{1}{2}}\Sigma_2^{-1}\Sigma_1^{\frac{1}{2}})_{ii}]^2\right]\non
&=&1-2M_1+M_2.
\end{eqnarray}

To calculate the value $\eta$, we write
\[
\frac{1}{p}\tr C^2=\frac{1}{p}\tr\Big[\frac{1}{s_{0n}}A-\frac{1}{m_{0n}}\Gamma^T B\Gamma\Big]^2=\frac{1}{s_{0n}^2}\frac{1}{p}\tr A^2-\frac{2}{s_{0n}m_{0n}}\frac{1}{p}\tr(A\Gamma^T B\Gamma)+\frac{1}{m_{0n}^2}\frac{1}{p}\tr(\Gamma^T B\Gamma)^2.
\]
Recall that $\frac{1}{p}\tr A^2\rightarrow \theta=s'_0$
from (\ref{yq10.6}). Next consider the limit for $\frac{1}{p}\tr(A\Gamma^T B\Gamma)$. Denote $\Lambda=(\Gamma^T B\Gamma)^{-1}=\Gamma^{-1}S_2^0(\Gamma^T)^{-1}$. Given $B$, we can view $\Lambda$ as a population covariance matrix and treat $(\Lambda^{\frac{1}{2}}A^{-1}\Lambda^{\frac{1}{2}})$ as the corresponding sample covariance matrix, then the CLT of LSS for sample covariance matrices in \cite{pan14} implies that
\[
\frac{1}{p}\tr(A\Gamma^T B\Gamma)=\frac{1}{p}\tr(\Lambda^{\frac{1}{2}}A^{-1}\Lambda^{\frac{1}{2}})^{-1}
= \frac{s_0}{p}\tr(\Lambda^{-1})+O_p(\frac{1}{p})=\frac{s_0}{p}\tr(\Gamma^T B\Gamma)+O_p(\frac{1}{p}).
\]
Note that $(\Gamma^T B\Gamma)^{-1}$ can again be treated as a sample covariance matrix. Therefore,
\[
\frac{s_0}{p}\tr(\Gamma^T B\Gamma)=\frac{s_0m_0}{p}\tr(\Gamma^T \Gamma)+O_p(\frac{1}{p})=\frac{s_0m_0}{p}\tr(\Sigma_1\Sigma_2^{-1})+O_p(\frac{1}{p}).
\]
Finally for the term $\frac{1}{p}\tr(\Gamma^T B\Gamma)^2$, we denote
\[
s_{\Gamma}(z)=\int\frac{1}{\lambda-z}dF^{\Lambda}(\lambda)
\]
as the Stieltjes transform of the sample covariance matrix
$
\Lambda=(\Gamma^T B\Gamma)^{-1}=(\Sigma_1^{-\frac{1}{2}}\Sigma_2^{\frac{1}{2}})S_2^0(\Sigma_2^{\frac{1}{2}}\Sigma_1^{-\frac{1}{2}})
$.
Then
\[
\frac{1}{p}\tr(\Gamma^T B\Gamma)^2= s'_0(\Gamma)+O_p(\frac{1}{p}),
\]
where $s'_0(\Gamma)$ is the first derivative of $s_{\Gamma}(z)$
at the point zero. By \cite{Bai04} and \cite{pan14}, for each $z\in \mathbb{C}^{+}\equiv\{z\in\mathbb{C}: \Im z>0\}$, the Stieltjes transform $s\equiv s_{\Gamma}(z)$ is the unique solution to
\[
s=\int\frac{1}{\lambda(1-c_2-c_2zs)-z}dF^{[\Sigma_1^{-\frac{1}{2}}\Sigma_2\Sigma_1^{-\frac{1}{2}}]}(\lambda).
\]
We thus can calculate that
\[
s_{\Gamma}(0)=\frac{1}{1-c_2}\cdot\frac{1}{p}\tr(\Sigma_1^{-\frac{1}{2}}\Sigma_2\Sigma_1^{-\frac{1}{2}})^{-1}=\frac{1}{1-c_2}\cdot\frac{1}{p}\tr(\Sigma_1\Sigma_2^{-1}),
\]
\begin{eqnarray*}
s'_0(\Gamma)&=&\frac{1}{(1-c_2)^2}\left[\frac{1}{p}\tr(\Sigma_1^{-\frac{1}{2}}\Sigma_2\Sigma_1^{-\frac{1}{2}})^{-2}+c_2s_{\Gamma}(0)\frac{1}{p}\tr(\Sigma_1^{-\frac{1}{2}}\Sigma_2\Sigma_1^{-\frac{1}{2}})^{-1}\right]\non
&=&m_0^2\left[\frac{1}{p}\tr(\Sigma_1\Sigma_2^{-1})^2+\frac{c_2}{1-c_2}\Big(\frac{1}{p}\tr(\Sigma_1\Sigma_2^{-1})\Big)^2\right].
\end{eqnarray*}

Combining the above three parts, we have
\begin{eqnarray}\label{1118.2}
\eta&=&\lim_{p\rightarrow \infty}\frac{1}{p}\tr C^2\non
&=& \lim_{p\rightarrow \infty}\left\{ \frac{s'_0}{s_0^2}-\frac{2}{p}\tr(\Sigma_1\Sigma_2^{-1})+
\left[\frac{1}{p}\tr(\Sigma_1\Sigma_2^{-1})^2+
\frac{c_2}{1-c_2}\Big(\frac{1}{p}\tr(\Sigma_1\Sigma_2^{-1})\Big)^2\right]\right\}\non
&=&\frac{1}{1-c_1}-2M_1+M_5+\frac{c_2}{1-c_2}M_1^2.
\end{eqnarray}
Therefore,
\begin{eqnarray}\label{y1201.2}
\psi^2&=&(m_4-1)\xi+2(\eta-\xi)\non
&=&(m_4-3)(1-2M_1+M_2)+2\left(\frac{1}{1-c_1}-2M_1+M_5+\frac{c_2}{1-c_2}M_1^2\right).
\end{eqnarray}

We claim that $\psi^2$ is a positive constant. To this end, one may check the following two aspects. First, the definitions of $\xi$ and $\eta$ imply that $\eta\geq\xi\geq 0$, thus $\psi^2\geq 0$. Secondly, the parameter $\eta$ has a positive lower bound:
\begin{eqnarray}\label{1126.7}
\eta&=&\lim_{p\rightarrow \infty}\left[\frac{s'_0}{s_0^2}-\frac{2}{p}\tr(\Sigma_1\Sigma_2^{-1})+\frac{1}{p}\tr(\Sigma_1\Sigma_2^{-1})^2+\frac{c_2}{1-c_2}\Big(\frac{1}{p}\tr(\Sigma_1\Sigma_2^{-1})\Big)^2\right]\non
&=&\frac{1}{1-c_1}+\lim_{p\rightarrow \infty}\left[-\frac{2}{p}\tr(\Sigma_1\Sigma_2^{-1})+\frac{1}{p}\tr(\Sigma_1\Sigma_2^{-1})^2+\frac{c_2}{1-c_2}\Big(\frac{1}{p}\tr(\Sigma_1\Sigma_2^{-1})\Big)^2\right]\non
&=&\frac{c_1}{1-c_1}+\lim_{p\rightarrow \infty}\left[\frac{1}{p}\tr(\bbI_p-\Sigma_1\Sigma_2^{-1})^2+\frac{c_2}{1-c_2}\Big(\frac{1}{p}\tr(\Sigma_1\Sigma_2^{-1})\Big)^2\right]>\frac{c_1}{1-c_1}.
\end{eqnarray}
According to these two points, one can see that $\psi^2>0$.
Moreover, as in (\ref{yq10.6}), we can find that
\begin{eqnarray}\label{1806.1}
\frac{1}{\sqrt{p}}\tr C&=&\frac{1}{s_{0n}}\frac{1}{\sqrt{p}}\tr A-\frac{1}{m_{0n}}\frac{1}{\sqrt{p}}\tr\Gamma^T B\Gamma=\sqrt{p}-\frac{1}{\sqrt{p}}\sum_{i=1}^p\gamma_i^T\gamma_i+O_p(\frac{1}{\sqrt{p}})\non
&=&\sqrt{p}-\frac{1}{\sqrt{p}}\tr\Gamma^T\Gamma+O_p(\frac{1}{\sqrt{p}})
=\frac{1}{\sqrt{p}}\tr(\bbI_p-\Sigma_1\Sigma_2^{-1})+O_p(\frac{1}{\sqrt{p}}).
\end{eqnarray}

Combing (\ref{y1031.3}), (\ref{y1201.2}) and (\ref{1806.1}), we have
\begin{equation}\label{1118.4}
\frac{1}{\sqrt{p}}Q_1-\frac{1}{\sqrt{p}}\tr(\bbI_p-\Sigma_1\Sigma_2^{-1})\xrightarrow{D}N(0,\psi^2),
\end{equation}
where
$\psi^2>0$ is given in (\ref{y1201.2}).

For $Q_2$, by applying similar methods in dealing with the terms $D_{12}$ and $D_{13}$ in (\ref{1030.4}) and (\ref{1030.5}), the following result holds
\begin{equation*}
\frac{1}{\sqrt{p}}Q_2=O_p(\frac{1}{\sqrt{p}}).
\end{equation*}

For $Q_3$,
we write
\[
\frac{1}{\sqrt{p}}Q_3=2(-\frac{1}{\sqrt{p}}Q_{31}+\frac{1}{\sqrt{p}}Q_{32}),
\]
where
\[
\frac{1}{\sqrt{p}}Q_{31}=\frac{1}{\sqrt{p}}\frac{1}{m_{0n}}(\Gamma\bbz^0)^TB\Sigma_2^{-\frac{1}{2}}(\bbmu_1-\bbmu_2),\quad
\frac{1}{\sqrt{p}}Q_{32}=\frac{1}{\sqrt{p}}\frac{1}{m_{0n}}(\bar\bby^0)^TB\Sigma_2^{-\frac{1}{2}}(\bbmu_1-\bbmu_2).
\]
Given $B$, we have $\E(\frac{1}{\sqrt{p}}Q_{31})=0$ and
\[
\E(\frac{1}{p}Q_{31}^2)=\frac{1}{pm_{0n}^2}(\bbmu_1-\bbmu_2)^T\Sigma_2^{-\frac{1}{2}}B\Gamma\Gamma^TB\Sigma_2^{-\frac{1}{2}}(\bbmu_1-\bbmu_2)
\leq \frac{\|\bbmu_1-\bbmu_2\|^2}{pm_{0n}^2}\lambda_1(\Sigma_2^{-\frac{1}{2}}B\Gamma\Gamma^TB\Sigma_2^{-\frac{1}{2}}).
\]
Thus
$
\frac{1}{\sqrt{p}}Q_{31}=O_p(\frac{\sqrt{|T_3|}}{p^{1/4}})
$, where $T_3$ is defined in (\ref{1126.6}) below.
Moreover, note that
\[
\left|\frac{1}{\sqrt{p}}Q_{32}\right|\leq\frac{1}{\sqrt{p}m_{0n}}\|\bar\bby^0\|\|B\|\|\Sigma_2^{-\frac{1}{2}}\|\|\bbmu_1-\bbmu_2\|
\]
and
\[
\E|(\bar\bby^0)^T\bar\bby^0|=\E(\bar\bby^0)^T\bar\bby^0=\frac{p}{n_2}=c_2.
\]
We get that $\|\bar\bby^0\|=O_p(1)$ and
$
\frac{1}{\sqrt{p}}Q_{32}=O_p(\frac{\sqrt{|T_3|}}{p^{1/4}})$. Therefore
\[
\frac{1}{\sqrt{p}}Q_{3}=O_p(\frac{\sqrt{|T_3|}}{p^{1/4}}).
\]

With the properties of $Q_1$, $Q_2$ and $Q_3$, now we are able to analyze the misclassification probability of classifying $\bbz$ to class 2. That is,
\begin{eqnarray*}
\pr_{2|1}^G
&=&\pr\left\{\frac{1}{s_{0n}}D_1(\bbz)+\log |S_1|-p l_{1n}>\frac{1}{m_{0n}}D_2(\bbz)+\log |S_2|-p l_{2n}\right\}\non
&=&\pr\bigg\{\frac{1}{\sqrt{p}}\Big(Q_1+ Q_2+Q_3+Q_4\Big)+\frac{1}{\sqrt{p}}\log|\Sigma_1|
-\frac{1}{\sqrt{p}}\log|\Sigma_2|\non
&&\quad\quad +\frac{1}{\sqrt{p}}\Big(\log|S_1^0|-pl_{1n}\Big)-\frac{1}{\sqrt{p}}\Big(\log|S_2^0|-pl_{2n}\Big)>0\bigg\}\non
&=&\pr\bigg\{T_1+T_2 + T_3+T_4 +\frac{1}{\sqrt{p}}Q_2+\frac{1}{\sqrt{p}}Q_3+\frac{1}{\sqrt{p}}\Big(\log|S_1^0|-pl_{1n}\Big)-\frac{1}{\sqrt{p}}\Big(\log|S_2^0|-pl_{2n}\Big)>0\bigg\},
\end{eqnarray*}
where
\begin{eqnarray}\label{1126.6}
&&T_1=\frac{1}{\sqrt{p}}Q_1-\frac{1}{\sqrt{p}}\tr(\bbI_p-\Sigma_1\Sigma_2^{-1}),\quad
T_2=\frac{1}{\sqrt{p}}\tr(\bbI_p-\Sigma_1\Sigma_2^{-1})+\frac{1}{\sqrt{p}}\log|\Sigma_1\Sigma_2^{-1}|,\non
&&T_3=-\frac{1}{\sqrt{p}}(\bbmu_1-\bbmu_2)^T\Sigma_2^{-1}(\bbmu_1-\bbmu_2),\quad T_4=\frac{1}{\sqrt{p}}Q_4-T_3.
\end{eqnarray}
According to (\ref{1118.4}), one can see that
\begin{equation}\label{1126.1}
T_1\xrightarrow{D}N(0,\psi^2).
\end{equation}
Moreover,
\begin{equation}\label{1126.4}
 |T_4|\prec\frac{1}{p}\|\bbmu_1-\bbmu_2\|^2, \quad T_4=O_p(\frac{|T_3|}{\sqrt{p}}).
\end{equation}
And according to the central limit distributions for $\log|S_1^0|$ and $\log|S_2^0|$ in \cite{pan14},
\begin{equation}\label{1126.5}
\frac{1}{\sqrt{p}}\Big(\log|S_i^0|-pl_{in}\Big)=O_p(\frac{1}{\sqrt{p}}),\quad i=1,2.
\end{equation}
The last point is on the terms $T_2$ and $T_3$. By the same argument for $\mathcal{E}_{2|1}^O$ in Section \ref{interp}, we know that as long as $\Sigma_1\neq \Sigma_2$, $T_2+T_3<0$.
Based on the above observations, the misclassification probability $\pr_{2|1}^G$ can be calculated as
\begin{eqnarray*}
\pr_{2|1}^G&=&\pr\bigg\{T_1+T_2 + T_3+T_4 +\frac{1}{\sqrt{p}}Q_2+\frac{1}{\sqrt{p}}Q_3+
\frac{1}{\sqrt{p}}\Big(\log|S_1^0|-pl_{1n}\Big)-\frac{1}{\sqrt{p}}\Big(\log|S_2^0|-pl_{2n}\Big)>0\bigg\}\non
&=&\pr\bigg\{T_1>-T_2-T_3+O_p(\frac{|T_3|}{\sqrt{p}})+O_p(\frac{\sqrt{|T_3|}}{p^{1/4}})+O_p(\frac{1}{\sqrt{p}})\bigg\}\non
&=&\pr\bigg\{\frac{T_1}{\psi}>\frac{1}{\psi}\big[-T_2-T_3+O_p(\frac{|T_3|}{\sqrt{p}})+O_p(\frac{\sqrt{|T_3|}}{p^{1/4}})+O_p(\frac{1}{\sqrt{p}})\big]\bigg\}
\xrightarrow{i.p} 1-\Phi\left(\frac{T}{\psi}\right),
\end{eqnarray*}
where $\psi$ is calculated in (\ref{y1201.2}), $T$ is defined in (\ref{1220.1}).

\subsection{Proof of Corollary \ref{coronew}}
In case $(i)$, if $\zeta_1\rightarrow \infty$, one can easily see that $-T_3\rightarrow\infty$ and $-\widetilde{T}_3\rightarrow\infty$. If $\zeta(\epsilon)\rightarrow \infty$,
according to the derivation for $\mathcal{E}_{2|1}^O$ in Section \ref{interp}, we get that $-T_2\rightarrow\infty$ and $-\widetilde{T}_2\rightarrow\infty$. Therefore, both $T\rightarrow\infty$ and $\widetilde{T}\rightarrow\infty$ in case $(i)$. Combining with Theorem \ref{themnew}, we have $R^G\xrightarrow{i.p}0$. This completes the proof of the first case. The other two cases can be verified similarly.

\subsection{Proof of Proposition \ref{prop1}}
Proposition \ref{prop1} can be proved in a similar way as Theorem \ref{themnew}. Below we only write down the key steps in the proof.
\begin{eqnarray*}
\pr_{2|1}^O&=&\pr\left\{d_1(\bbz)+\log |\Sigma_1|>d_2(\bbz)+\log |\Sigma_2|\right\}\non
&=&\pr\bigg\{(\bbz^0)^T\bbz^0+\log |\Sigma_1|>(\bbz^0)^T(\Sigma_1^{\frac{1}{2}}\Sigma_2^{-1}\Sigma_1^{\frac{1}{2}})\bbz^0
+2(\bbmu_1-\bbmu_2)^T\Sigma_2^{-1}\Sigma_1^{\frac{1}{2}}\bbz^0+\non
&&\qquad(\bbmu_1-\bbmu_2)^T\Sigma_2^{-1}(\bbmu_1-\bbmu_2)+\log |\Sigma_2|\bigg\}\non
&=&\pr\bigg\{T_1^0+T_2 + T_3+\frac{1}{\sqrt{p}}Q_3^0>0\bigg\},
\end{eqnarray*}
where
\begin{eqnarray}
&&T_1^0=\frac{1}{\sqrt{p}}(\bbz^0)^T(\bbI_p-\Sigma_1^{\frac{1}{2}}\Sigma_2^{-1}\Sigma_1^{\frac{1}{2}})\bbz^0
-\frac{1}{\sqrt{p}}\tr(\bbI_p-\Sigma_1\Sigma_2^{-1}),\quad Q_3^0=-2(\bbmu_1-\bbmu_2)^T\Sigma_2^{-1}\Sigma_1^{\frac{1}{2}}\bbz^0,\non
&&T_2=\frac{1}{\sqrt{p}}\tr(\bbI_p-\Sigma_1\Sigma_2^{-1})+\frac{1}{\sqrt{p}}\log|\Sigma_1\Sigma_2^{-1}|,\quad
T_3=-\frac{1}{\sqrt{p}}(\bbmu_1-\bbmu_2)^T\Sigma_2^{-1}(\bbmu_1-\bbmu_2).
\end{eqnarray}
Note that $T_2$ and $T_3$ here are the same as in the proof of Lemma \ref{lema2}. By applying the same arguments in dealing with $\frac{1}{\sqrt{p}}Q_3$, one can see that
$
\frac{1}{\sqrt{p}}Q_{3}^0=O_p(\frac{\sqrt{|T_3|}}{p^{1/4}}).
$
Following the same steps in deriving $T_1$, we can conclude that
\begin{equation}
T_1^0\xrightarrow{D}N(0,\psi_0^2),
\end{equation}
where $\psi_0^2=(m_4-1)\xi_0+2(\eta_0-\xi_0)$ and
\begin{eqnarray*}
&&\xi_0=\lim_{p\rightarrow \infty}\frac{1}{p}\sum_{i=1}^p (\bbI_p-\Sigma_1^{\frac{1}{2}}\Sigma_2^{-1}\Sigma_1^{\frac{1}{2}})_{ii}^2=1-2M_1+M_2=\xi,\non &&\eta_0=\lim_{p\rightarrow \infty}\frac{1}{p}\tr (\bbI_p-\Sigma_1^{\frac{1}{2}}\Sigma_2^{-1}\Sigma_1^{\frac{1}{2}})^2=1-2M_1+M_5.
\end{eqnarray*}
Therefore,
\begin{eqnarray*}
\pr_{2|1}^O&=&\pr\bigg\{T_1^0+T_2 + T_3+\frac{1}{\sqrt{p}}Q_3^0>0\bigg\}=\pr\bigg\{T_1^0>-T_2-T_3+O_p(\frac{\sqrt{|T_3|}}{p^{1/4}})\bigg\}\non
&=&\pr\bigg\{\frac{T_1^0}{\psi_0}>\frac{1}{\psi_0}\big[-T_2-T_3+O_p(\frac{\sqrt{|T_3|}}{p^{1/4}})\big]\bigg\}
\xrightarrow{i.p} 1-\Phi\left(\frac{T}{\psi_0}\right).
\end{eqnarray*}
The conclusion that $\pr_{1|2}^O\xrightarrow{i.p} 1-\Phi\left(\frac{\widetilde{T}}{\widetilde{\psi}_0}\right)$
can be shown similarly. Therefore we complete the proof of Proposition \ref{prop1}.

\subsection{Proof of Corollary \ref{coro3} and Theorem \ref{them3}}
Corollary \ref{coro3} is a direct implication of Proposition \ref{prop1} and Theorem \ref{them3} is easy to see from Corollary \ref{coronew} and  Corollary \ref{coro3}. We take the case when $\zeta_1\rightarrow \infty$ or $\zeta(\epsilon)\rightarrow  \infty$ as an example. Under this assumption, from the arguments in the proof of case $(i)$ in Corollary \ref{coronew}, we have \{$-T_2\rightarrow\infty$ and $ -\widetilde{T}_2\rightarrow\infty$\} or \{$-T_3\rightarrow\infty$ and $-\widetilde{T}_3\rightarrow\infty$\}.  Theorem \ref{themnew} tells that $R^G\xrightarrow{i.p}0$ and Proposition \ref{prop1} indicates that $R^O\xrightarrow{i.p}0$. Therefore, $\diff\xrightarrow{i.p} 0$. The other cases can be proved similarly and thus we omit them here.

\end{appendices}

\newpage

\setcounter{section}{0}
\setcounter{equation}{0}
\setcounter{figure}{0}
\setcounter{table}{0}
\setcounter{page}{1}
\makeatletter
\renewcommand{\thesection}{S\arabic{section}}
\renewcommand{\thefigure}{S\arabic{figure}}

\begin{center}
{\Large\bf Supplementary Material to ``Quadratic Discriminant Analysis under Moderate Dimension''}
\end{center}

\section{Modification for the ``hard'' case -- Divide-and-conquer over samples (negative result)}\label{negativ}
In this section, we investigate the conventional divide-and-conquer method over samples in the ``hard'' case.
For each class,
the samples are divided into $H$ non-overlapping groups, each group\footnote{We assume equal group size for simplicity, with a straightforward extension to unequal size.} with $m_1=\lfloor\frac{n_1}{H}\rfloor$ and $m_2=\lfloor\frac{n_2}{H}\rfloor$ subsamples
respectively. For $k=1,\cdots, H$, define
\begin{eqnarray}\label{0322.1}
\D_k&=&\left[\frac{1}{s_{0(k)}}(\bbz-\bar\bbx_{(k)})^TS_{1(k)}^{-1}(\bbz-\bar\bbx_{(k)})
+\log|S_{1(k)}|-pl_{1(k)}\right]\non
&-&\left[\frac{1}{m_{0(k)}}(\bbz-\bar\bby_{(k)})^TS_{2(k)}^{-1}(\bbz-\bar\bby_{(k)})
+\log|S_{2(k)}|-pl_{2(k)}\right],
\end{eqnarray}
where the terms with the subscript ``$(k)$'' indicate the corresponding values in the $k$-th group. For example, \{$\bbx_{1(k)},\cdots,\bbx_{m_1(k)}$\} are the $m_1$ samples in the $k$-the group of class 1,
$
\bar\bbx_{(k)}=\frac{1}{m_1}\sum_{i=1}^{m_1}\bbx_{i{(k)}}$,
$S_{1{(k)}}=\frac{1}{m_1-1}\sum_{i=1}^{m_1}(\bbx_{i{(k)}}-\bar\bbx_{(k)})(\bbx_{i{(k)}}-\bar\bbx_{(k)})^T$, $s_{0(k)}=\frac{1}{1-p/m_1}$, and $m_{0(k)}$, $l_{1(k)}$, $l_{2(k)}$ are defined similarly by replacing $n_1$ ($n_2$) with $m_1$ $(m_2)$. With the generalized QDA statistics $\D_k$ from these subsamples, we consider two classification rules -- weighted voting and majority voting.

\vspace{10pt}

{\bf Weighted voting:}

A weighted divide-and-conquer version of the generalized QDA rule (\ref{11.1}) is defined as:
classify $\bbz$ to class 1 if and only if
\begin{equation}\label{0322.2}
\D<0,\qquad \text{where}\qquad \D=\frac{1}{H}\sum_{k=1}^H \D_k.
\end{equation}

\begin{figure}[!htp]
    \centering
    \includegraphics[width=5.5in]{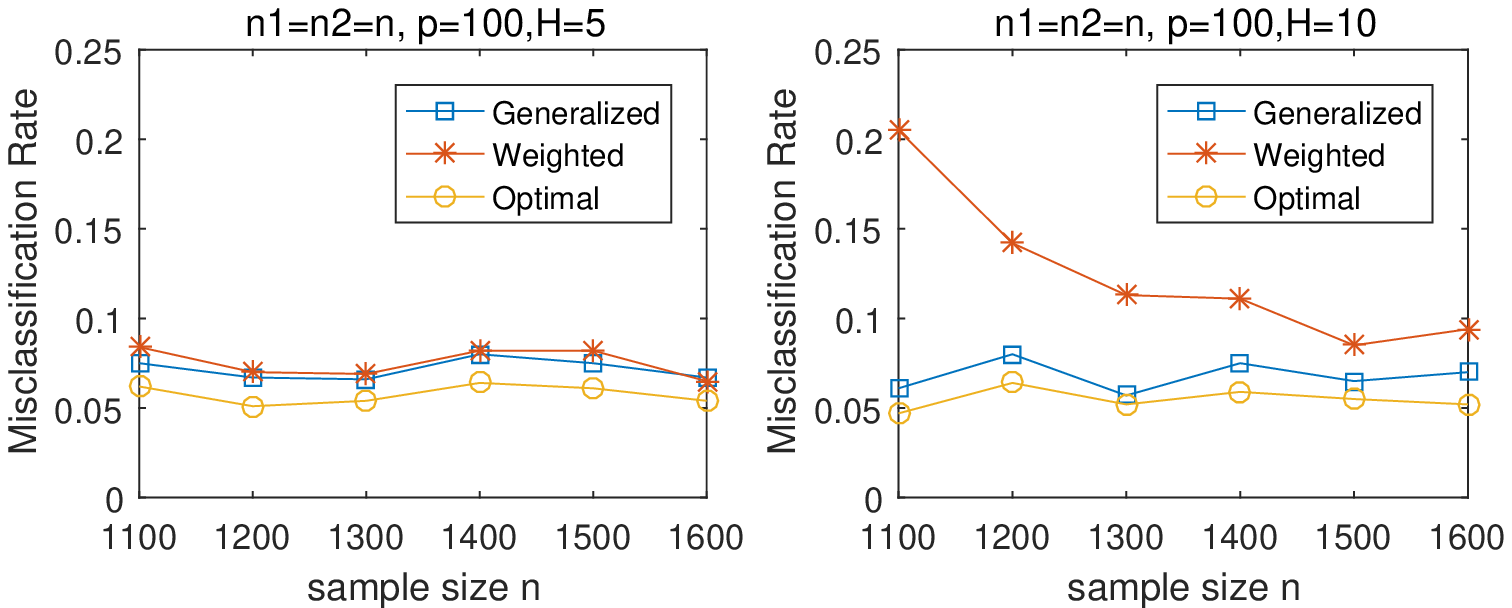}
    \caption{\small{\it Comparison of empirical misclassification rates $R^G$, $R^{\D}$ (Left: $H=5$; Right: $H=10$) and $R^O$ under Case 5 (``hard'' case). $\bbmu_1=\bbmu_2=0$, $n_1=n_2=n$, $p=100$ and 1000 replications.}}
    \label{dc1}
  \end{figure}

Despite of its success in the literature, we claim that divide-and-conquer over samples fails to boost our generalized QDA, both theoretically and empirically. Theoretically,
its misclassification rate $R^{\D}$ is established in Theorem \ref{them4}.
\begin{thm}\label{them4}Assume $p/m_1\rightarrow c_{1(H)}=c_1H\in(0,1)$ and $p/m_2\rightarrow c_{2(H)}=c_2H\in(0,1)$.
Under Conditions \ref{cond1}-\ref{cond3}, the misclassification rate of the modified rule (\ref{0322.2})
\[
R^{\D}=\frac{1}{2}\left[\pr_{2|1}^{\D}+\pr_{1|2}^{\D}\right]\xrightarrow{i.p}
1-\frac{1}{2}\left[\Phi\left(\frac{T}{\psi_{\D}}\right)+\Phi\left(\frac{\widetilde{T}}{\widetilde{\psi}_{\D}}\right)\right],
\]
where $T$ and $\widetilde{T}$ are the same as in Theorem \ref{themnew}.
The parameters $\psi_{\D}$ and $\widetilde{\psi}_{\D}$ are  positive constants given by
\begin{eqnarray*}
\psi_{\D}^2
&=&\psi_0^2+\frac{2}{H}\left[\frac{c_{1(H)}}{1-c_{1(H)}}+\frac{c_{2(H)}}{1-c_{2(H)}}M_1^2\right],\non
\widetilde{\psi}_{\D}^2&=&\widetilde{\psi}_0^2+\frac{2}{H}\left[\frac{c_{2(H)}}{1-c_{2(H)}}+\frac{c_{1(H)}}{1-c_{1(H)}}M_3^2\right].
\end{eqnarray*}
\end{thm}
Comparing the expressions of $\psi_0^2$, $\psi^2$ and $\psi_{\D}^2$ (similar phenomenon in $\widetilde{\psi}_0^2$, $\widetilde{\psi}^2$ and $\widetilde{\psi}_{\D}^2$)
we can calculate that
\begin{eqnarray*}
\psi_{\D}^2-\psi^2&=&\frac{2}{H}\left[\frac{c_{1(H)}}{1-c_{1(H)}}+\frac{c_{2(H)}}{1-c_{2(H)}}M_1^2\right]
-2\left[\frac{c_1}{1-c_1}+\frac{c_2}{1-c_2}M_1^2\right]\non
&=&2
\left[\frac{c_1}{1-c_1H}-\frac{c_1}{1-c_1}+\Big(\frac{c_2}{1-c_2H}-\frac{c_2}{1-c_2}\Big)M_1^2\right]\geq 0.
\end{eqnarray*}
Therefore $R^{\D}\geq R^{G}$ and ``$=$'' holds if and only if when $H=1$, i.e. no divide.
Empirically, Figure \ref{dc1} compares the performance after applying \eqref{0322.2} to the hard case 5 and it shows that this approach works worse than the original generalized QDA. Moreover, this inferiority becomes more obvious for larger $H$.


\vspace{10pt}

{\bf Majority voting:}

\begin{figure}[!htp]
    \centering
    \includegraphics[width=5.5in]{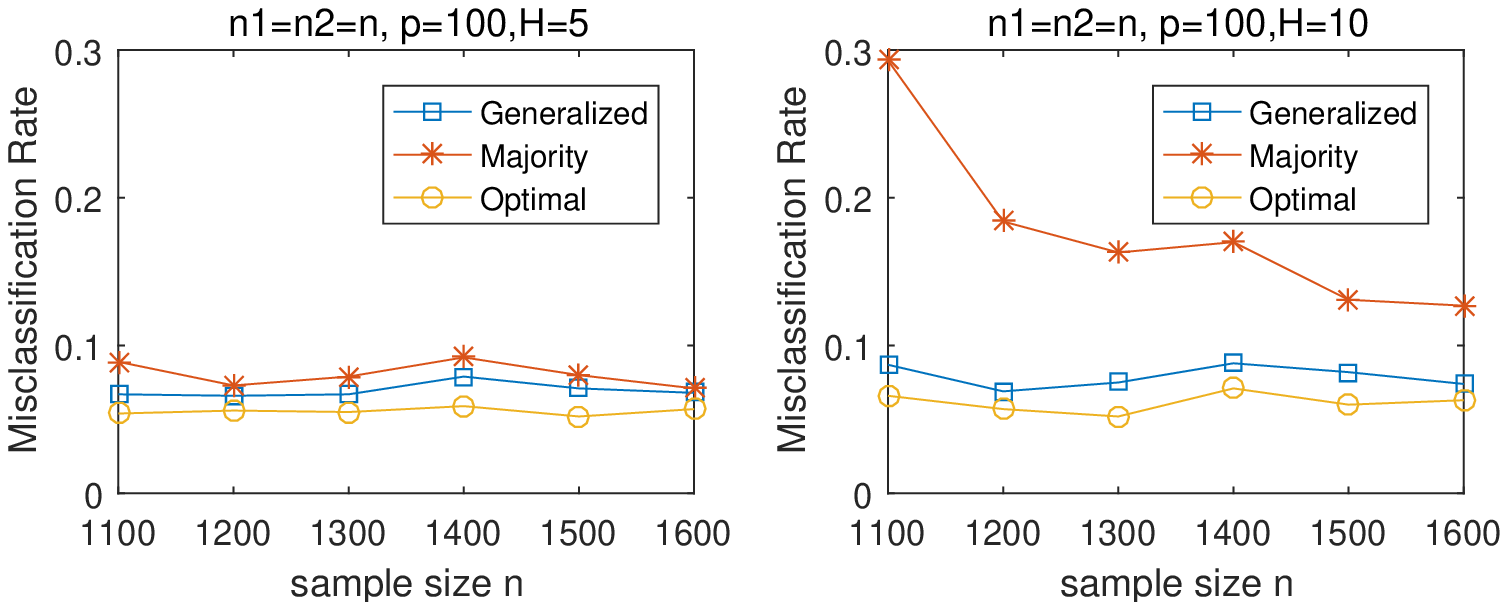}
    \caption{\small{\it Comparison of empirical misclassification rates under Case 5. $\bbmu_1=\bbmu_2=0$, $n_1=n_2=n$, $p=100$ and 1000 replications. ``Majority'' is the rule \eqref{dc3} (Left: $H=5$; Right: $H=10$). }}
    \label{dc2}
  \end{figure}

Another approach is via the well-known majority voting over the $H$ groups, that is, classify $\bbz$ to class 1 if and only if
\begin{equation}\label{dc3}
\#\left\{1\leq k\leq H: \D_k<0\right\}> H/2.
\end{equation}
Its empirical performance is recorded in Figure \ref{dc2}, which points to the same conclusion as the weighted voting \eqref{0322.2}.

\subsection{Proof of Theorem \ref{them4}}
Recall the notation $\Gamma=\Sigma_2^{-\frac{1}{2}}\Sigma_1^{\frac{1}{2}}=(\gamma_1,\cdots,\gamma_p)$ and let
\[
A_{(k)}=\left[\frac{1}{m_1-1}(\bbX_{(k)}^0-\bar\bbX_{(k)}^0)(\bbX_{(k)}^0-\bar\bbX_{(k)}^0)^T
\right]^{-1}=(a_{ij}^{(k)})_{p\times p},
\]
\[
B_{(k)}=\left[\frac{1}{m_2-1}(\bbY_{(k)}^0-\bar\bbY_{(k)}^0)(\bbY_{(k)}^0-\bar\bbY_{(k)}^0)^T
\right]^{-1}=(b_{ij}^{(k)})_{p\times p}.
\]
Then
\[
S_{1(k)}^{-1}=\Sigma_1^{-\frac{1}{2}}A_{(k)}\Sigma_1^{-\frac{1}{2}},\qquad
S_{2(k)}^{-1}=\Sigma_2^{-\frac{1}{2}}B_{(k)}\Sigma_2^{-\frac{1}{2}}.
\]
In order to calculate the value $\pr_{2|1}^{\D}$, we rewrite $\D$ as
\begin{eqnarray*}
\D&=&\frac{1}{H}\sum_{k=1}^H\left[\frac{1}{s_{0(k)}}(\bbz-\bar\bbx_{(k)})^TS_{1(k)}^{-1}(\bbz-\bar\bbx_{(k)})
+\log|S_{1(k)}|-pl_{1(k)}\right]\non
&-&\frac{1}{H}\sum_{k=1}^H\left[\frac{1}{m_{0(k)}}(\bbz-\bar\bby_{(k)})^TS_{2(k)}^{-1}(\bbz-\bar\bby_{(k)})
+\log|S_{2(k)}|-pl_{2(k)}\right]\non
&=&Q_1^{\D}+Q_2^{\D}+Q_3^{\D}+\sqrt{p}(T_3+T_4^{\D})+\log|\Sigma_1\Sigma_2^{-1}|+Q_4^{\D},
\end{eqnarray*}
where
{\small
\begin{eqnarray*}
&&Q_1^{\D}=(\bbz^0)^T\frac{1}{H}\sum_{k=1}^H\bigg(\frac{1}{s_{0(k)}}A_{(k)}-\frac{1}{m_{0(k)}}\Gamma^T B_{(k)}\Gamma\bigg)\bbz^0,\non
&&Q_2^{\D}=\frac{1}{H}\sum_{k=1}^H\left[-\frac{2}{s_{0(k)}}(\bar\bbx_{(k)}^0)^TA_{(k)}\bbz^0
+\frac{1}{s_{0(k)}}(\bar\bbx_{(k)}^0)^TA_{(k)}\bar\bbx_{(k)}^0
-\frac{2}{m_{0(k)}}(\bar\bby_{(k)}^0)^TB_{(k)}\Gamma\bbz^0+\frac{1}{m_{0(k)}}(\bar\bby_{(k)}^0)^TB_{(k)}\bar\bby_{(k)}^0\right],\non
&&Q_3^{\D}=\frac{1}{H}\sum_{k=1}^H\left[-\frac{2}{m_{0(k)}}(\Gamma\bbz^0-\bar\bby_{(k)}^0)^TB_{(k)}\Sigma_2^{-\frac{1}{2}}(\bbmu_1-\bbmu_2)\right],\non
&&T_3=-\frac{1}{\sqrt{p}}(\bbmu_1-\bbmu_2)^T\Sigma_2^{-1}(\bbmu_1-\bbmu_2),\quad T_4^{\D}=-\frac{1}{\sqrt{p}}\frac{1}{H}\sum_{k=1}^H\left[\frac{1}{m_{0(k)}}(\bbmu_1-\bbmu_2)^T\Sigma_2^{-\frac{1}{2}}B_{(k)}\Sigma_2^{-\frac{1}{2}}(\bbmu_1-\bbmu_2)\right]-T_3,\non
&&Q_4^{\D}=\frac{1}{H}\sum_{k=1}^H\left[(\log|S_{1(k)}^0|-pl_{1(k)})-(\log|S_{2(k)}^0|-pl_{2(k)})\right].
\end{eqnarray*}
}
From the derivation of $Q_1$ in Lemma \ref{lema2}, we know that
\[
\frac{1}{\sqrt{p}}Q_1^{D}-\frac{1}{\sqrt{p}}\tr C_{\D} \xrightarrow{D}N(0,\psi_{\D}^2),
\]
where $C_{\D}=\frac{1}{H}\sum\limits_{k=1}^H\bigg(\frac{1}{s_{0(k)}}A_{(k)}-\frac{1}{m_{0(k)}}\Gamma^T B_{(k)}\Gamma\bigg)=(c_{ij}^{\D})_{p\times p}$, $\psi_{\D}^2=(m_4-1)\xi_{\D}+2(\eta_{\D}-\xi_{\D})$ and
\begin{equation*}
\xi_{\D}=\lim_{p\rightarrow \infty}\frac{1}{p}\sum_{i=1}^p (c_{ii}^{\D})^2,\qquad \eta_{\D}=\lim_{p\rightarrow \infty}\frac{1}{p}\tr C_{\D}^2.
\end{equation*}
Moreover,
\begin{eqnarray*}
\frac{1}{\sqrt{p}}\tr C_{\D}&=&
\frac{1}{\sqrt{p}}\frac{1}{H}\sum\limits_{k=1}^H\bigg(\frac{1}{s_{0(k)}}\tr A_{(k)}-\frac{1}{m_{0(k)}}\tr \Gamma^T B_{(k)}\Gamma\bigg)=\frac{1}{H}\sum\limits_{k=1}^H\bigg(\sqrt{p}-\frac{1}{\sqrt{p}}\sum_{i=1}^p\gamma_i^T\gamma_i
+O_p(\frac{1}{\sqrt{p}})\bigg)\non
&=&\frac{1}{\sqrt{p}}\tr(\bbI_p-\Sigma_1\Sigma_2^{-1})+O_p(\frac{1}{\sqrt{p}})
\end{eqnarray*}
and
\[\xi_{\D}=
\lim_{p\rightarrow \infty}\frac{1}{p}\sum_{i=1}^p (c_{ii}^{\D})^2
=\lim_{p\rightarrow \infty}\frac{1}{p}\sum_{i=1}^p (1-\gamma_i^T\gamma_i)^2
=\lim_{p\rightarrow \infty}\left[1-\frac{2}{p}\tr(\Sigma_1\Sigma_2^{-1})+\frac{1}{p}\sum_{i=1}^p[(\Sigma_1^{\frac{1}{2}}\Sigma_2^{-1}\Sigma_1^{\frac{1}{2}})_{ii}]^2\right]
=\xi_0.
\]
Next we look at the value $\eta_{\D}$. To this end, we use the fact that  $s_{0(k)}=s_{0(1)}, m_{0(k)}=m_{0(1)}, k=1,\cdots, H$ and denote
$\frac{1}{p}\tr C_{\D}^2=A_{\D}-2B_{\D}+E_{\D}$,
where $A_{\D}=\frac{1}{s_{0(1)}^2}\frac{1}{p}\tr\left[\frac{1}{H}\sum\limits_{k=1}^H A_{(k)}\right]^2$,
\[
B_{\D}=\frac{1}{s_{0(1)}m_{0(1)}}\frac{1}{p}\tr\left[\frac{1}{H}\sum_{k=1}^H A_{(k)}\right]\left[\frac{1}{H}\sum_{k=1}^H \Gamma^T B_{(k)}\Gamma\right],\quad
E_{\D}=\frac{1}{m_{0(1)}^2}\frac{1}{p}\tr\left[\frac{1}{H}\sum_{k=1}^H \Gamma^T B_{(k)}\Gamma\right]^2.
\]
We study the three terms one by one below. First write
\begin{eqnarray*}
A_{\D}=\frac{1}{s_{0(1)}^2}\frac{1}{pH^2}\left[\sum\limits_{k=1}^H \tr A_{(k)}^2+\mathop{\sum^H\sum^H}_{k_1\neq k_2}\tr A_{(k_1)}A_{(k_2)}\right].
\end{eqnarray*}
By (\ref{yq10.6}), it can be seen that
$\frac{1}{p}\tr A_{(k)}^2=s'_{0(k)}+O_p(\frac{1}{p})$, where $s'_{0(k)}=\frac{1}{(1-c_{1(H)})^3}$. Then
\begin{equation}\label{1806.3}
\frac{1}{s_{0(1)}^2}\frac{1}{pH^2}\sum\limits_{k=1}^H \tr A_{(k)}^2=\frac{1}{H}\cdot\frac{1}{1-c_{1(H)}}+O_p(\frac{1}{p}).
\end{equation}
By similar methods in deriving $\frac{1}{p}\tr(A\Gamma^TB\Gamma)$ in the proof of Lemma \ref{lema2}, we can get that
\[
\frac{1}{p}\tr A_{(k_1)}A_{(k_2)}=s_{0(1)}^2+O_p(\frac{1}{p}),\quad
\frac{1}{p}\tr (A_{(k_1)}\Gamma^T B_{(k_2)}\Gamma)=s_{0(1)}m_{0(1)}\frac{1}{p}\tr(\Sigma_1\Sigma_2^{-1})+O_p(\frac{1}{p}).
\]
 Then
\begin{eqnarray}\label{1806.4}
&&\frac{1}{s_{0(1)}^2}\frac{1}{pH^2}\mathop{\sum^H\sum^H}_{k_1\neq k_2}\tr A_{(k_1)}A_{(k_2)}=\frac{H-1}{H}+O_p(\frac{1}{p}),\non
&&\frac{1}{s_{0(1)}m_{0(1)}}\frac{1}{pH^2}\sum_{k_1=1}^H\sum_{k_2=1}^H\tr (A_{(k_1)}\Gamma^T B_{(k_2)}\Gamma)=\frac{1}{p}\tr(\Sigma_1\Sigma_2^{-1})+O_p(\frac{1}{p}).
\end{eqnarray}
Equalities (\ref{1806.3}) and (\ref{1806.4}) indicate that \[A_{\D}=\frac{1}{H}\frac{1}{1-c_{1(H)}}+\frac{H-1}{H}+O_p(\frac{1}{p})=\frac{1}{H}\cdot\frac{c_{1(H)}}{1-c_{1(H)}}+1+O_p(\frac{1}{p}),
\quad B_{\D}=\frac{1}{p}\tr(\Sigma_1\Sigma_2^{-1})+O_p(\frac{1}{p}).
\]
Next look at the term $E_{\D}$. By the derivation of $\frac{1}{p}\tr(\Gamma^T B\Gamma)^2$ in the proof of Lemma \ref{lema2}, we have $\frac{1}{p}\tr(\Gamma^T B_{(k)}\Gamma)^2=m_{0(1)}^2\left[\frac{1}{p}\tr(\Sigma_1\Sigma_2^{-1})^2+\frac{c_{2(H)}}{1-c_{2(H)}}\Big(\frac{1}{p}\tr(\Sigma_1\Sigma_2^{-1})\Big)^2\right]
+O_p(\frac{1}{p})$. Moreover, when $k_1\neq k_2$, $\frac{1}{p}\tr(\Gamma^T B_{(k_1)}\Gamma)(\Gamma^T B_{(k_2)}\Gamma)
=\frac{m_{0(1)}^2}{p}\tr(\Sigma_1\Sigma_2^{-1})^2+O_p(\frac{1}{p})$. Therefore,
\begin{eqnarray*}
E_{\D}&=&\frac{1}{m_{0(1)}^2}\frac{1}{pH^2}\left[\sum\limits_{k=1}^H \tr (\Gamma^T B_{(k)}\Gamma)^2
+\mathop{\sum^H\sum^H}_{k_1\neq k_2}\tr(\Gamma^T B_{(k_1)}\Gamma)(\Gamma^T B_{(k_2)}\Gamma)\right]\non
&=&\frac{1}{H}\left[\frac{1}{p}\tr(\Sigma_1\Sigma_2^{-1})^2+\frac{c_{2(H)}}{1-c_{2(H)}}\Big(\frac{1}{p}\tr(\Sigma_1\Sigma_2^{-1})\Big)^2\right]
+\frac{H-1}{H}\frac{1}{p}\tr(\Sigma_1\Sigma_2^{-1})^2+O_p(\frac{1}{p})\non
&=&\frac{1}{H}\frac{c_{2(H)}}{1-c_{2(H)}}\Big(\frac{1}{p}\tr(\Sigma_1\Sigma_2^{-1})\Big)^2+\frac{1}{p}\tr(\Sigma_1\Sigma_2^{-1})^2+O_p(\frac{1}{p}).
\end{eqnarray*}
Combing the results for the three terms $A_{\D}$, $B_{\D}$ and $E_{\D}$, we have
\begin{eqnarray*}
\eta_{\D}&=&\lim_{p\rightarrow \infty}(A_{\D}-2B_{\D}+E_{\D})\non
&=&
\lim_{p\rightarrow \infty}\left[1-\frac{2}{p}\tr(\Sigma_1\Sigma_2^{-1})+\frac{1}{p}\tr(\Sigma_1\Sigma_2^{-1})^2\right]+
\frac{1}{H}\lim_{p\rightarrow \infty}\left[\frac{c_{1(H)}}{1-c_{1(H)}}+\frac{c_{2(H)}}{1-c_{2(H)}}\Big(\frac{1}{p}\tr(\Sigma_1\Sigma_2^{-1})\Big)^2\right]\non
&=&\eta_0+\frac{1}{H}\left[\frac{c_{1(H)}}{1-c_{1(H)}}+\frac{c_{2(H)}}{1-c_{2(H)}}M_1^2\right].
\end{eqnarray*}
Therefore,
$\psi_{\D}^2=(m_4-1)\xi_{\D}+2(\eta_{\D}-\xi_{\D})=\psi_0^2+\frac{2}{H}\left[\frac{c_{1(H)}}{1-c_{1(H)}}+\frac{c_{2(H)}}{1-c_{2(H)}}M_1^2\right]$ and we get
$$\pr_{2|1}^{\D}=\pr\{\D>0\}\xrightarrow{i.p} 1-\Phi\left(\frac{T}{\psi_{\D}}\right).$$
The probability $\pr_{1|2}^{\D}$ can be derived in a similar way and the proof of Theorem \ref{them4} is done.

\section{Comparison with Sample QDA}\label{compare2}
As a byproduct of Theorem \ref{themnew}'s proof, we can also theoretically analyze the asymptotic misclassification rate of the sample QDA under moderate dimension.

\begin{figure}[!htp]
    \centering
    \includegraphics[height=1.8in, width=2.5in]{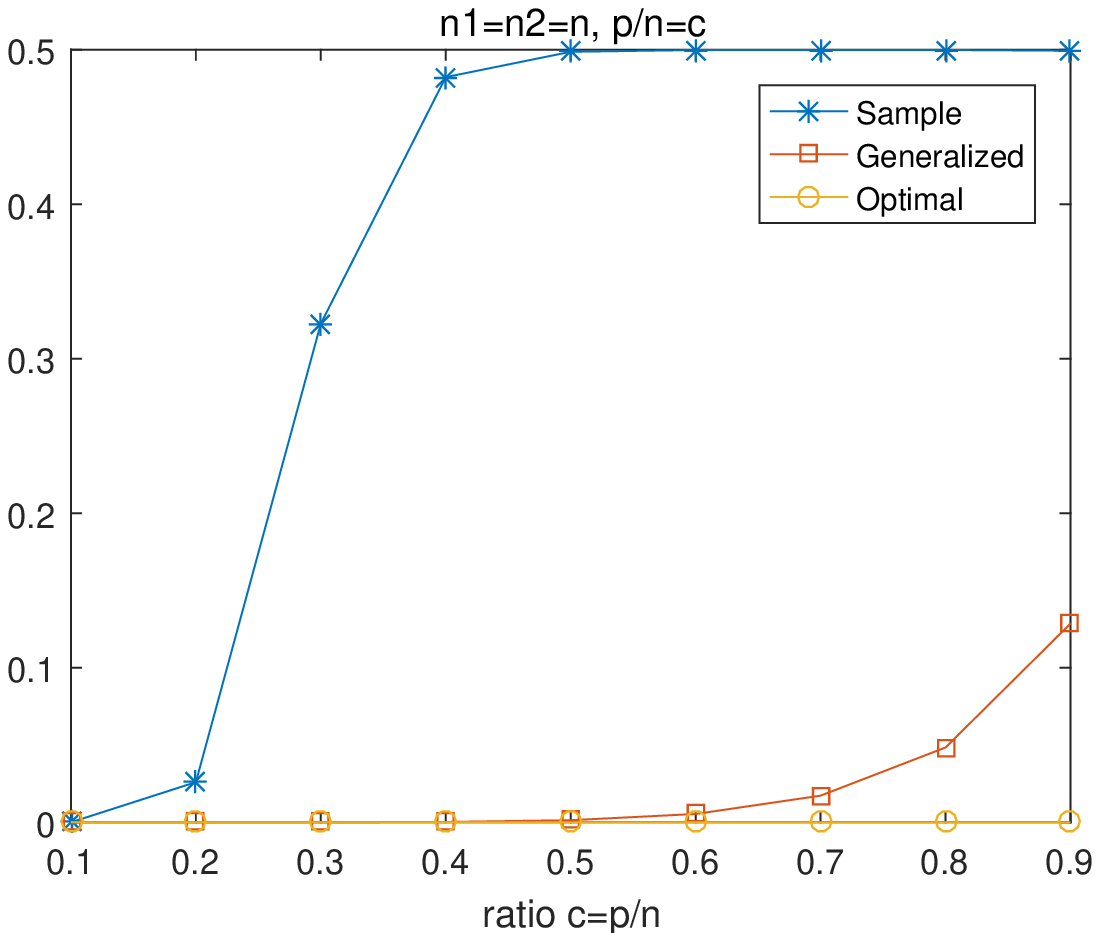}\quad
    \includegraphics[height=1.8in, width=2.5in]{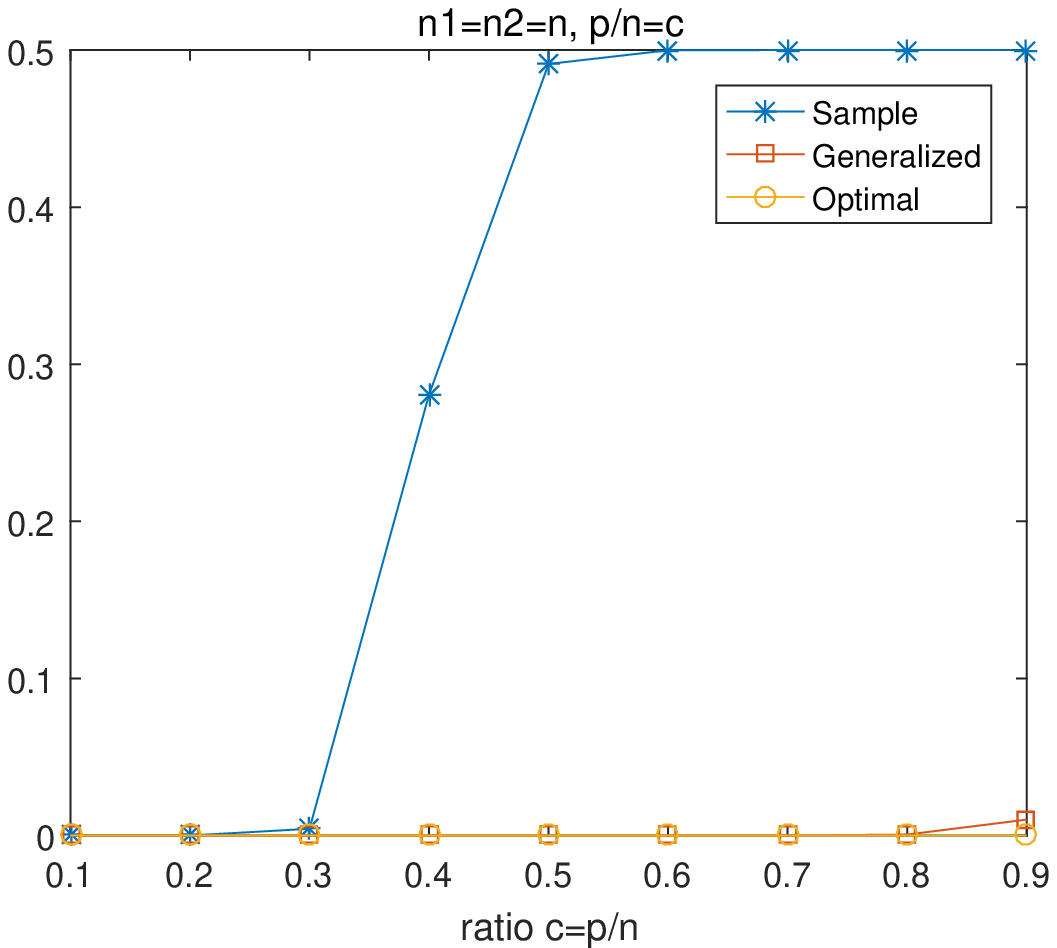}
    \caption{\small{\it Plots of the limits of $R^S$, $R^G$ and $R^O$. $p=1000$,  $\bbmu_1=\bbmu_2$, $n_1=n_2$ and $m_4=3$.  Left: $\Sigma_2=0.5*\Sigma_1$; Right:  $\Sigma_2=3*\Sigma_1$. }}
    \label{compar}
  \end{figure}

\begin{prop}\label{prop2}
Under Conditions \ref{cond1}-\ref{cond3}, the misclassification rate of the sample QDA (\ref{yq4})
\[
R^S=\frac{1}{2}\left[\pr_{2|1}^S+\pr_{1|2}^S\right]\xrightarrow{i.p}
1-\frac{1}{2}\left[\Phi\left(\frac{T_S}{\psi_S}\right)
+\Phi\left(\frac{\widetilde{T}_S}{\widetilde{\psi}_S}\right)\right],
\]
where
{\small
\[
T_S=\lim_{p\rightarrow\infty}\Big\{-\frac{1}{\sqrt{p}}[s_{0n}\tr\bbI_p-m_{0n}\tr\Sigma_1\Sigma_2^{-1}]
-\frac{1}{\sqrt{p}}\log|\Sigma_1\Sigma_2^{-1}|+\sqrt{p}(l_{2n}-l_{1n})
+\frac{m_{0n}}{\sqrt{p}}(\bbmu_1-\bbmu_2)^T\Sigma_2^{-1}(\bbmu_1-\bbmu_2)\Big\},
\]
\[
\widetilde{T}_S=\lim_{p\rightarrow\infty}\Big\{-\frac{1}{\sqrt{p}}[m_{0n}\tr\bbI_p-s_{0n}\tr\Sigma_2\Sigma_1^{-1}]
-\frac{1}{\sqrt{p}}\log|\Sigma_2\Sigma_1^{-1}|+\sqrt{p}(l_{1n}-l_{2n})
+\frac{s_{0n}}{\sqrt{p}}(\bbmu_1-\bbmu_2)^T\Sigma_1^{-1}(\bbmu_1-\bbmu_2)\Big\}.
\]
}
The parameters $\psi_S$ and $\widetilde{\psi}_S$ are positive constants given by
\[
\psi_S^2=(m_4-3)(s_0^2-2s_0m_0M_1+m_0^2M_2)+2\left[\frac{1}{(1-c_1)^3}-2s_0m_0M_1+m_0^2\Big(M_5+\frac{c_2}{1-c_2}M_1^2\Big)\right],
\]
\[
\widetilde{\psi}_S^2=(m_4-3)(m_0^2-2s_0m_0M_3+s_0^2M_4)+2\left[\frac{1}{(1-c_2)^3}-2s_0m_0M_3+s_0^2\Big(M_5+\frac{c_1}{1-c_1}M_3^2\Big)\right].
\]
\end{prop}
The limit of $R^S$ above involves too many uncertain parameters and we are unable to get an analogue of Corollary \ref{coronew} and Corollary \ref{coro3} under the three cases. Instead,  in Figure \ref{compar}, we plot the limit of $R^S$ under some specific settings -- ($\bbmu_1=\bbmu_2$ and $\Sigma_2=\kappa*\Sigma_1$) -- satisfying case (i).  One may observe that different from the conclusion $R^G\rightarrow 0$ and $R^O\rightarrow 0$ in Corollary \ref{coronew} and Corollary \ref{coro3}, $R^S$ could be significantly larger than zero and even behaves like random guessing for large ratio $c$.

\subsection{Proof of Proposition \ref{prop2}}
Below we only give the derivation of $\pr_{2|1}^S$, the other one $\pr_{1|2}^S$ can be calculated in a similar way. Write
\begin{eqnarray}\label{1807.1}
\pr_{2|1}^S
&=&\pr\left\{D_1(\bbz)+\log |S_1|>D_2(\bbz)+\log |S_2|\right\}\non
&=&\pr\bigg\{\frac{1}{\sqrt{p}}\Big(Q_1^S+ Q_2^S+Q_3^S+Q_4^S\Big)+\frac{1}{\sqrt{p}}\log|\Sigma_1\Sigma_2^{-1}|+\sqrt{p}(l_{1n}-l_{2n})+T_5^S>0\bigg\},\non
&&
\end{eqnarray}
where
\begin{eqnarray*}
&&Q_1^S=(\bbz^0)^T\big(A-\Gamma^T B\Gamma\big)\bbz^0,\quad Q_2^S=-2(\bar\bbx^0)^TA\bbz^0+(\bar\bbx^0)^TA\bar\bbx^0
-2(\bar\bby^0)^TB\Gamma\bbz^0+(\bar\bby^0)^TB\bar\bby^0,\non
&&Q_3^S=-2(\Gamma\bbz^0-\bar\bby^0)^TB\Sigma_2^{-\frac{1}{2}}(\bbmu_1-\bbmu_2),\quad
Q_4^S=-(\bbmu_1-\bbmu_2)^T\Sigma_2^{-\frac{1}{2}}B\Sigma_2^{-\frac{1}{2}}(\bbmu_1-\bbmu_2),\non
&&T_5^S=\frac{1}{\sqrt{p}}\Big(\log|S_1^0|-pl_{1n}\Big)-\frac{1}{\sqrt{p}}\Big(\log|S_2^0|-pl_{2n}\Big).
\end{eqnarray*}
Let $T_3^S=-\frac{m_{0n}}{\sqrt{p}}(\bbmu_1-\bbmu_2)^T\Sigma_2^{-1}(\bbmu_1-\bbmu_2)$ and $T_4^S=\frac{1}{\sqrt{p}}Q_4^S-T_3^S$.
Applying similar arguments as in Section \ref{proflem}, we know that
\begin{equation}\label{1807.2}
\frac{1}{\sqrt{p}}Q_2^S=O_p(\frac{1}{\sqrt{p}}),\quad
\frac{1}{\sqrt{p}}Q_{3}^S=O_p(\frac{\sqrt{|T_3^S|}}{p^{1/4}}),\quad
T_5^S=O_p(\frac{1}{\sqrt{p}}),\quad
T_4^S=O_p(\frac{|T_3^S|}{\sqrt{p}}).
\end{equation}
Denote $C_S=A-\Gamma^T B\Gamma=(c_{ij}^S)_{p\times p}$. Then
\begin{equation}\label{1807.3}
\frac{1}{\sqrt{p}}(Q_1^S-\tr C_S)\xrightarrow{D} N(0,\psi_S^2),
\end{equation}
where $\psi_S^2=(m_4-1)\xi_S+2(\eta_S-\xi_S)$,
$\xi_S=\lim\limits_{p\rightarrow \infty}\frac{1}{p}\sum\limits_{i=1}^p (c_{ii}^S)^2$ and
$\eta_S=\lim\limits_{p\rightarrow \infty}\frac{1}{p}\tr C_S^2$.
As in deriving the terms in Section \ref{proflem}, we can get that
\begin{eqnarray*}
\frac{1}{\sqrt{p}}\tr C_S&=&\frac{1}{\sqrt{p}}\tr A-\frac{1}{\sqrt{p}}\tr\Gamma^T B\Gamma=\sqrt{p}s_{0n}-\frac{m_{0n}}{\sqrt{p}}\sum_{i=1}^p\gamma_i^T\gamma_i+O_p(\frac{1}{\sqrt{p}})\non
&=&\sqrt{p}s_{0n}-\frac{m_{0n}}{\sqrt{p}}\tr\Gamma^T\Gamma+O_p(\frac{1}{\sqrt{p}})
=\frac{1}{\sqrt{p}}[s_{0n}\tr\bbI_p-m_{0n}\tr\Sigma_1\Sigma_2^{-1}]+O_p(\frac{1}{\sqrt{p}}),
\end{eqnarray*}
$$
\xi_S=\lim_{p\rightarrow \infty}\frac{1}{p}\sum\limits_{i=1}^p (c_{ii}^S)^2
=\lim_{p\rightarrow\infty}\frac{1}{p}\sum\limits_{i=1}^p(s_{0n}-m_{0n}\gamma_i^T\gamma_i)^2
=s_0^2-2s_0m_0M_1+m_0^2M_2,
$$
and
\begin{eqnarray*}
\eta_S&=&\lim_{p\rightarrow \infty}\frac{1}{p}\tr (A-\Gamma^T B\Gamma)^2\non
&=&s'_0+\lim_{p\rightarrow\infty}\left[-\frac{2s_{0n}m_{0n}}{p}\tr(\Sigma_1\Sigma_2^{-1})+
m_{0n}^2\Big[\frac{1}{p}\tr(\Sigma_1\Sigma_2^{-1})^2+\frac{c_{2n}}{1-c_{2n}}\Big(\frac{1}{p}\tr(\Sigma_1\Sigma_2^{-1})\Big)^2\Big]\right]\non
&=&s'_0-2s_0m_0M_1+m_0^2\Big(M_5+\frac{c_2}{1-c_2}M_1^2\Big).
\end{eqnarray*}
Moreover, $\eta_S>s'_0-s_0^2=\frac{c_1}{(1-c_1)^3}>0$. Therefore
\[
\psi_S^2=(m_4-3)(s_0^2-2s_0m_0M_1+m_0^2M_2)+2\left[s'_0-2s_0m_0M_1+m_0^2\Big(M_5+\frac{c_2}{1-c_2}M_1^2\Big)\right]>0.
\]
Then the probability $\pr_{2|1}^S$ can be derived by
\eqref{1807.1}, \eqref{1807.2} and \eqref{1807.3}:
$\pr_{2|1}^S\xrightarrow{i.p}1-\Phi\left(\frac{T_S}{\psi_S}\right)$, where
$
T_S=\lim\limits_{p\rightarrow\infty}\Big\{-\frac{1}{\sqrt{p}}[s_{0n}\tr\bbI_p-m_{0n}\tr\Sigma_1\Sigma_2^{-1}]
-\frac{1}{\sqrt{p}}\log|\Sigma_1\Sigma_2^{-1}|+\sqrt{p}(l_{2n}-l_{1n})
+\frac{m_{0n}}{\sqrt{p}}(\bbmu_1-\bbmu_2)^T\Sigma_2^{-1}(\bbmu_1-\bbmu_2)\Big\}$.

\end{document}